\pdfoutput=1
\documentclass[10pt]{article}

\usepackage{amsmath, amsfonts, amsbsy}
\usepackage{multirow, makecell, booktabs, array}
\usepackage{tikz-cd}
\usepackage{xcolor}
\usepackage{subcaption}
\usepackage{graphicx,url}
\usepackage{dblfloatfix}
\usepackage{pifont, marvosym, nicefrac, microtype}
\usepackage{authblk}
\usepackage{listings}
\usepackage[margin=1in]{geometry}
\usepackage[hidelinks]{hyperref}
\usepackage{xurl}

\usepackage{xr-hyper}
\newcommand{\xref}[2]{#1}

\newcommand{\R}{\mathbb{R}}
\newcommand{\PDEgym}{\textsc{PDEgym}}
\newcommand{\greencheck}{{\color{green}\ding{51}}}
\newcommand{\redcross}{{\color{red}\ding{55}}}
\newcommand\verB[1]{#1}
\newcommand\verC[1]{#1}
\newcommand{\keywords}[1]{\par\medskip\noindent\textbf{Keywords: }#1\par\vspace{\baselineskip}}
\newcommand{\tablefont}{\small}
\begin{document}
\title{PDEformer-2: a versatile foundation model for two-dimensional partial differential equations}
\author[1]{Zhanhong Ye}
\author[2]{Zining Liu}
\author[2]{Bingyang Wu}
\author[2]{Hongjie Jiang}
\author[1]{Leheng Chen}
\author[3]{Minyan Zhang}
\author[4]{Xiang Huang}
\author[5]{Qinghe Meng}
\author[5]{Jingyuan Zou}
\author[3]{Hongsheng Liu}
\author[6,7,8]{Bin Dong\thanks{Corresponding author: dongbin@math.pku.edu.cn}}
\affil[1]{Beijing International Center for Mathematical Research, Peking University, Beijing 100871, China}
\affil[2]{School of Mathematical Sciences, Peking University, Beijing 100871, China}
\affil[3]{Central Software Institute, Huawei Technologies Co., Ltd., Shanghai 201206, China}
\affil[4]{Central Software Institute, Huawei Technologies Co., Ltd., Hangzhou 310053, China}
\affil[5]{Computing Product Line, Huawei Technologies Co., Ltd., Hangzhou 310051, China}
\affil[6]{Beijing International Center for Mathematical Research and the New Cornerstone Science Laboratory, Peking University, Beijing 100871, China}
\affil[7]{Center for Machine Learning Research, Peking University, Beijing 100871, China}
\affil[8]{Zhongguancun Academy, Beijing 100094, China}
\date{}
\twocolumn[
\begin{@twocolumnfalse}
\maketitle
\begin{abstract}
	Partial differential equations (PDEs) are fundamental to modeling physical systems.
	Traditional PDE solvers are often computationally expensive.
	Neural operators require independent training for different PDEs, and tend to be data-greedy when adapted to new scenarios.
	To this end,
	we introduce PDEformer-2, a versatile foundation model for solving diverse two-dimensional PDEs.
	The model receives the PDE form as input through a flexible computational graph representation and
	outputs solutions \verB{that can be queried} at any spatio-temporal point.
	Pretrained on a 40 TB dataset spanning varied equation types, domains, boundary conditions,
	variable counts, and time dependence, PDEformer-2 enables instant zero-shot predictions for familiar PDEs
	and adapts faster to unseen ones than specialized baselines,
	achieving lower errors with limited samples.
	It also facilitates efficient inverse problem solving, successfully recovering unknown PDE parameters,
	including scalar coefficients and fields, in experiments.
\end{abstract}
\end{@twocolumnfalse}
]

\section{Introduction}\label{sec:intro}
Partial differential equations (PDEs) underpin countless scientific and engineering disciplines, from aircraft design and climate modeling to medical imaging and energy systems.
PDE solvers are thus critical for related design, control, and inverse problems.
Practical applications demand solvers with the following attributes:
(1) Speed: Solutions should be generated rapidly within limited resources,
especially in real-time or multi-query settings.
An offline preparation stage is often acceptable.
(2) Adequate accuracy: Prediction errors should lie within a reasonable range to support downstream use.
Meanwhile, small errors can often be tolerated when rapid prototyping is desired or when physical mechanisms are themselves imperfectly known.
(3) Differentiability: Derivatives of solutions with respect to parameters should be accessible
to support gradient-based optimization in design and inverse problems.
(4) Versatility: Solvers should handle diverse PDE types and coupled systems
without case-specific adjustments, avoiding workflow fragmentation and simplifying exploration of new scenarios.

Traditional numerical methods (e.g., finite difference, finite element)
have played a pivotal role over decades of development.
Built upon rigorous mathematical foundations, these solvers are versatile, reliable, and highly precise.
However, they are computationally expensive for applications where efficiency outweighs accuracy.
Although adjoint methods enable gradient computation, they introduce complexity and overhead.
Their adaptability to new problems often requires expert tuning and discretization adjustments.

Neural operators, such as Fourier Neural Operator (FNO)~\cite{FNO} and DeepONet~\cite{DeepONet}, learn solution mappings from data to provide fast and differentiable approximations.
Their reduced accuracy and reliability limit their use in safety-critical settings, but
can be acceptable for applications tolerant of small errors.
While offline data preparation and training can be amortized over repeated online evaluations,
this costly procedure has to be repeated for each new PDE because neural operators are equation-specific solvers.
Achieving desired accuracy also demands large datasets without cross-equation knowledge, which may impede adoption for data-scarce scenarios.

The specialized nature of existing neural operators motivates the development of PDE foundation models---unified architectures pretrained on diverse PDE families to achieve broad applicability~\cite{DD-FEM}.
Such models would combine the speed and differentiability of neural operators with greater generality.
Ideally, they would also offer:
(1) Discretization flexibility: Handling arbitrary discretization to support various geometries,
and generating mesh-free solutions.
(2) Zero-shot capability: Direct prediction for familiar physical settings without additional solver data.
(3) Efficient adaptation: Rapid finetuning for unseen equations with minimal data and computational cost.

Many efforts have been made to develop such a general-purpose foundation model for PDEs.
Some approaches include the symbolic form of the PDE as explicit model input,
using textual representations~\cite{PITT,ICON-LM,Unisolver,UPS,multimodal,explainfive} or
Polish notation of expression trees~\cite{PROSE,PROSE-PDE,LeMON-PROSE,PROSE-FD}.
These methods may face challenges in capturing the intricate interplay between the PDE form and numerical information (coefficient fields, boundary values, etc.),
and their generalization across diverse PDEs remains unverified.
Other strategies implicitly inform the model about the PDE,
using parameter-solution pairs~\cite{ICON,ICON-flux,VICON,GenICON}
or solution-history snapshots~\cite{MPP,DPOT,OmniArch,BCAT,PDE-Transformer,PhysiX,GPhyT,FMT},
but they could be inadequate for highly varied PDE classes
and require a traditional solver to prepare additional solution examples.
A third category ignores PDE information in pretraining, necessitating finetuning for any specific system
and precluding zero-shot prediction~\cite{towards-foundation-models,CoDA-NO,Poseidon,MOFS,MORPH}.

\begin{figure*}[htpb]
	\centering
	\includegraphics[width=400pt]{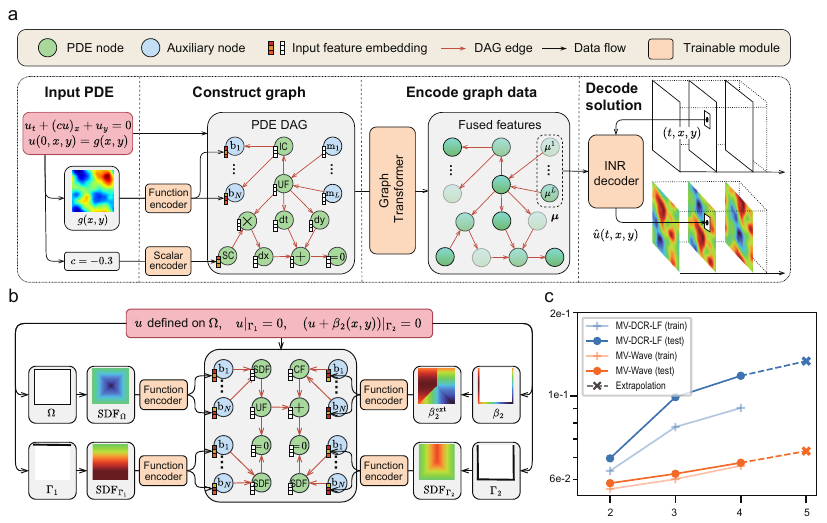}
	\caption{PDEformer-2 architecture and illustrative results.
		(a) Architecture of PDEformer-2, illustrated for an equation on $\Omega=[0,1]^2$ with periodic boundaries.
		(b) Computational graph example for Dirichlet boundary conditions.
		Function $\beta_2(\boldsymbol{r})$, defined on boundary $\Gamma_2$, is extended to the full square domain before encoding.
		(c) PDEformer-2-base accuracy on MV-DCR with linear fluxes (LF) and MV-Wave equations.
		The model generalizes zero-shot to systems with five variables, exceeding the pretraining bound of four.
	}%
	\label{fig:PDEformerV2Arch}
\end{figure*}

In this work, we introduce PDEformer-2, a foundation model for solving two-dimensional PDEs.
Building on our prior PDEformer-1~\cite{PDEformer1} for 1D systems,
we extend the computational graph representation to handle various domains and boundary conditions in 2D,
and pretrain the model on a much more diverse dataset.
Notable features of PDEformer-2 include:
\begin{itemize}
	\item Architectural flexibility: The model flexibly encodes the symbolic formulation
		of the PDE with a computational graph (Fig.~\ref{fig:PDEformerV2Arch}a),
		laying a good foundation for generality.
	\item Versatility: We generate a large pretraining dataset (40 TB),
		allowing PDEformer-2 to solve diverse PDEs with varying equation types, domains, boundary conditions,
		numbers of equations and variables, and time dependence.
	\item \verB{Mesh-independent solutions: Despite the fixed-grid encoding of function inputs, the output}
		solutions are given by an implicit neural representation (INR), queryable at arbitrary coordinates.
	\item Unified space-time prediction: The model predicts full solutions directly,
		which avoids time-consuming and error-accumulating rollouts and simplifies gradient computation.
	\item Zero-shot inference: Explicit incorporation of the PDE form enables zero-shot inference
		for systems covered in pretraining.
	\item Efficient finetuning: For PDEs unseen during pretraining, the model achieves better accuracy
		than specialized baselines with limited data and time.
	\item Downstream task support: Fast and differentiable inference facilitates use in inverse problems
		to recover unknown PDE parameters, including scalar coefficients and fields.
\end{itemize}
The paper is structured as follows.
We first present experimental results on pretraining, adaptation to new PDEs, and inverse problems.
We then conclude with a summary of findings and future directions, and outline our method.
Additional implementation details are provided in the Supplementary Information.

\section{Results}\label{sec:results}
We begin by detailing the pretraining dataset for PDEformer-2 and summarizing its performance on the test set.
The pretrained model is then evaluated on specific PDEs, either through zero-shot inference or after finetuning, and is further applied to inverse problems.

\subsection{Pretraining results}\label{sec:results_pretrain}
\paragraph{Datasets}
\begin{table}[tbp]
	\centering
\scriptsize
	\caption{Summary of pretraining dataset composition and PDEformer-2's relative errors on the test set, including two model variants.
	The dataset comprises 8 generic PDE types, each yielding diverse PDE formulations through randomized coefficient selection.}
	\label{tab:pretrain_data_overview}
	\begin{tabular*}{\columnwidth}{@{\extracolsep\fill}ccccc@{\extracolsep\fill}}
		\toprule
		\makecell{\textbf{PDE} \\\textbf{type}} & \makecell{\textbf{Training}\\\textbf{samples}} & \makecell{\textbf{Data size} \\\textbf{(TB)}} & \makecell{\textbf{nRMSE} \\\textbf{(base)\%}} & \makecell{\textbf{nRMSE} \\\textbf{(fast)\%}} \\
		\midrule
		DCR & 745k & 4.85 & 6.9 & 8.4 \\
		Wave & 605k & 4.34 & 6.5 & 8.6 \\
		MV-DCR & 660k & 12.24 & 12.6 & 13.5 \\
		DC-DCR & 400k & 7.49 & 17.0 & 19.7 \\
		MV-Wave & 320k & 6.20 & 5.3 & 6.5 \\
		DC-Wave & 200k & 2.62 & 9.5 & 12.2 \\
		G-SWE & 100k & 1.93 & 13.0 & 14.3 \\
		Elasticity & 255k & 0.38 & 3.8 & 3.5 \\
		\addlinespace[0.4em]
		Total & 3.29M & 40 & 8.7 & 10.3 \\
		\bottomrule
	\end{tabular*}
\end{table}
To pretrain a general-purpose PDE model, we construct a large and diverse dataset encompassing multiple PDE families.
Each family is defined by a generic equation form containing numerous terms.
Randomly setting coefficients to zero or one yields distinct specific PDE instances.
We include 8 generic forms: diffusion-convection-reaction (DCR), wave, multi-variable DCR (MV-DCR), divergence-constrained DCR (DC-DCR), MV-Wave, DC-Wave, generalized shallow-water equation (G-SWE), and steady-state elasticity.
The dataset also varies in domain geometry, boundary conditions (BCs),
number of solution variables (also called ``solution components''), and time dependence.
Data were generated using spectral (Dedalus~V3~\cite{Dedalus}) and finite-element (FEniCSx~\cite{FEniCSx}) solvers.
Table~\ref{tab:pretrain_data_overview} provides an overview of the pretraining dataset,
and more details are given in Supplementary Note~\xref{S3.1}{app:dataset_pretrain}.
Note that this dataset, in which each sample corresponds to a different PDE, is unsuitable for standard neural operators designed for fixed equations.

\paragraph{Test set predictions}
\begin{figure*}[pht]
	\centering
	\includegraphics[width=\textwidth]{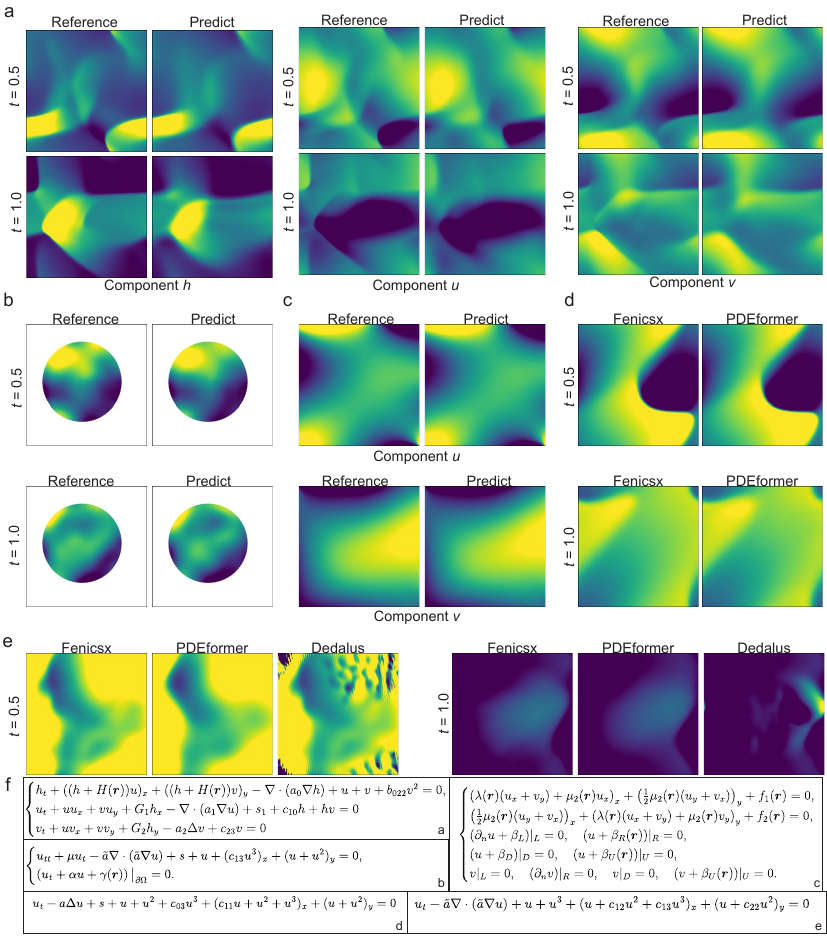}
	\caption{Test-set predictions of PDEformer-2-base.
		(a) Generalized shallow-water equation (with \verB{periodic} boundaries).
		(b) Wave equation in a disk.
		(c) Steady-state elasticity equation.
		(d) Comparison between the solutions produced by PDEformer-2 and FEniCSx for a DCR equation
		in a selected case where Dedalus fails to solve.
		(e) Comparison between the solutions produced by PDEformer-2, FEniCSx, and Dedalus, in a case where Dedalus exhibits apparent artifacts.
		(f) PDE formulations for the cases shown above. Note that these PDEs are randomly generated
		and may lack direct physical interpretation.
	}%
	\label{figs:results_forward}
\end{figure*}
We pretrained two model variants, PDEformer-2-base and PDEformer-2-fast.
The base model's hyperparameters mainly follow the XL configuration of PDEformer-1~\cite{PDEformer1},
and full spatio-temporal inference is relatively slow.
PDEformer-2-fast trades a slight accuracy loss for improved efficiency (Table~\ref{tab:inference_time}).
Figure~\ref{figs:results_forward}a-c shows test set predictions of PDEformer-2-base.
These cases differ in PDE types, domain shapes, boundary conditions, variable counts, and time dependence, demonstrating the model's flexibility.

\paragraph{Extrapolating to PDEs with more variables}
A key expectation for foundation PDE models is generalization beyond the pretraining distribution.
Here we evaluate PDEformer-2's ability to extrapolate to PDEs involving more variables than encountered during pretraining,
which corresponds to additional \verb|UF| (unknown field) nodes in the computational graphs.
While pretraining data includes at most four interacting variables, the model achieves reasonable zero-shot accuracy on systems with five variables, as shown in Fig.~\ref{fig:PDEformerV2Arch}c,
indicating its preliminary \verB{compositional generalization capability across PDEs} through computational graph representations.
Additional forms of extrapolation, such as unseen PDE formulations or coefficient field distributions,
are examined in the next subsection on adaptation to specific PDEs, with both zero-shot and few-shot settings.

\paragraph{Robustness to non-physical artifacts}
A small fraction of pretraining data contains non-physical solutions (e.g., spurious oscillations)
due to occasional instability of the Dedalus solver.
To assess model robustness, we regenerated reference solutions for selected DCR equations using FEniCSx, and compared them against Dedalus solutions and PDEformer-2 predictions.
As shown in Fig.~\ref{figs:results_forward}e, PDEformer-2-base prediction aligns with the physically consistent FEniCSx solution when Dedalus produces oscillatory artifacts.
This suggests the model effectively disregards erroneous pretraining samples and generalizes from physically plausible data.
PDEformer-2 outperforms its teacher solver in this specific sense.

\paragraph{Performance on failure cases of teacher solver}
In addition to producing non-physical solutions, Dedalus may also fail for certain PDEs,
resulting in \texttt{NaN} outputs.
These invalid solutions are excluded from the pretraining dataset.
We test PDEformer-2-base on such failure cases using FEniCSx-generated references,
and find that the model can occasionally produce reasonable predictions, as shown in Fig.~\ref{figs:results_forward}d.
However, the advantage over the teacher solver is limited in this setting because accuracy is
unsatisfactory in many cases.
A representative example (Supplementary Fig.~\xref{S12}{fig:dedalus_failed_sample_avg}) shows that although the model captures broad physical trends, its imperfect handling of shock evolution leads to observable errors.

\subsection{Adaptation to specific PDEs}\label{sec:results_forward}
This section evaluates PDEformer-2's performance on forward problems for solving fixed PDEs.
We compare against several specialized neural operator baselines that are trained from scratch for each PDE,
including DeepONet~\cite{DeepONet} with different branch networks, U-Net~\cite{U-Net}, FNO~\cite{FNO}, and Geo-FNO~\cite{Geo-FNO}.
All models predict full space-time solutions, avoiding autoregressive error accumulation.
Time is treated either as extended channels (2D) or an extra spatial dimension (3D).
We also include two existing PDE foundation models, MPP~\cite{MPP} (rollout-based) and Poseidon~\cite{Poseidon} (per-step prediction).
Further model and dataset details are provided in Supplementary Notes~\xref{S1.3}{app:model_baseline}
and~\xref{S3.2}{app:dataset_forward}.

\paragraph{Zero-shot and few-shot performance}
\begin{figure*}[!t]
	\centering
	\includegraphics[width=\textwidth]{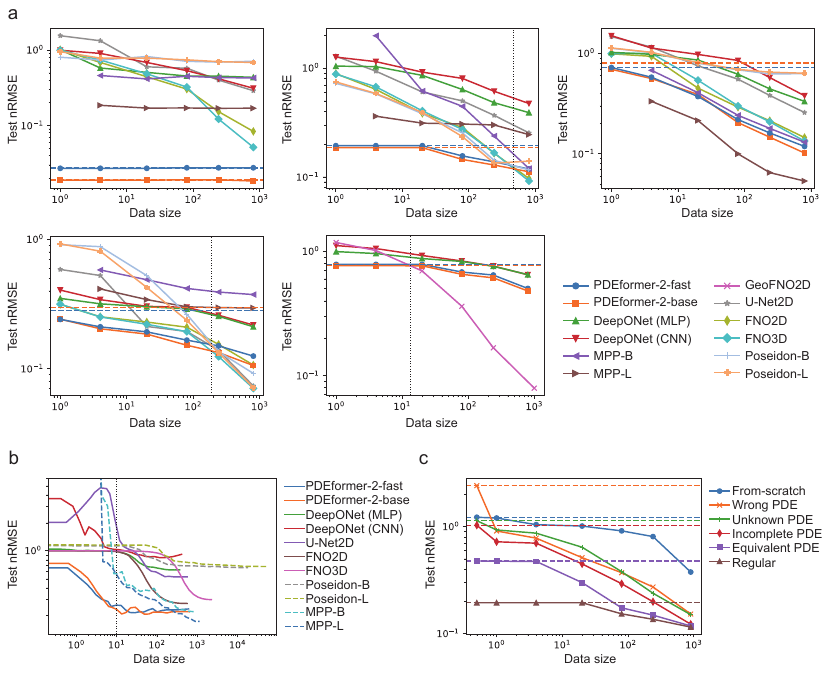}
	\caption{Adaptation performance of PDEformer-2 across data availability,
		adaptation speed, and PDE input quality.
		(a) Test nRMSE versus training data size.
		Horizontal dashed lines: zero-shot performance of PDEformer-2.
		Vertical dotted lines: data size beyond which PDEformer-2 is outperformed.
		Most baselines require Cartesian grids, and are incompatible with irregular discretization in Wave-C-Sines unless we interpolate the solution.
		(b) Finetuning speed for 1,000 epochs using 80 INS-Pipe samples.
		\verB{Poseidon-B/L and MPP-B/L (dashed lines) are executed on PyTorch+V100, and their time axis is platform-rescaled (Supplementary Note~\xref{S4.3}{app:prediction_speed}).}
		(c) Ablation of PDEformer-2-fast on INS-Tracer with varying PDE input quality or from-scratch training.
		Dashed lines and the left endpoints indicate zero-shot performance.
	}
	\label{figs:ablation}
\end{figure*}
Intuitively, PDEformer-2 will achieve high zero-shot accuracy on in-distribution PDEs.
For out-of-distribution cases, finetuning improves performance,
and can outperform specialized models by leveraging prior knowledge from pretraining, especially when data is limited.
Such transfer efficiency may also depend on similarity to pretraining data.

Experimental results across three custom and two public datasets confirm this intuition.
In Fig.~\ref{figs:ablation}a, we report test normalized root mean square error (nRMSE) versus training data size after finetuning foundation models or training specialized baselines from scratch.
On the in-distribution Sine-Gordon dataset, PDEformer-2 outperforms all baselines even without finetuning.
For INS-Tracer (slight shift, new combination of known mechanisms), zero-shot accuracy remains reasonable,
and finetuning helps gradually.
On INS-Pipe (new combination of physics and boundaries), zero-shot performance is poor,
but finetuning efficiently restores superiority 
\verB{over all baselines except MPP-L.
The advantage of MPP-L possibly originates from its much larger capacity and fluid-centered pretraining.}
For Wave-Gauss~\cite{Poseidon} (different initial conditions and boundary treatments), PDEformer-2
shows a marginal low-data advantage but loses its lead with more samples.
On the disk-domain Wave-C-Sines~\cite{RIGNO} dataset, PDEformer-2 handles irregular
\verB{solution queries} via its INR decoder but underperforms relative to Geo-FNO,
possibly due to the lack of high-frequency initial conditions and insufficient non-square domains in pretraining data.

\paragraph{Adaptation speed}
Beyond data efficiency, prior knowledge in PDEformer-2 reduces computational demands during adaptation.
As shown in Fig.~\ref{figs:ablation}b, our models converge within ten seconds, over an order of magnitude faster than baselines.
\verB{Although MPP-L ultimately achieves the highest accuracy, its advantage over PDEformer-2 emerges only after $\sim 200$ seconds (platform-rescaled; see Supplementary Note~\xref{S4.3}{app:prediction_speed}).}

\paragraph{Prediction speed}
\begin{table*}[t]
	\centering
	\caption{Comparison of parameter counts and per-sample inference time across solvers.
		Neural models are evaluated on a single NPU (batch size 1) by default, and Dedalus runs on a CPU.
		PDEformer-2 base and fast variants are merged when solely preparing an INR,
		as their model sizes differ only in the INR.
	}
	\label{tab:inference_time}
	\tablefont
	\begin{tabular*}{\textwidth}{@{\extracolsep\fill}ll@{\hspace{1em}}ccc@{\hspace{1em}}ccc@{\extracolsep\fill}}
		\toprule
		\multicolumn{2}{c}{\textbf{Model}}
			& \multicolumn{3}{c}{\textbf{Sine-Gordon}}
			& \multicolumn{3}{c}{\textbf{INS-Pipe}} \\
		& & Param. & Time (s) & Error & Param. & Time (s) & Error \\
		\midrule
		\multicolumn{8}{l}{\textbf{Full spatio-temporal output}} \\
		& DeepONet (CNN) & 2.56M & 0.085 & \verB{0.312} & 3.09M & 0.136 & \verB{0.385} \\
		& DeepONet (MLP) & 3.68M & 0.086 & \verB{0.435} & 6.30M & 0.138 & \verB{0.336} \\
		& FNO2D & 0.94M & 0.015 & \verB{0.083} & 0.95M & 0.014 & \verB{0.145} \\
		& FNO3D & 22.12M & 0.062 & \verB{0.051} & 22.12M & 0.062 & \verB{0.134} \\
		& U-Net2D & 13.40M & 0.013 & \verB{0.289} & 13.41M & 0.012 & \verB{0.258} \\
		& MPP-B & 115.6M & \verB{3.367*} & \verB{0.430} & 115.6M & \verB{2.800*} & \verB{0.134} \\
		& \verB{MPP-L} & \verB{407.1M} & \verB{11.10*} & \verB{0.169} & \verB{407.1M} & \verB{11.01*} & \verB{0.0536} \\
		& Poseidon-B & 158M & \verB{0.241}* & \verB{0.712} & 158M & \verB{0.242}* & \verB{0.629} \\  
		& \verB{Poseidon-L} & \verB{629M} & \verB{0.506}* & \verB{0.685} & \verB{629M} & \verB{0.506}* & \verB{0.635} \\  
		& PDEformer-2-base & 82.65M & 0.631 & 0.018 & 82.65M & 1.262 & \verB{0.102} \\
		& PDEformer-2-fast & 71.07M & 0.202 & 0.027 & 71.07M & 0.406 & \verB{0.118} \\
		\addlinespace[0.5em]
		\multicolumn{8}{l}{\textbf{Final time step only}} \\
		& DeepONet (MLP) & 3.68M & 0.005 & --- & 6.30M & 0.006 & --- \\
		& Poseidon-B & 158M & \verB{0.072*} & --- & 158M & \verB{0.072*} & --- \\  
		& \verB{Poseidon-L} & \verB{629M} & \verB{0.074*} & \verB{---} & \verB{629M} & \verB{0.072*} & \verB{---} \\  
		& PDEformer-2-base & 82.65M & 0.033 & --- & 82.65M & 0.052 & --- \\
		& PDEformer-2-fast & 71.07M & 0.033 & --- & 71.07M & 0.050 & --- \\
		\addlinespace[0.5em]
		\multicolumn{8}{l}{\textbf{Solely preparing an INR}} \\
		& PDEformer-2 & --- & 0.020 & --- & --- & 0.022 & --- \\
		\addlinespace[0.5em]
		\multicolumn{8}{l}{\textbf{Traditional solvers}} \\
		& Dedalus & --- & 2.398 & --- & --- & 24.300 & --- \\
		& Dedalus (resolution 64) & --- & 0.936 & 0.029 & --- & 10.131 & 0.026 \\
		& Dedalus (resolution 32) & --- & 0.436 & 0.085 & --- & 5.222 & 0.079 \\
		& Dedalus (resolution 16) & --- & 0.207 & 0.197 & --- & 2.569 & 0.216 \\
		& Dedalus (resolution 8) & --- & 0.092 & 0.396 & --- & 1.058 & 0.459 \\
		\bottomrule
	\end{tabular*}
	\par
	{\raggedright\verB{\footnotesize *Platform-rescaled estimates (MindSpore+Ascend~910B NPU equivalent).
	See Supplementary Note~\xref{S4.3}{app:prediction_speed} for the calibration methodology and original PyTorch+V100 GPU measurements.}\par}
\end{table*}
A fast solver is essential for real-time and multi-query applications.
Although PDEformer-2 sacrifices some speed for flexibility in handling irregular domains and discretizations
through its \verB{arbitrary-coordinate-query} design, as well as the capacity to learn diverse PDEs concurrently,
it remains faster than traditional solvers.
Parameter counts, inference time, and nRMSE for the Sine-Gordon and INS-Pipe datasets are reported in Table~\ref{tab:inference_time}.
We distinguish three prediction modes: full spatio-temporal output ($128^2\times 101$),
final time step only ($128^2\times 1$), and solely preparing an INR for future spatio-temporal queries.
Prediction errors are reported only for full spatio-temporal output.

\paragraph{Dealing with imperfect PDE specifications}
PDEs may be imperfectly known in many practical applications.
Unlike standard neural operators, PDEformer-2 accepts an explicit PDE form as input, enabling it to incorporate partial physical knowledge while adapting to data.
We ablate the PDE input quality by: (1) replacing conservation form by convection (equivalent),
(2) removing viscosity, pressure, and incompressible condition (incomplete),
(3) marking the PDE as unknown, (4) providing Lorenz system (incorrect).
As shown in Fig.~\ref{figs:ablation}c, zero-shot accuracy declines with reduced PDE information, but finetuning with data compensates effectively.
Performance gaps diminish with increasing samples.
Training from scratch yields poor results even with correct PDE input, underscoring the necessity of pretraining.
We note that PDEformer-2 demonstrates the ability to integrate known physics, scarce data, and cross-system knowledge:
with incomplete PDE specifications, it outperforms models using unknown PDE specifications;
finetuning with scarce data improves upon zero-shot results obtained with the same incomplete specifications;
and its pretrained cross-system knowledge enables it to outperform models trained from scratch.
This capability is valuable for modeling partially known systems.

\subsection{Application to inverse problems}\label{sec:results_inverse}
\begin{figure*}[!t]
	\centering
	\includegraphics[width=400pt]{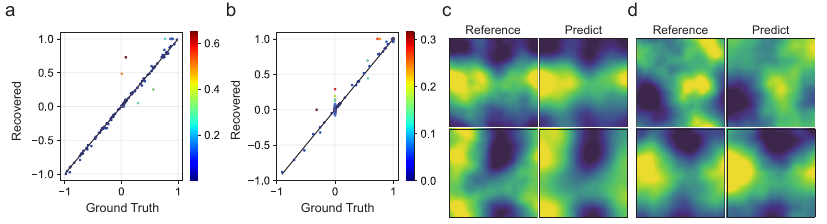}
	\caption{Results of inverse problems using PDEformer-2 as the surrogate solver.
		(a) Recovery of scalar coefficients in DCR equations.
		Each coefficient is depicted as a point in the figure, with the horizontal axis representing the ground-truth value and the \verB{vertical} axis the recovered value.
		The recovery is considered successful when the points are closely aligned with the diagonal line $y=x$.
		(b) System identification, with the same plotting scheme as (a).
		(c) Recovery of source fields in DCR equations, showing two independent cases.
		(d) Recovery of wave velocity fields in Wave equations, showing two independent cases.
	}%
	\label{figs:results_inverse}
\end{figure*}
To demonstrate the applicability of our pretrained model to downstream tasks, we conduct inverse-problem experiments as case studies,
with PDEformer-2-base serving as the surrogate solver.
The tasks involve recovering unknown scalar coefficients or source fields in DCR equations,
as well as velocity fields in Wave equations.
In these scenarios, multiple solutions of a PDE are obtained with different initial conditions.
Rather than having access to complete solution fields, we assume access only to noisy observations at a limited number of spatio-temporal points, and the available initial conditions are also noisy.
The model input comprises the equation form, noisy initial conditions, and the current estimates of the unknown PDE parameters, whether scalar coefficients or fields.
We compare the predicted solutions against observations to compute the nRMSE,
and solve an optimization problem to obtain the estimate.
\paragraph{Recovery of scalar coefficients}
We assume the forms of the PDEs are known, and aim to recover the unknown coefficient values in each equation from noisy observations at random spatio-temporal locations.
PDEs are generated as in the DCR dataset used for pretraining.
For example, for a PDE of the form
\begin{equation}\label{eq:inv_dcr_eg}
	\begin{aligned}
		u_t -a\Delta u + s(\boldsymbol{r}) +c_{01}u + c_{02}u^2 + c_{03}u^3\\
		+ (c_{12}u^2+u^3)_x + (u + c_{23}u^3)_y=0,
	\end{aligned}
\end{equation}
the goal is to recover coefficients $a$, $c_{01}$, $c_{02}$, $c_{03}$, $c_{12}$, and $c_{23}$.
For each PDE, solutions from $25$ different initial conditions are used.
Observations are limited to $128$ random spatial points, and the accessible timesteps at these locations are also random, with an average of $20$ time points per location.
Both observations and initial values are corrupted with 1\% additive noise.
Given the non-convex loss landscape, Particle Swarm Optimization (PSO)~\cite{PSO} is employed.
Figure~\ref{figs:results_inverse}a compares the recovered and true coefficients for $40$ distinct PDEs, showing successful recovery for most cases.
The number of coefficients per PDE varies from 1 to 10, hence the number of points in the figure exceeds $40$.
Additional experiments with varying noise levels, numbers of observations, and numbers of initial conditions are presented in Supplementary Note~\xref{S4.4.1}{subsec:inv-coef}.


\paragraph{System identification}
System identification is a structured approach used in control engineering and signal processing to develop mathematical models of dynamic systems by analyzing observed data.
By assuming knowledge of a complete dictionary of all potential terms in the PDE,
we treat system identification as a particular instance of coefficient recovery,
and need to identify which coefficients in the complete PDE are zero.
In our experiment, the complete PDE takes the form
\begin{equation*}
	\begin{aligned}
		u_t - a\Delta u + s(\boldsymbol{r}) + c_{01}u + c_{02}u^2 + c_{03}u^3 & \\
		+ \left(c_{11} u + c_{12}u^2 + c_{13}u^3 \right)_x & \\
		+ \left(c_{21} u + c_{22}u^2 + c_{23}u^3 \right)_y &= 0
	,\end{aligned}
\end{equation*}
and we recover the coefficients
$a,\ c_{ij}\ (i=0,1,2; j=1,2,3)$ to identify the underlying dynamic system.
The numbers of initial conditions and observations, as well as the noise levels, are the same as those for scalar coefficient recovery.
Results for 10 DCR cases are shown in Fig.~\ref{figs:results_inverse}b.
Without explicit sparsity constraints or penalties, the method correctly identifies most zero coefficients; only 6 out of $29$ recovered values have an absolute magnitude exceeding $0.1$.
Detailed results are presented in Supplementary Table~\xref{S8}{tab:sys_id}.

\paragraph{Recovery of source fields}
We recover the source field $s(\boldsymbol{r})$ in DCR equations from noisy solution measurements.
Equation~\eqref{eq:inv_dcr_eg} gives one such DCR equation.
The numbers of initial conditions and observations, as well as the noise levels, are the same as those for scalar coefficient recovery.
Let $\verB{G_\theta}$ denote the forward operator mapping sources to PDE solutions with PDEformer-2-base as the backbone,
$S^{\text{obs}}$ denote the observable spatio-temporal locations, and $u_{\text{obs}}$ denote the observed solutions.
We solve the optimization problem
\begin{equation}\label{eq:inverse-opt}
	\operatorname*{arg\,min}_{s}\mathcal{L}(s) =  \Bigl\|u_{\text{obs}} - \verB{G_\theta}(s)|_{S^{\text{obs}}}\Bigr\|_2^2 + \lambda \mathcal{R}(s)
\end{equation}
using Adam~\cite{Adam},
where $\mathcal{R}$ represents a regularization term that discourages over-oscillatory patterns
(details in Supplementary Note~\xref{S4.4.3}{subsec:inv-field}).
The average nRMSE over $40$ test sources is $0.159$.
Two example recoveries are shown in Fig.~\ref{figs:results_inverse}c.

\paragraph{Recovery of wave velocity fields}
We address a simplified full waveform inversion task: estimating the wave velocity field $c(\boldsymbol{r})$ (or $a(\boldsymbol{r})=c(\boldsymbol{r})^2$) from wave propagation data.
The PDEs are generated as in the Wave dataset used for pretraining, with inhomogeneous velocities.
For instance, if the PDE takes the form
\begin{equation*}
	\begin{aligned}
		u_{tt} + \mu(\boldsymbol{r})u_t - \nabla \cdot(a(\boldsymbol{r}) \nabla u) + s(\boldsymbol{r}) + c_{02}u^2 \\
		+ (u + u^2)_x + (c_{21}u + u^3)_y = 0,
	\end{aligned}
\end{equation*}
we recover $a(\boldsymbol{r})$ based on observations and known initial conditions.
For each velocity field, solutions are generated for $100$ different combinations of initial conditions $u(0,\boldsymbol{r}),u_t(0,\boldsymbol{r})$ and source fields $s(\boldsymbol{r})$.
The number of observations and the noise levels are the same as in the preceding experiments.
The inverse problem is formulated similarly to~\eqref{eq:inverse-opt}.
To mitigate local minima, we employ multi-start optimization with Adam.
See Supplementary Notes~\xref{S3.3.2}{app:inv_wave} and~\xref{S4.4.4}{subsec:fwi} for details.
Two reconstructed examples are provided in Fig.~\ref{figs:results_inverse}d.

\section{Conclusion}\label{sec:conclusion}
We present PDEformer-2, a foundation model for 2D PDEs.
The flexible model architecture incorporates both symbolic expression and numerical information (initial conditions, coefficients, etc.), and \verB{supports mesh-free solution queries} in full space-time.
Pretraining on a large, diverse dataset enables its broad applicability across varying equation types, domain geometries, boundary conditions, equation and variable counts, as well as time dependence.

PDEformer-2 exhibits preliminary \verB{compositional generalization capability across PDEs}, demonstrated through zero-shot predictions on equations beyond its training distribution,
including those with more variables (Fig.~\ref{fig:PDEformerV2Arch}c)
or novel combinations of physical mechanisms (Fig.~\ref{figs:ablation}a, INS-Tracer).
It avoids inheriting artifacts present in solver-generated pretraining data (Fig.~\ref{figs:results_forward}e).
When adapted to new PDEs, the model outperforms specialized baselines in low-data regimes (Fig.~\ref{figs:ablation}a), indicating competence in data-scarce scenarios.
It also serves as an efficient surrogate for inverse applications that recover unknown PDE parameters,
including scalar coefficients and fields (Fig.~\ref{figs:results_inverse}).

Several limitations remain and require future improvements.
Full-solution inference speed lags behind some baselines \verB{(Table~\ref{tab:inference_time})}, partly owing to its \verB{arbitrary-coordinate-query} design and large capacity.
This may be addressed via more efficient architectures or knowledge distillation into smaller models.
Solution accuracy declines when PDEs deviate substantially from pretraining, making the model inferior to data-sufficient baselines and unsuitable for applications
\verB{(Fig.~\ref{figs:ablation}a, non-square and oscillatory Wave-C-Sines)}.
Performance could be enhanced through more diverse pretraining data, advanced architectures with multi-scale capabilities, and improved finetuning strategies.
Future work will also extend to more complicated domains, time-dependent coefficients and boundaries, as well as 3D settings.
We expect that integrating data from more traditional solvers can encapsulate a wider spectrum of PDE-solving expertise,
thereby reducing the need for engineers to switch solvers and tune parameters.
Validation on real-world problems and applications in design and control remain open challenges.
\verB{Practical deployment also calls for assessing prediction reliability in the absence of reference solutions from conventional solvers.}

\section{Method}\label{sec:method}
\subsection{Problem settings}
We aim to develop a foundation surrogate model for solving almost arbitrary two-dimensional PDEs.
The model input encompasses all information needed to specify a PDE problem.
This includes the symbolic formulations of the PDE and boundary conditions (BCs), the computational domain shape,
coefficients in scalar or field form, and the values or functions prescribed as initial or boundary conditions.
The output is the complete solution that can be queried at any spatio-temporal point.
The number of equations, BCs, coefficients, and solution components may vary arbitrarily across PDEs.
\verB{See Supplementary Note~\xref{S1.1.1}{app:problem_settings} for details.}

\subsection{Model overview}
Figure~\ref{fig:PDEformerV2Arch}a illustrates the overall network architecture of PDEformer-2,
which is based on our previous PDEformer-1~\cite{PDEformer1} model.
As shown in the figure, PDEformer-2 first formulates the symbolic expression of the PDE as a computational graph,
and employs a scalar encoder and a function encoder to embed the numerical information of the PDE into the graph node features.
A graph Transformer processes the graph, and the resulting latent embeddings are decoded by an INR, producing predicted values of each solution component at specific spatio-temporal coordinates.
In our implementation, the scalar encoder is a multi-layer perceptron (MLP), the function encoder is a convolutional neural network (CNN), the graph Transformer is based on the Graphormer architecture~\cite{graphormer}, and the INR decoder is an adapted version of Poly-INR~\cite{Poly-INR}.
Further architectural details are provided in Supplementary Note~\xref{S1}{app:model}.

\subsection{Domain shapes and boundary conditions}
To handle arbitrary domain geometries in two dimensions, PDEformer-2 represents domains and boundary locations as signed distance functions (SDFs).
For a subset $S\subseteq[0,1]^2$, the SDF is defined as
\[\mathrm{SDF}_S(\boldsymbol{r})=\begin{cases}
	\inf_{\boldsymbol{r}'\in S}d(\boldsymbol{r},\boldsymbol{r}'),&\boldsymbol{r}\notin S,\\
	-\inf_{\boldsymbol{r}'\notin S}d(\boldsymbol{r},\boldsymbol{r}'),&\boldsymbol{r}\in S
.\end{cases}\]
Functions defined on $S$ are extended to $[0,1]^2$ via ``nearest neighbor'':
\begin{equation}\label{eq:field_extension}
	f^\text{ext}(\boldsymbol{r})=f\left(\operatorname*{arg\,min}_{\boldsymbol{r}'\in S}d(\boldsymbol{r},\boldsymbol{r}')\right)
,\end{equation}
enabling encoding through the function encoder.
(At certain locations, further technical conventions may be required to make this extension well-defined.
We omit them because they rarely cause practical issues in discretized computation.)
Figure~\ref{fig:PDEformerV2Arch}b shows an example computational graph of Dirichlet boundary conditions.
Additional examples can be found in Supplementary Note~\xref{S2.1}{app:DAG_BC}.

\section*{Data availability}
The pretraining dataset, three custom adaptation datasets (Sine-Gordon, INS-Tracer, INS-Pipe),
and the inverse-problem datasets
can be accessed via \url{https://github.com/functoreality/pdefoundry-2}.
The two publicly available datasets include Wave-Gauss at \url{https://huggingface.co/datasets/camlab-ethz/Wave-Gauss}
and Wave-C-Sines at \url{https://zenodo.org/records/14765453/files/Wave-C-Sines.nc}.
\section*{Code availability}
The source code of this paper is available in the GitHub repository at \url{https://github.com/functoreality/pdeformer-2}.
The pretrained model weights for PDEformer-2-base are available at \url{https://gitee.com/functoreality/PDEformer2-L},
and those for PDEformer-2-fast are available at \url{https://gitee.com/functoreality/PDEformer2-M}.

\section*{Acknowledgements}
Z.Y.'s contribution to this work comprised two phases, the latter of which was completed during an internship
at Huawei Technologies Co., Ltd.
We thank Wendao Wu for early-stage analysis of the non-physical Dedalus solutions and constructive discussions.

\section*{Funding}
This work is supported in part by the National Science and Technology Major Project (2022ZD0117804).
Bin Dong is supported in part by the New Cornerstone Investigator Program.

\section*{Author contributions}
Z.Y.: computational graph representation design \& implementation, data generation---Dedalus, exploration for architecture \& pretraining \& finetuning, mainstream pretraining, adaptation---major, writing---major
\\Z.L.: exploration for architecture \& pretraining \& finetuning, inverse problems, prediction speed evaluation, solution visualization, writing---major
\\B.W.: data generation---FEniCSx, exploration for architecture \& pretraining \& finetuning, analysis related to Dedalus artifacts \& failure cases, adaptation and speed evaluation of Poseidon \& Geo-FNO, writing---related sections
\\H.J.: data generation---Dedalus MV-Wave, exploration for architecture \& pretraining \& finetuning, baseline investigation, adaptation and speed evaluation of MPP, writing---related sections
\\L.C.: early-stage exploration for architecture \& training, investigation for dataset design
\\M.Z.: data generation---solver investigation \& supporting early-stage exploration \& improving implementation, baseline investigation, software \& hardware infrastructure maintenance
\\X.H.: early-stage exploration for architecture \& training, early-stage implementation
\\Q.M.: baseline investigation, improving implementation for parallel training, software \& hardware infrastructure maintenance
\\J.Z.: assisting results analysis, software \& hardware infrastructure maintenance
\\H.L.: project administration, design of code \& data repository structure
\\B.D.: project administration, writing---overall structure \& improving text logic \& review

\noindent\textbf{\textit{Conflict of interest statement.}} None declared.
\bibliographystyle{unsrt}
\bibliography{Ref}

\clearpage
\onecolumn
\appendix
\renewcommand{\thesection}{S\arabic{section}}
\renewcommand{\thesubsection}{\thesection.\arabic{subsection}}
\renewcommand{\thesubsubsection}{\thesubsection.\arabic{subsubsection}}
\renewcommand{\theparagraph}{\thesubsubsection.\arabic{paragraph}}
\renewcommand{\thefigure}{S\arabic{figure}}
\renewcommand{\thetable}{S\arabic{table}}
\setcounter{figure}{0}
\setcounter{table}{0}
\section*{Supplementary Information}
\addcontentsline{toc}{section}{Supplementary Information}
This supplementary information provides additional details on several key aspects of our work.
Supplementary Note~\ref{app:model} outlines the architecture of PDEformer-2 and the baseline models used in our experiments.
Supplementary Note~\ref{app:PDEDAGdetails} describes the construction of DAGs used as inputs to PDEformer-2, derived from various PDE forms.
Supplementary Note~\ref{app:dataset} presents all datasets employed in the experiments, including the pretraining dataset and those used for finetuning.
Supplementary Note~\ref{app:training} details the training settings, including procedures, hyperparameters, and further results on inverse problems.
Supplementary Note~\ref{app:code_example} provides code examples demonstrating how to implement the PDEformer-2 model for solving different types of PDEs.
Finally, Supplementary Note~\ref{app:visual} offers additional visualizations comparing the predictions of PDEformer-2 with those of other baseline models on the PDEs discussed in the main text.

\section{Detailed model architecture}\label{app:model}
We shall present the detailed model architectures in this section, including PDEformer-2 and the baseline models.
PDEformer-2 retains the overall architecture of PDEformer-1 \cite{PDEformer1}, with modifications made only to the input dimension, the function encoder architecture, the activation function in INR, as well as some network hyperparameters (width and depth).
However, for the convenience of readers, we provide the complete model architecture of PDEformer-2 here.

\subsection{Overall pipeline of PDEformer-2}
\subsubsection{Detailed problem settings}\label{app:problem_settings}
We consider two-dimensional PDEs defined on $(t,\boldsymbol{r})\in[0,1]\times\Omega$ of the general form
\begin{equation}\label{eq:general_pde}\begin{split}
	\mathcal{F}(u_1,u_2,\dots,c_1,c_2,\dots,s_1(\boldsymbol{r}),s_2(\boldsymbol{r}),\dots)&=0\quad\text{in }\Omega,\\
	\mathcal{B}_i(u_1,u_2,\dots,c_{i1},c_{i2},\dots,s_{i1}(\boldsymbol{r}),s_{i2}(\boldsymbol{r}),\dots)&=0\quad\text{on }\Gamma_i,
\end{split}\end{equation}
where $\boldsymbol{r}=(x,y)\in\Omega\subseteq[0,1]^2$ is the spatial coordinate, $c_1,c_2,\dots,c_{11},c_{12},\dots \in \mathbb{R}$ are real-valued scalar coefficients, $s_1(\boldsymbol{r}),s_2(\boldsymbol{r})\dots,s_{11}(\boldsymbol{r}),\dots$ are scalar-valued functions (which may serve as initial conditions, boundary values or coefficient fields in the equation), and $u_1,u_2,\dots:[0,1]\times\Omega\to\R$ are unknown field variables to be solved in the equation.
The boundary conditions are indexed by $i=1,2,\dots$.
Here, we assume that each of the operators $\mathcal{F},\mathcal{B}_1,\mathcal{B}_2,\dots$ admits a symbolic expression, which may involve differential and algebraic operations.
The goal of PDEformer-2 is to construct a surrogate model of the solution mapping
\[(\Omega,\mathcal{F},c_1,\dots,s_1(\boldsymbol{r}),\dots,\Gamma_1,\mathcal{B}_1,c_{11},\dots,s_{11}(\boldsymbol{r}),\dots)\mapsto(u_1,u_2,\dots).\]
The input of this solution mapping includes the location of the computational domain $\Omega$ and the boundaries $\Gamma_1,\Gamma_2,\dots$,
the symbolic expressions of the interior operator $\mathcal{F}$ and boundary operators $\mathcal{B}_1,\mathcal{B}_2,\dots$,
as well as the numeric information $c_1,\dots,c_{11},\dots,s_1(\boldsymbol{r}),\dots,s_{11}(\boldsymbol{r}),\dots$ involved.
The output includes all components of the predicted solution, i.e., $u_1,u_2,\dots:[0,1]\times\Omega\to\R$.
For time-independent steady-state equations, we restrict the problem to $t=0$ for simplicity, which requires only the prediction of $u_1,u_2,\dots:\{0\}\times\Omega\to\R$.

\verB{
\paragraph{Difference from a conditioned operator}
We emphasize that the above problem setting should not be viewed as a surrogate operator conditioned on the PDE form.
In a conditioned neural operator framework in the usual sense, the solver would typically take the form
\begin{equation*}
	G_{\mathcal{F},\mathcal{B}_1,\dots}\colon\R^{d_c}\times \mathcal{A}^{d_s}\to \mathcal{U}^{d_u},\quad (c_1,\dots,s_1(\boldsymbol{r}),\dots)\mapsto(u_1,u_2,\dots),
\end{equation*}
where $\mathcal{A}=\mathcal{A}(\Omega)$ and $\mathcal{U}=\mathcal{U}([0,1]\times\Omega)$ are certain function spaces,
and the mapping from numeric inputs to solutions depends on the PDE form $\mathcal{F},\mathcal{B}_1,\dots$.
However, in our setting, the input and output dimensions $d_c$, $d_s$, and $d_u$ directly depend on the symbolic operators $\mathcal{F},\mathcal{B}_1,\dots$, rendering such a formulation impractical for standard neural operator architectures.
Moreover, the semantic interpretation of each input and output component varies with the specific PDE form, and the ordering of these components lacks a canonical convention that can be decoupled from the symbolic information.

\begin{figure}[htbp]
	\centering
	\includegraphics[width=0.6\linewidth]{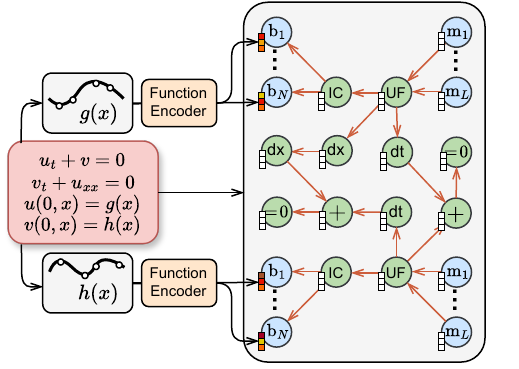}
	\caption{Computational graph example of a multi-variable PDE, showing a 1D version for simplicity.
		Two different nodes of type \texttt{UF} have been introduced to represent the two variables $u$ and $v$, respectively.
		Adapted from our prior work~\cite{PDEformer1}.
	}
	\label{fig:DAG_mCompn}
\end{figure}
Our approach directly addresses these challenges by encoding the symbolic structure alongside the numeric information within a unified computational graph representation.
Such a computational graph plays a central role in our architecture, which constitutes the structural backbone that governs how PDE information flows through the network.
The introduction of additional scalar coefficients, coefficient fields, or unknown variables is reflected by additional nodes in such a computational graph,
enabling the model to handle PDEs involving varying $d_c,d_s$, and $d_u$ without a pre-selected upper bound.
See Fig.~\ref{fig:DAG_mCompn} for a computational graph example with more than one unknown variable.
This design also supports extrapolation to a greater number of output components $d_u$ than seen during training, as demonstrated in our multi-variable generalization experiments (Fig.~\ref{fig:PDEformerV2Arch}(c) in the main text).
}

\subsubsection{Construct graph}
We first represent $\mathcal{F},\mathcal{B}_1,\mathcal{B}_2,\dots$  in~\eqref{eq:general_pde}, i.e., the symbolic information specifying the form of the PDE and BCs, as a computational graph.
In such a computational graph, a node may stand for a mathematical object or an operation,
and a directed edge can be used to specify the operands involved in an operation.
These nodes and edges constitute a directed acyclic graph (DAG) with heterogeneous nodes and homogeneous edges.
Supplementary Note~\ref{app:PDEDAGdetails} provides more details and examples of such a computational graph.

Then, in order to include the numeric information, we endow each graph node with an input feature embedding vector $\xi\in\R^{d_e}$.
This would involve three cases as follows:
(1)
For a scalar coefficient $c$, we input this numeric value into a scalar encoder, and use the $d_e$-dimensional output as the input features for the corresponding \verb|SC| node.
(2)
In terms of a function $s(\boldsymbol{r})$ that typically contains more information, we shall introduce $N$ new ``branch'' nodes, with types denoted as $\mathtt{b}_1,\mathtt{b}_2,\dots,\mathtt{b}_N$, respectively.
These new nodes are then connected to the \verb|CF|, \verb|IC|, or \verb|SDF| node corresponding to $s(\boldsymbol{r})$.
We feed $s(\boldsymbol{r})$ into a function encoder, and use the output to generate the $N$ input features of these branch nodes.
(3)
As for the remaining nodes, zero input features can be assigned conceptually, as they have no numeric information attached.
To simplify implementation, we actually feed zero into the scalar encoder, and take its output as the input features of these nodes.

To improve model efficiency, the function $s(\boldsymbol{r})$ is not represented as scattered points as in PDEformer-1.
Instead, we assume it is discretized on the $128\times 128$ uniform grid of the entire domain $[0,1]^2$, and this may require function extension (according to~\eqref{eq:field_extension} in the main text) and interpolation in \verB{practice}.
The grid values are fed into a function encoder with a convolutional neural network (CNN) backbone, generating a single output vector with dimension $d_eN$.
We then split this output vector into $N$ parts, which form all the $d_e$-dimensional input features of the branch nodes.
The value $N=4$ is found to be sufficient in experiments, and is used in our PDEformer-2-base and PDEformer-2-fast models.

Moreover, for each unknown field variable to be solved, we introduce $L$ additional nodes with types $\mathtt{m}_1,\mathtt{m}_2,\dots,\mathtt{m}_L$, respectively, and connect them to the corresponding \verb|UF| node.
These nodes will be used to decode the predicted solution, as will be explained in the following text.
The specific value of $L$, being determined by the depth of the INR decoder, equals $11$ for the PDEformer-2-base and PDEformer-2-fast models.

\subsubsection{Encode graph data}\label{app:arch_graph_encode}
The graph data obtained in the previous step contains both symbolic and numeric information inherent in the PDE.
We shall integrate the information from the graph data, and generate a latent code that represents the solution for each field variable $u_j$ to be solved.
Each latent code takes the form
\[\pmb{\mu}_j = [{\mu}^1_j, \dots, {\mu}^L_j]^{\mathrm{T}} \in\R^{L \times d_e}.\]

The fusion process leverages a graph Transformer,
a powerful type of graph neural network based on the Transformer architecture,
which is skilled at capturing and expressing complex graph structural information.
In practice, we adopt the Graphormer architecture~\cite{graphormer}, and make some adjustments to adapt it specifically for encoding PDEs.
For each variable index $j$ and INR layer index $\ell$, the \verb`UF` node representing $u_j$ is connected with another node with type $\mathtt{m}_\ell$.
We let $\mu_j^\ell\in\R^{d_e}$ be the embedding vector assigned to this $\mathtt{m}_\ell$ node in the output layer of the graph Transformer.

\verB{
\paragraph{Necessity of a ``graph'' Transformer}
One might ask whether a standard sequence Transformer could replace the graph Transformer by simply flattening the graph nodes into a node sequence.
However, such a scheme would lose the directed edge information, which is essential for fully specifying the PDE form.
For instance, consider the two equations $u_t+cu_x=0$ and $u_t+u_x+cu=0$: their computational graphs share the same set of nodes in the same order
(\texttt{UF},\texttt{dt},\texttt{dx},\texttt{SC},$\times,+,=0$),
differing only in how directed edges connect them.
A sequence Transformer, having no access to edge connectivity, would be fundamentally unable to distinguish between such PDEs that are mathematically distinct.

Beyond encoding directed edges, the graph Transformer enjoys an additional benefit: it is invariant to the ordering of nodes.
This naturally aligns with several important inductive biases in PDE modeling.
Addition and multiplication are commutative and associative; the order of unknown variables in a system of equations is arbitrary; and the order of equations within a system is similarly inconsequential.
A sequence Transformer, being inherently order-sensitive, lacks these inductive biases at the architectural level.
Instead, it needs to learn them from data alone, which would require explicit data augmentation (e.g., permuting terms or equations) during training.
}

\subsubsection{Decode PDE solution}
We employ an INR that takes the spatio-temporal coordinate $(t,x,y)$ as input, and produces the mesh-free prediction $\hat{u}_j(t,x,y)$ according to $\pmb{\mu}_j$ for each field variable.
In consistency with PDEformer-1~\cite{PDEformer1}, we utilize an adapted version of Poly-INR~\cite{PolyINR} with $L=11$ hidden layers,
and the modulations of the $\ell$-th hidden layer are generated according to $\mu_j^\ell$.

\subsection{Detailed PDEformer-2 architecture}
\subsubsection{Scalar and function encoder}
The scalar encoder is a multi-layer perceptron (MLP) that has two hidden layers with $256$ neurons each.
The input has dimension one (i.e., a scalar input), and the output dimension equals $d_e=768$, the embedding dimension of the graph Transformer.

To encode functions defined on 2D uniform grids, we employ a lightweight convolutional neural network (CNN).
The input has a single channel and a resolution of $128\times 128$.
Our CNN consists of the following layers:
\begin{enumerate}
	\item 2D convolution with $32$ output channels, stride $4$, kernel size $4\times 4$, with bias included.
	\item Rectified linear unit (ReLU) activation function.
	\item 2D convolution with $128$ output channels, stride $4$, kernel size $4\times 4$, with bias included.
	\item ReLU activation function.
	\item 2D convolution with $768$ output channels, stride $4$, kernel size $4\times 4$, with bias included.
	\item Flattening operation, applied to the $2\times 2$ feature map with $768$ channels.
	\item Vector \verB{splitting} into $N=4$ feature embeddings assigned to the branch nodes.
\end{enumerate}

\subsubsection{Graph Transformer}
\paragraph{Initial embedding vector}
In the graph Transformer, the initial embedding of node $i$ is given as
\[h_i^{(0)} = x_{\text{type}(i)} + \xi_i + z^-_{\text{deg}^-(i)} + z^+_{\text{deg}^+(i)},\]
where $\xi_i\in\R^{d_e}$ is the input feature embedding vector of node $i$, and $x,z^-,z^+\in\R^{d_e}$ are learnable vectors specified by node type $\text{type}(i)$, in-degree $\text{deg}^-(i)$ and out-degree $\text{deg}^+(i)$, respectively.

\paragraph{Attention bias}
We denote $\phi(i,j)$ to be the shortest path length from node $i$ to node $j$.
If such a path does not exist or has a length greater than $14$, we set $\phi(i,j) = 14$.
For each attention head involved in the graph Transformer, the attention bias according to the node pair $(i,j)$ is given as
\begin{equation}\label{eq:attn_bias}
	B_{ij} = b^+_{\phi(i,j)} + b^-_{\phi(j,i)} + d_{ij}
.\end{equation}
Here, $b^+_{\phi(i,j)},b^-_{\phi(i,j)}$ are learnable scalars indexed by $\phi(i,j)$ and $\phi(j,i)$, respectively, and are shared across all layers.
The additional term $d_{ij}$, which does not appear in the original Graphormer, is introduced to mask out attention between disconnected node pairs.
More specifically, when node $i$ and node $j$ are connected in the graph, i.e., there exists a path either from $i$ to $j$ or from $j$ to $i$, we take $d_{ij}=0$, and set $d_{ij}=-\infty$ otherwise.
Moreover, since our graph has homogeneous edges, we do not include the edge encoding term that appears in the original Graphormer.

\paragraph{Graph Transformer layer}
The structure of the graph Transformer layer is the same as the original Graphormer, and we include it here for convenience to the readers.
Each layer takes the form
\begin{align}
	\bar h^{(l)} &= \text{Attn}(\text{LN}(h^{(l-1)})) + h^{(l-1)}\\
	h^{(l)} &= \text{FFN}(\text{LN}(\bar h^{(l)})) + \bar h^{(l)}
,\end{align}
where $\text{FFN}$ represents a position-wise feed-forward network with a single hidden layer and GeLU activation function, and $\text{LN}$ stands for layer normalization.
In terms of the self-attention block $\text{Attn}$, we shall follow the convention in the original Graphormer paper and only present the single-head case for simplicity.
Let $H = [h_1', \cdots, h_n']^\mathrm{T}\in\R^{n\times d_e}$ denote the input of the self-attention module involving $n$ graph nodes, the self-attention is computed as
\[\begin{aligned}
	Q = HW_Q,\quad K = HW_K,\quad V = HW_V,\\
	A = \frac{QK^\mathrm{T}}{\sqrt{d_e}} + B, \quad \text{Attn}(H) = \text{softmax}(A)V,
\end{aligned}\]
where $W_Q,W_K,W_V\in\R^{d_e\times d_e}$ are the projection matrices, and $B$ is the attention bias given in~\eqref{eq:attn_bias}.
The extension to the multi-head attention is standard and straightforward.

\paragraph{Hyperparameters}
PDEformer-2, including base and fast variants, uses $12$ graph Transformer layers with embedding dimension $d_e=768$, FFN embedding dimension 1536, and $32$ attention heads.
This configuration is selected based on our previous experience of PDEformer-1~\cite{PDEformer1}.

\subsubsection{INR}
\begin{figure}[htbp]
	\centering
	\includegraphics[width=\textwidth]{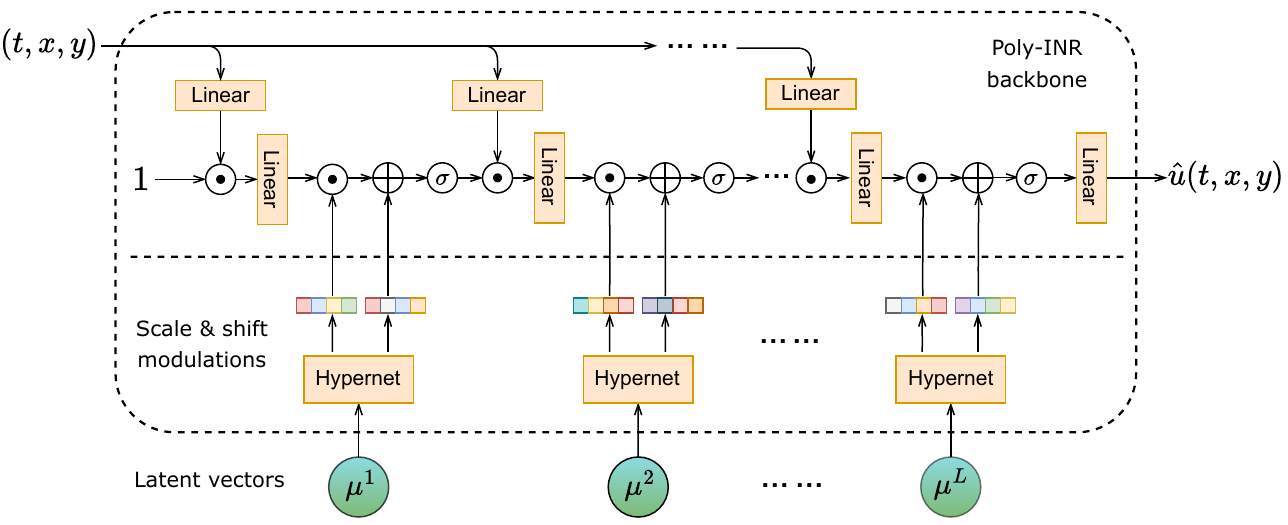}
	\caption{INR decoder architecture of PDEformer-2.}
	\label{fig:INR-decoder}
\end{figure}
In the realm of implicit neural representation (INR), data objects are interpreted as coordinate-based functions, where each function accepts a coordinate $(t,x,y)$ as input and yields an approximated function value $\hat{u}(t,x,y)$ at that specific coordinate point.
Among various INR architectures proposed in the literature,
PDEformer-2 employs an adapted version of Poly-INR~\cite{PolyINR}.
We make a modification to the $L$ hypernets in the Poly-INR structure, in which the $\ell$-th hypernet takes $\mu^\ell\in\R^{d_e}$ as its input, and generates the scale- and shift-modulations for the $\ell$-th hidden layer of our Poly-INR.

The architecture of our INR decoder is illustrated in Fig.~\ref{fig:INR-decoder}, with the mathematical framework detailed below.
We take $h_0=\mathbf{1}$ to be the vector with all entries being equal to one.
For $\ell=1,2,\dots,L$, we compute
\[\begin{aligned}
	g_\ell &= W^{\text{in}}_{\ell} \begin{bmatrix} t\\x\\y \end{bmatrix} + b^\text{in}_\ell, \quad
	s_{\ell}^{\text{scale}} = \text{MLP}_\ell^{\text{scale}}(\mu^\ell), \quad
	s_{\ell}^{\text{shift}} = \text{MLP}_\ell^{\text{shift}}(\mu^\ell), \\
	q_{\ell} &= s_{\ell}^{\text{scale}} \odot \left( W^{\text{h}}_{\ell} \left(h_{\ell-1} \odot g_\ell\right) + b_\ell^\text{h} \right) + s_{\ell}^{\text{shift}}, \quad
	h_\ell = \sigma\left(q_{\ell}\right),
\end{aligned}\]
and the network output is given as $\hat{u}(t,x,y) = W^{\text{Last}}h_{L} + b^{\text{Last}}$.
Here, the activation function $\sigma(\cdot)$ is a sine function with output scaled by $\sqrt 2$, which performs slightly better than the clipped leaky-ReLU operation used in PDEformer-1.
The hypernets correspond to $\text{MLP}_\ell^{\text{scale}}$ and $\text{MLP}_\ell^{\text{shift}}$.
Note that in the original Poly-INR, the hypernets are utilized to generate $W^\text{in}_\ell$ and $b^\text{in}_\ell$.
Compared with our practice of generating $s^\text{scale}_\ell$ and $s^\text{shift}_\ell$, this method makes training slower, and is not adopted in our implementation.

PDEformer-2-base uses an INR containing $L=11$ hidden layers with $768$ neurons each.
Discretization of 2D time-dependent PDEs often involves many grid points (typically $101\times 128\times 128$ in our pretraining dataset).
The INR deals with each grid point independently, which may be computationally intensive for inference.
This can be seen in the difference of the inference time between full prediction and Encoder-Only in Table~\ref{tab:inference_time} in the main text.
To accelerate model inference, PDEformer-2-fast employs a relatively smaller INR, with the width reduced to $256$ from $768$.
This is the only difference between the base and the fast variants.
We keep the number of hidden layers $L$ so as not to affect the computational graph structures.

We focus on 2D PDEs in this work, and the coordinates should only involve $(t,x,y)$.
However, to make our code compatible for future extension to 3D PDEs, the $z$ coordinate is also included in the INR input.
This coordinate remains zero in all the experiments, and does not affect the results.
The same convention applies to the baseline DeepONet trunk networks in the experiments.

\subsection{Baseline models}\label{app:model_baseline}
\paragraph{DeepONet}
We implement two variants of DeepONet~\cite{DeepONet}, each comprising a trunk network and a branch network.
The trunk network is a six-layer MLP of width $256$, and accepts spatio-temporal coordinates as network input.
The branch network accepts the coefficient scalars and fields in PDE as network input, and generates a vector of dimension $n_o=2048$.
In consideration of the multi-output case with $n_v>1$ variables to predict,
we adopt approach 4 mentioned in Section~3.1.6 of~\cite{FAIR},
and split the trunk network output into $n_v$ groups, each with $n_o$ neurons.
For $k=1,2,\dots,n_v$, we select the $k$-th group of the trunk network output, and take inner product with the entire branch network output to get the prediction of the $k$-th variable.
The case of $n_v=1$ would degenerate to the vanilla DeepONet initially proposed in~\cite{DeepONet}.

The direct variant, designated as ``DeepONet (MLP)'', employs a six-layer MLP branch network of width $256$.
The input of the MLP contains the variable coefficient scalars, as well as the grid values of all coefficient fields involved.
The coefficient fields in the datasets are typically discretized on a uniform grid with resolution $128\times 128$.
Such a fine resolution may introduce too many grid points, leading to a high input dimension of the branch MLP.
This can make our NPU run out of memory in the experiments, especially when a number of (for example, four for the INS-Pipe dataset) such coefficient fields are involved.
Therefore, we downsample the coefficient fields to a resolution of $64\times 64$ before feeding them into the MLP branch network.

When DeepONet was initially proposed in~\cite{DeepONet}, the authors have emphasized that the architectures of the branch and trunk networks can be arbitrary.
To better account for the nature of 2D PDEs, we introduce another DeepONet variant with the branch network replaced by a CNN, referred to as ``DeepONet (CNN)''.
The CNN possesses six convolutional layers followed by two fully connected layers.
Each PDE coefficient scalar is treated as a constant coefficient field, and form an independent channel of the CNN input.

\paragraph{U-Net}
The U-Net~\cite{U-Net} model features a symmetric design with five downsampling and five upsampling stages.
The number of channels along the downsampling pathway equals $64$, $128$, $256$, $512$, and $512$, respectively,
while upsampling utilizes $1024$, $512$, $256$, $128$, and $64$.
A linear projection then applies to the feature map with $64$ channels to get the final output.
All layers employ ReLU activation functions.
The network input is arranged in the same way as the branch network of DeepONet (CNN).
As a 2D neural operator, U-Net treats the $n_t$ temporal steps as extended channels, leading to a total of $n_tn_v$ output channels.

\paragraph{FNO}
Our implementation includes both two-dimensional and three-dimensional variants of the Fourier Neural Operator (FNO)~\cite{FNO}.
The input and output of the 2D variant is formulated in the same way as U-Net.
The 3D variant treats the temporal axis as an additional spatial dimension with resolution $n_t$.
The network input, containing both scalar and field channels, is replicated along this new dimension to fit into the network.
Both variants employ four Fourier layers with $20$ hidden channels, and retain $12$ Fourier modes in each spatial dimension.
The final output is generated by a MLP with a single hidden layer of width $128$.
The Wave-Gauss dataset provided by {\PDEgym}~\cite{Poseidon} has only $n_t=16$ timesteps,
and we reduce the number of Fourier modes to $6$ along the temporal axis when applying FNO3D to this dataset.

\paragraph{Geo-FNO}
The Geometry-aware Fourier Neural Operator (Geo-FNO)~\cite{Geo-FNO} extends the standard FNO framework to handle
irregular geometries by deforming physical domains into a uniform latent space compatible with FFT operations.
Our implementation replicates the official Geo-FNO architecture, which consists of a deformation mapping network
and a FNO backbone, which are jointly optimized in an end-to-end manner.

The FNO component employs five Fourier layers with $32$ hidden channels and $12$ Fourier modes per spatial dimension.
The deformation mapping network is a MLP consisting of four hidden layers, with $32$, $128$, $128$, and $128$ neurons, respectively.

The model takes two input tensors: (i) the input function values of size $(n_{xy}, c_\text{in})$, where $c_\text{in}$ is the number of input variables, and
(ii) the coordinates of mesh points in the physical domain of size $(n_{xy}, 2)$.
The output tensors are of size $(n_{xy}, c_\text{out})$, with $c_\text{out} = n_t \times n_v$ accounting for $n_t$ temporal steps and $n_v$ target variables.

\paragraph{Poseidon}
Poseidon~\cite{Poseidon} is a foundation model for learning the solution operators of PDEs,
leveraging a transformer-based architecture to capture complex spatio-temporal dynamics.
The model is designed to approximate the underlying solution operator $\mathcal{S}:[0, T] \times \mathcal{X} \to \mathcal{X}$,
such that $u(t) = \mathcal{S}(t, u_0)$ is the solution of the PDE at any time $t \in [0, T]$ given the initial condition $u_0$.
Additional coefficient scalars and fields (such as variable viscosity and source terms) are treated as additional components (channels) of $u$,
and are excluded from the computation of finetuning loss.
\verB{We evaluate pretrained Poseidon-B and Poseidon-L models in the data-scaling experiments and inference-time comparison.
Poseidon-B is the size-matched comparison to our PDEformer-2-base and PDEformer-2-fast models, while Poseidon-L is additionally reported to show the larger model variant.}

\paragraph{Multiple Physics Pretraining (MPP)}
Multiple Physics Pretraining (MPP)~\cite{MPP} is a pretraining strategy designed for learning from multiple physical systems simultaneously.
It embeds data from different dynamics into a shared representation space while maintaining separation between systems at certain network layers.
This allows the model to learn transferable representations across equations, while still preserving task-specific information.
Similar to Poseidon-B, we choose the MPP-AViT-B model provided by the official MPP repository in our experiments.

During inference, MPP predicts the solution at the next time step only.
To obtain a full trajectory from the initial conditions, predictions must be generated sequentially in an autoregressive manner.
This is very different from PDEformer-2, which can make direct predictions at arbitrary spatio-temporal coordinates, without relying on step-by-step recursive evaluation.

A total of $16$ consequent time snapshots were taken as the model input when MPP was pretrained.
Since the prediction task in our experiments requires generating the whole solution according to the initial state,
we only take a single time snapshot as the model input during our finetuning process.

MPP distinguishes between different physical equations by assigning distinct channels in specific intermediate layers for each equation.
The number of such channels is fixed, which limits the total number of variables that can be accommodated in the pretrained model,
unless we introduce additional network parameters with random initialization.
To enable direct finetuning on additional equation types beyond those in the public model weights,
we selectively enable channel sharing between new equations and those with similar physical characteristics in its original pretraining dataset.

MPP is not designed to handle equations with varying input coefficient fields.
To address this, we treat the coefficient fields as additional components (channels) of the PDE solution.
During prediction, the coefficient field is included in the model output but ignored during the autoregressive updating process.
Its prediction errors are also excluded from the evaluation metrics.

\section{Computational graph representation of PDEs}\label{app:PDEDAGdetails}
The overall idea of representing PDEs as computational graphs is mainly consistent with our previous PDEformer-1~\cite{PDEformer1}.
However, some new ingredients exist, especially when we switch from one-dimensional PDEs to two-dimensional ones.
We shall explain these in detail in this section.

\subsection{Domain shape and boundary conditions}\label{app:DAG_BC}
In order to inform the model of the domain shape $\Omega$, we first convert it to a function $\mathrm{SDF}_\Omega(\boldsymbol{r})$.
Another possible method for such conversion would be using characteristic functions, which takes value one inside $\Omega$ and zero outside $\Omega$.
However, we believe that SDFs have the following advantages:
\begin{itemize}
	\item SDFs contain richer information.
		Even when it is discretized on a grid with relatively low resolution, we can still recover the approximate boundary of $\Omega$.
	\item SDFs are more sensitive to topological changes.
		For example, if we dig a very small hole inside $\Omega$, its SDF changes a lot, but the characteristic function does not change as much, especially in the $L^2$ sense.
		We want to inform the model of the change, since the corresponding solution may change significantly in this case.
	\item Besides domains, edges and line segments with zero area can also be described by SDFs.
		If characteristic functions are used, we may have to first expand the lines to a certain level of thickness, especially when the characteristic function is discretized on a relatively low-resolution grid.
		Specifying boundary conditions requires declaration of the boundary location that are typically lines,
		and SDF provides a unified representation of these domains and boundaries.
\end{itemize}

\begin{figure}[htpb]
	\centering
	\includegraphics[width=0.7\linewidth]{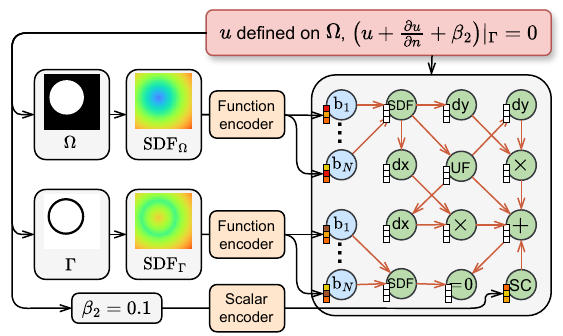}
	\caption{Computational graph example of non-periodic boundary conditions in a PDE.
	The domain $\Omega$ is a disk, and inhomogeneous Robin boundary condition is imposed on its boundary circle.}%
	\label{fig:DAG_BC_Robin}
\end{figure}
\begin{figure}[htpb]
	\centering
	\includegraphics[width=0.8\linewidth]{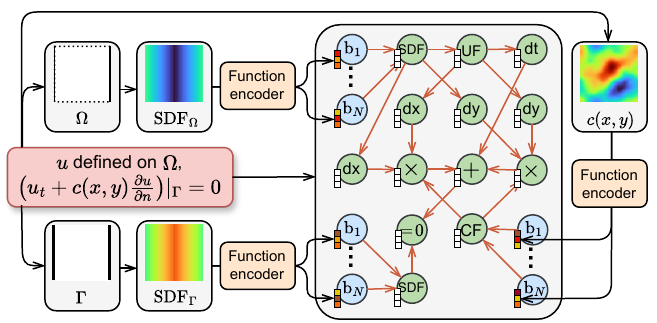}
	\caption{Computational graph example of non-periodic boundary conditions in a PDE.
	The square domain $\Omega$ is assumed to be periodic along the $y$-axis, and the wave-absorbing Mur boundary condition is imposed on the left and right edges.}%
	\label{fig:DAG_BC_Mur}
\end{figure}

In addition to Fig.~\ref{fig:PDEformerV2Arch}(b) in the main text,
Fig.~\ref{fig:DAG_BC_Robin} and~\ref{fig:DAG_BC_Mur} give two additional examples of how domain shapes and boundary conditions can be represented as computational graphs.
We also have the following remarks on the representation:
\begin{itemize}
	\item When the domain $\Omega=[0,1]\times[0,1]$ is periodic along both axes, it is deemed equivalent to the two-dimensional torus $S^1\times S^1$, and the SDF cannot be defined.
		The \verb|UF| nodes do not have a \verb|SDF| predecessor in the computational graph.
		An example is the advection equation in Fig.~\ref{fig:PDEformerV2Arch}(a) in the main text.
	\item If $\Omega$ is periodic along one axis, say $x$ without loss of generality, we shall require%
		\footnote{All datasets in our experiments in this work satisfy $a=0,b=1$ for simplicity.}
		$\Omega=[0,1]\times[a,b]$ for some $0\le a<b\le 1$.
		By identifying $\Omega$ with $S^1\times[a,b]$ as a metric space, the definition of SDF only takes into account the distance along the $y$-axis.
	\item Only non-periodic boundary conditions need to be specified in the computational graph.
		Periodic BCs are assumed by default, and introduce no additional graph nodes.
	\item For $\partial_{\boldsymbol{n}}=\frac{\partial}{\partial\boldsymbol{n}}$,
		i.e., derivative along the outward normal direction on the boundary of $\Omega$,
		we shall use the fact $\boldsymbol{n}(\boldsymbol{r})=\nabla\mathrm{SDF}_\Omega(\boldsymbol{r})$ for $\boldsymbol{r}\in\partial\Omega$.
		This indicates
		\[\partial_{\boldsymbol{n}}u=\boldsymbol{n}\cdot\nabla u=\nabla\mathrm{SDF}_\Omega(\boldsymbol{r})\cdot\nabla u
		=u_x\cdot\partial_x\mathrm{SDF}_\Omega(\boldsymbol{r})+u_y\cdot\partial_y\mathrm{SDF}_\Omega(\boldsymbol{r}),\]
		and the corresponding computational graph representation follows.
\end{itemize}

\subsection{Terms with unknown formulations}\label{app:dag_unknown}
\begin{figure}[htpb]
	\centering
	\includegraphics[width=0.7\linewidth]{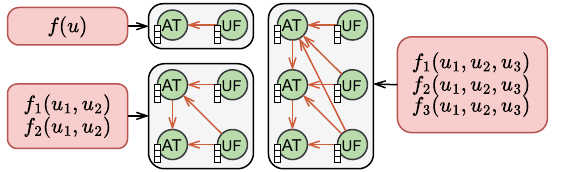}
	\caption{Computational graph examples of an arbitrary transformation $\boldsymbol{f}:\R^m\to\R^n$ that may occur in some PDEs, showing $n=1,2,3$.}%
	\label{fig:unknown_func_dag}
\end{figure}
\begin{figure}[htbp!]
	\centering
	\begin{subfigure}[b]{0.33\textwidth}
		\centering
		\includegraphics[width=\textwidth]{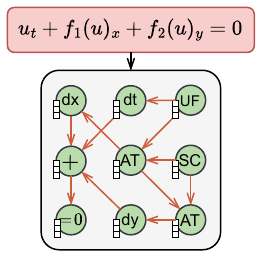}
		\caption{}
		\label{fig:unknown_func_conserv}
	\end{subfigure}%
	~ 
	\begin{subfigure}[b]{0.3\textwidth}
		\centering
		\includegraphics[width=\textwidth]{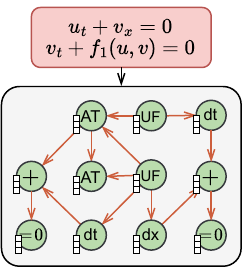}
		\caption{}
		\label{fig:unknown_func_reac}
	\end{subfigure}%
	\caption{Two computational graph examples of PDEs involving terms without known formulations.
	An extra scalar coefficient (\texttt{SC}) node is introduced to complete the computational graph in (a).}
\end{figure}

In some specific physical scenarios encountered in real applications, the underlying physical mechanism might not be completely clear, and the corresponding PDE may contain interaction terms without a known symbolic expression.
When we finetune PDEformer to adapt to such scenarios using collected sensor data, it would seem better to inform the model of these extra interactions.
Such interaction terms would typically involve an arbitrary transformation $\boldsymbol{f}:\R^m\to\R^n$ with input values $u_1,u_2,\dots,u_m$ and output values $f_1,f_2,\dots,f_n$.
This subsection interprets the conventions we made to represent this transformation using computational graphs.

We first consider the case $m=n$.
The input values $u_i$ should already have their corresponding nodes in the computational graph, and we create a node with type \verb|AT| for each output value $f_i$.
The $f_1$ node receives edges from all input nodes $u_1,\dots,u_n$.
For $i=2,3,\dots,n$, the $f_i$ node receives edges from $f_{i-1}$ and $u_i,u_{i+1},\dots,u_n$.
Such an edge convention is introduced only to indicate the non-equivalence of the input and the output nodes.
Fig.~\ref{fig:unknown_func_dag} illustrates the cases of $m=n=1,2,3$.

For the case $m<n$, we include $(n-m)$ additional scalar coefficients $c_{m+1},\dots,c_n$ to the input values,
and the augmented transformation $\bar{\boldsymbol{f}}:(u_1,\dots,u_m,c_{m+1},\dots,c_n)\mapsto(f_1,\dots,f_n)$ can be represented as discussed earlier.
The specific values of these coefficients can be any fixed real numbers, and zeros are typically used if only one such transformation is included in the PDE.
Fig.~\ref{fig:unknown_func_conserv} illustrates the computational graph of the PDE $u_t+f_1(u)_x+f_2(u)_y=0$,
which includes a transformation with $m=1$ and $n=2$ that represents the unknown flux functions in the conservation law.

For the case $m>n$, we include instead $(m-n)$ additional output values $f_{n+1},\dots,f_m$,
resulting in an augmented transformation $\bar{\boldsymbol{f}}:(u_1,\dots,u_m)\mapsto(f_1,\dots,f_m)$ represented as before.
The \verb|AT| nodes corresponding to these additional outputs will be created, but will not be used, i.e., they emit no edges towards the other nodes in the graph.
Fig.~\ref{fig:unknown_func_reac} illustrates the computational graph of the PDE $u_t+v_x=0$, $v_t+f_1(u,v)=0$,
which includes a transformation with $m=2$ and $n=1$.

Note that such a convention is only designed for the task of solving forward problems, in which we want to make predictions according to collected data.
Some scientific and engineering scenarios may also require system identification, which should recover the explicit mathematical formulation of the additional term.
If the additional term can be represented as a linear combination of a preselected set of candidate terms, this could be accomplished in a way similar to Supplementary Note~\ref{subsec:sys-id}.
The more general case is left for future work.

\subsection{Summary of node types}\label{app:dag_node_types}
Only a small number of node types are involved in the computational graph.
We list them as follows:
\begin{description}
	\item[Mathematical objects:]
		\verb|UF| (unknown field variable),
		\verb|SC| (scalar coefficient),
		\verb|CF| (coefficient field),
		\verb|VC| (varying coefficient),
		\verb|IC| (initial condition),
		\verb|SDF| (signed distance function).
	\item[Differential operations:]
		$\mathrm{dt,dx,dy}$. 
	\item[Algebraic operations:]
		$+$ (sum), $\times$ (product), $-$ (negation), $(\cdot)^2$ (square),
		$=0$ (being equal to zero).
	\item[Special functions:]
		\verb|sin|, \verb|cos| (trigonometric functions),
		\verb|exp10|, \verb|log10| (exponential and logarithm of base $10$),
		\verb|AT| (arbitrary transformation).
	\item[Auxiliary nodes:]
		$\mathtt{b}_1,\mathtt{b}_2,\dots,\mathtt{b}_N$ (branch nodes to receive input from function encoder),
		$\mathtt{m}_1,\mathtt{m}_2,\dots,\mathtt{m}_L$ (modulation nodes to generate output into INR decoder),
		\verb|pad| (padding nodes to enable batch training).
\end{description}
Note that the exponential and logarithm operations are of base $10$ instead of $e$.
We adopt such a convention to \verB{suppress} the magnitude of the input into scalar and function encoders,
since typical neural networks work better for $O(1)$ inputs and outputs.
For example, $\log_{10}10^{-3}$ has a smaller absolute value than $\log_e10^{-3}$.

\section{Datasets}\label{app:dataset}
\subsection{Pretraining}\label{app:dataset_pretrain}
To generate random PDEs that constitute our dataset, we randomly set the coefficients in a generic PDE form to obtain randomized specific forms of PDEs.
In the specific form, zero terms (including terms with zero coefficients) are removed.
If a term has coefficient one, this coefficient no longer appears in the specific form, for example ``$+1u$'' would degenerate to ``$+u$''.
Coefficient fields that have a constant value are treated as scalar coefficients.
Such a specific form simplifies the mathematical expression of the PDE as well as the corresponding computational graph.

Indeed, not all such specific PDEs admit a valid solution on the whole time interval $t\in[0,1]$.
Even if it does, the traditional solver may also fail if the solver parameters are incompatible with the current PDE.
Lacking a practical criterion to automatically identify these cases and tune solver parameters,
we choose the na\"ive method to discard the sampled PDEs that our traditional solvers failed to solve,
and only record the successful results in our datasets.
The rate of a successful solution differs across different datasets, ranging from 10\% to 90\%.

A total of $8$ generic PDE forms are included in our pretraining dataset,
\verB{and their canonical equation templates and which terms are optional are summarized in Table~\ref{tab:pretrain_pde_templates}.
We shall introduce them} in detail in the following text.

\verB{
\begin{table}[htpb]
	\centering
	\caption{Canonical equation templates for the eight PDE categories in the pretraining dataset.
		Common distributions of coefficients and source terms across all categories are listed in Table~\ref{tab:pretrain_coef_dist}.
		See the corresponding subsections for category-specific details.
	}
	\label{tab:pretrain_pde_templates}
	\small
	\renewcommand{\arraystretch}{1.35}
	\begin{tabular}{c p{6cm} p{7cm}}
		\toprule
		\textbf{Category} & \textbf{Equation template} & \textbf{Variable components} \\
		\midrule
		DCR &
		\(u_t + Lu + f_0(u) + s(\boldsymbol{r})\) \newline
		\(\quad {} + \partial_x f_1(u) + \partial_y f_2(u) = 0\) &
		$L$: diffusion (3 alternative forms);\newline
		$f_i(u)$: polynomial ($+$ up to $J$ sinusoidal terms);\newline
		BC: periodic or (generalized) Robin \\
		\midrule
		Wave &
		\(u_{tt} + \mu(\boldsymbol{r})u_t + Lu + f_0(u) + s(\boldsymbol{r})\) \newline
		\(\quad {} + \partial_x f_1(u) + \partial_y f_2(u) = 0\) &
		$\mu(\boldsymbol{r})$: damping (optional);\newline
		BC: additionally generalized Mur\\
		\midrule
		MV-DCR &
		\(\begin{aligned}[t]
			&\partial_t u_i + L_i u_i + \boldsymbol{f}_0(\boldsymbol{u})_i + s_i \\
			&\quad {} + \partial_x\boldsymbol{f}_1(\boldsymbol{u})_i + \partial_y\boldsymbol{f}_2(\boldsymbol{u})_i = 0
		\end{aligned}\) \newline
		(for $d_u\in\{2,3,4\}$ variables) &
		$\boldsymbol{f}_l(\boldsymbol{u})_i = \sum a_{lij}u_j + \sum b_{lijk}u_ju_k$ (sparse);\newline
		BC: periodic only \\
		\midrule
		DC-DCR &
		MV-DCR form ($d_u{=}2$) with pressure terms $-c_i p + \partial_{x_i}p$;
		divergence constraint: \newline
		\(\partial_x u_0+\partial_y u_1 + c_0 u_0 + c_1 u_1 + c_2 = 0\) &
		$c_0,c_1,c_2$: constraint coefficients;\newline
		includes INS as special case;\newline
		BC: periodic only \\
		\midrule
		MV-Wave &
		\(\begin{aligned}[t]
			&\partial_{tt} u_i + \mu_i\partial_t u_i + L_i u_i + \boldsymbol{f}_0(\boldsymbol{u})_i + s_i \\
			&\quad {} + \partial_x\boldsymbol{f}_1(\boldsymbol{u})_i + \partial_y\boldsymbol{f}_2(\boldsymbol{u})_i = 0
		\end{aligned}\) \newline
		(for $d_u\in\{2,3,4\}$ variables) &
		$\mu_i$: damping (optional);\newline
		otherwise same as MV-DCR;\newline
		BC: periodic only \\
		\midrule
		DC-Wave &
		MV-Wave form ($d_u{=}2$) $+$ pressure $+$ divergence constraint &
		Same as DC-DCR but with $u_{tt}$ and $\mu_i$;\newline
		BC: periodic only \\
		\midrule
		G-SWE &
		\(\begin{aligned}[t]
			h_t &+ L_0h + \boldsymbol{f}_0([h;u;v])_0 + s_0 \\
			&+ ((h+H)u)_x + ((h+H)v)_y = 0,\\
			u_t &+ L_1u + \boldsymbol{f}_0([h;u;v])_1 + s_1 \\
			&+ uu_x + vu_y + G_1 h_x = 0,\\
			v_t &+ L_2v + \boldsymbol{f}_0([h;u;v])_2 + s_2 \\
			&+ uv_x + vv_y + G_2 h_y = 0
		\end{aligned}\) &
		$H(\boldsymbol{r})$: zero, one, random scalar, or random field;\newline
		$G_1,G_2$: random positive scalars;\newline
		$L_i,s_i,\boldsymbol{f}_0$: as DCR (each optional);\newline
		BC: periodic only \\
		\midrule
		Elasticity &
		\(\begin{aligned}[t]
			&(\lambda^*(u_x+v_y)+2\mu u_x)_x \\
			&\quad {} + (\mu(u_y+v_x))_y + f_1 = 0,\\
			&(\lambda^*(u_x+v_y)+2\mu v_y)_y \\
			&\quad {} + (\mu(u_y+v_x))_x + f_2 = 0
		\end{aligned}\) &
		$\lambda^*(\boldsymbol{r}),\mu(\boldsymbol{r})$ from $E(\boldsymbol{r}),\nu(\boldsymbol{r})$ (scalar or random field);\newline
		$f_1,f_2$: zero, scalar, or random field;\newline
		BC: Dirichlet/Neumann mix per variable per edge \\
		\bottomrule
	\end{tabular}
\end{table}

\begin{table}[htpb]
	\centering
	\caption{Common distributions of coefficients, source terms, and initial conditions across all pretraining PDE categories.}
	\label{tab:pretrain_coef_dist}
	\small
	\begin{tabular}{ll}
		\toprule
		\textbf{Component} & \textbf{Distribution} \\
		\midrule
		Scalar coefficients & Zero, one, or scalar drawn from $U([-1,1])$ \\
		Initial conditions & GRFs~\eqref{eq:grf} \\
		Source $s(\boldsymbol{r})$, damping $\mu(\boldsymbol{r})$ & Zero, one, $U([-1,1])$ scalar, or GRF~\eqref{eq:grf} \\
		Scalar $a$ in $Lu$ & Scalar: $\log_{10}a\sim U([-3,-2])$ (DCR) or $U([-2,\log_{10}4])$ (Wave)\\
		Field $a(\boldsymbol{r})$ in $Lu$ & $\log_{10}a(\boldsymbol{r})=$ GRF~\eqref{eq:grf} rescaled to a random subinterval\\
		\bottomrule
	\end{tabular}
\end{table}
}

\begin{table}[htpb]
	\centering
	\caption{Summary of the pretraining datasets.
	The number of testing samples is fixed to 1k for all datasets, and not shown in the table.}
	\label{tab:pretrain_datasets_summary}
	\begin{tabular}{cccc}
		\toprule
		\textbf{PDE type} & \textbf{Dataset name} & \textbf{Training samples} & \textbf{Recorded variables} \\
		\midrule
		\multirow{6}{*}{DCR}
			& \verb|dcr_base| & 110k & $u$ \\
			& \verb|dcr_npX| & 105k & $u$ \\
			& \verb|dcr_npY| & 105k & $u$ \\
			& \verb|dcr_disk| & 200k & $u$ \\
			& \verb|dcr_sJ3| & 125k & $u$ \\
			& \verb|dcr_inhom| & 100k & $u$ \\
		\midrule
		\multirow{5}{*}{Wave}
			& \verb|wave_base| & 125k & $u$ \\
			& \verb|wave_npX| & 90k & $u$ \\
			& \verb|wave_npY| & 90k & $u$ \\
			& \verb|wave_disk| & 200k & $u$ \\
			& \verb|wave_inhom| & 100k & $u$ \\
		\midrule
		\multirow{4}{*}{MV-DCR}
			& \verb|mvdcr_2| & 160k & $u_0,u_1$ \\
			& \verb|mvdcr_2_0| & 100k & $u_0,u_1$ \\
			& \verb|mvdcr_3_1| & 200k & $u_0,u_1,u_2$ \\
			& \verb|mvdcr_4_0| & 200k & $u_0,u_1,u_2,u_3$ \\
		\midrule
		\multirow{2}{*}{DC-DCR}
			& \verb|dcdcr_icA| & 200k & $u_0,u_1,p$ \\
			& \verb|dcdcr_icV| & 200k & $u_0,u_1,p$ \\
		\midrule
		\multirow{3}{*}{MV-Wave}
			& \verb|mvwave_2| & 120k & $u_0,u_1$ \\
			& \verb|mvwave_3| & 100k & $u_0,u_1,u_2$ \\
			& \verb|mvwave_4| & 100k & $u_0,u_1,u_2,u_3$ \\
		\midrule
		\multirow{2}{*}{DC-Wave}
			& \verb|dcwave_icA| & 100k & $u_0,u_1$ \\
			& \verb|dcwave_icV| & 100k & $u_0,u_1$ \\
		\midrule
		G-SWE & \verb|swe| & 100k & $h,u,v$ \\
		\midrule
		Elasticity & \verb|elasticity| & 255k & $u,v$ \\
		\bottomrule
	\end{tabular}
\end{table}

\subsubsection{Diffusion-convection-reaction (DCR) equation}
\label{app:data_pretrain_dcr}
The PDE takes the form
\[u_t+Lu+f_0(u)+s(\boldsymbol{r})+f_1(u)_x+f_2(u)_y=0,\]
\[u(0,\boldsymbol{r})=g(\boldsymbol{r})\]
for $t\in[0,1]$, $\boldsymbol{r}=(x,y)\in\Omega\subseteq[0,1]^2$.
On each edge $\Gamma_i$, the boundary condition imposed is either periodic or of
the general form $B_iu(\boldsymbol{r})=0$ for $\boldsymbol{r}\in\Gamma_i$.

The initial value is generated as a Gaussian random field (GRF) according to
\begin{equation}\label{eq:grf}
	g(\boldsymbol{r}) \sim \mathcal{N}(0,\ 50^2(-\Delta + 5^2 I)^{-3})
.\end{equation}
The source term $s(\boldsymbol{r})$ is taken to be zero, one, a random scalar from $U([-1,1])$, or a GRF as in~\eqref{eq:grf}, with probability $5:1:5:5$.
For the convection flux and reaction term,
we take $f_i(u) = \sum_{k=1}^3c_{i0k}u^k + \sum_{j=1}^{J_i}c_{ij0}h_{ij}(c_{ij1}u+c_{ij2}u^2)$
for $i=0,1,2$, where the non-negative integers $J_0+J_1+J_2\le J$ are randomly generated.
Each scalar coefficient $c_{ijk}$ involved is taken to be zero, one, or a random scalar drawn from $U([-1,1])$.

The spatial second-order term (diffusion for DCR) $Lu$ is randomly selected from
the non-divergence form $Lu=-a(\boldsymbol{r})\Delta u$, the factored form
\verB{$Lu=-\tilde{a}(\boldsymbol{r})\nabla\cdot(\tilde{a}(\boldsymbol{r})\nabla u)$ with $\tilde{a}(\boldsymbol{r})=\sqrt{a(\boldsymbol{r})}$,}
and the divergence form $Lu=-\nabla\cdot(a(\boldsymbol{r})\nabla u)$ with equal probability.
For the homogeneous case, $\log_{10}a(\boldsymbol{r})\sim U([-3,-2])$ is a random scalar.
If the flux functions $f_1(u),f_2(u)$ are both linear, we set $a(\boldsymbol{r})$ to zero with probability $0.2$.
For the inhomogeneous case, we generate a GRF as in~\eqref{eq:grf}, rescale and shift to make its range span a random subinterval of $[-3,-1]$,
and then take exponential with base $10$ to get $a(\boldsymbol{r})$.

Each boundary operator $B_iu$ is taken to be Robin with type I
$B_iu = u + \beta_i(\boldsymbol{r})\partial u/\partial n + \gamma_i(\boldsymbol{r})$, or Robin with type II
$B_iu = \alpha_i(\boldsymbol{r})u + \partial u/\partial n + \gamma_i(\boldsymbol{r})$, with equal probability.
Each of the coefficient fields $\alpha_i(\boldsymbol{r}),\beta_i(\boldsymbol{r}),\gamma_i(\boldsymbol{r})$ is taken to be zero,
one, a random scalar from $U([-1,1])$, or a GRF in~\eqref{eq:grf}, with probability $5:1:5:5$.
Note that when $\alpha_i(\boldsymbol{r})$ or $\beta_i(\boldsymbol{r})$ equals zero, the boundary condition would
degenerate to the Dirichlet type or the Neumann type. We may also set
$\gamma_i(\boldsymbol{r})$ to be compatible with the initial condition with certain probability.

As part of our PDE, the generic boundary condition reduces to specific forms after these coefficients are selected.
If the boundary conditions share the same specific form across multiple edges $\Gamma_i$,
we merge them into a single edge, and the \verB{computational} graph is further simplified.

We have prepared several datasets for DCR equation, each containing at least 100k training samples.
The datasets are listed as follows:
\begin{itemize}
	\item \verb|dcr_base|: Basic DCR dataset, being periodic along both axes (involving no boundary operators) on $\Omega=[0,1]^2$, with $J=0$ (no sinusoidal terms).
		Diffusion is homogeneous, i.e., $a(\boldsymbol{r})$ is a constant field, degenerating to a random scalar.
	\item \verb|dcr_npX|: Non-periodic along the $x$-axis of $\Omega=[0,1]^2$, and the boundary conditions at the left and right edges need to be generated.
		(We only mark the difference from \verb|dcr_base|. This applies to the other datasets.)
	\item \verb|dcr_npY|: Non-periodic along the $y$-axis.
	\item \verb|dcr_disk|: The domain is a random disk $\Omega=D(\boldsymbol{r}_0;R)\subseteq[0,1]^2$,
		with $R\sim U([0.3,0.5])$, $\boldsymbol{r}_0\sim U([R,1-R]^2)$.
		The boundary condition on the outer boundary circle needs to be generated.
	\item \verb|dcr_sJ3|: The PDE contains up to $J=3$ sinusoidal terms.
	\item \verb|dcr_inhom|: Diffusion is inhomogeneous, and $a(\boldsymbol{r})$ has spatial dependency.
\end{itemize}

\subsubsection{Wave equation}\label{app:data_pretrain_wave}
The PDE takes the form
\[u_{tt}+\mu(\boldsymbol{r})u_t+Lu+f_0(u)+s(\boldsymbol{r})+f_1(u)_x+f_2(u)_y=0,\]
\[u(0,\boldsymbol{r})=g(\boldsymbol{r}),\quad u_t(0,\boldsymbol{r})=h(\boldsymbol{r})\]
for $t\in[0,1]$, $\boldsymbol{r}=(x,y)\in\Omega\subseteq[0,1]^2$.
On each edge $\Gamma_i$, the boundary condition imposed is either periodic or of
the general form $B_iu(\boldsymbol{r})=0$ for $\boldsymbol{r}\in\Gamma_i$.

The initial values $g(\boldsymbol{r}),h(\boldsymbol{r})$ are generated as independent GRFs given in~\eqref{eq:grf}.
The source term $s(\boldsymbol{r})$ is the same as in DCR equation,
and the damping coefficient $\mu(\boldsymbol{r})$ is generated similarly and independently.
For the homogeneous case, the squared velocity field $a(\boldsymbol{r})$ in the spatial second-order term $Lu$
is taken to be zero, one, or a random positive scalar satisfying $\log_{10}a(\boldsymbol{r})\sim U([-2,\log_{10}4])$.
The probability for the three cases is $5:2:18$ if $f_1(u),f_2(u)$ are both zeros, and $0:1:9$ otherwise.
The inhomogeneous case is similar to that of DCR equation,
in which we let $\log_{10}a(\boldsymbol{r})$ span a subinterval of $[-2,\log_{10}4]$.

Each boundary operator $B_iu$ is taken to be Robin with type I
$B_iu = u + \beta_i(\boldsymbol{r})\partial u/\partial n + \gamma_i(\boldsymbol{r})$, Robin with type II
$B_iu = \alpha_i(\boldsymbol{r})u + \partial u/\partial n + \gamma_i(\boldsymbol{r})$, or generalized Mur type
$B_iu = u_t + \alpha_i(\boldsymbol{r})u + \beta_i(\boldsymbol{r})\partial u/\partial n + \gamma_i(\boldsymbol{r})$, with equal probability.
Note that when $\alpha_i(\boldsymbol{r})=0,\beta_i(\boldsymbol{r})=c(\boldsymbol{r})$ is the wave velocity, and $\gamma_i(\boldsymbol{r})=0$,
the generalized Mur boundary condition becomes the ordinary Mur boundary condition that can absorb waves.

The remaining settings are the same as the DCR equation.

List of datasets:
\begin{itemize}
	\item \verb|wave_base|: Basic wave dataset, being periodic along both axes of $\Omega=[0,1]^2$, with $J=3$ (up to 3 sinusoidal terms).
		Wave velocity is homogeneous.
	\item \verb|wave_npX|: Non-periodic along the $x$-axis of $\Omega=[0,1]^2$.
		(We only mark the difference from \verb|wave_base|.)
	\item \verb|wave_npY|: Non-periodic along the $y$-axis.
	\item \verb|wave_disk|: The domain is a random disk $\Omega=D(\boldsymbol{r}_0;R)\subseteq[0,1]^2$ like \verb|dcr_disk|.
		We do not introduce any sinusoidal terms ($J=0$) for simplicity.
	\item \verb|wave_inhom|: Wave velocity is inhomogeneous, and has spatial dependency.
\end{itemize}

\subsubsection{Multi-variable DCR equation (MV-DCR)}\label{app:data_pretrain_mvdcr}
The PDE system includes $d_u$ equations and initial conditions, in which the one indexed by $i\in\{0,\dots,d_u-1\}$ takes the form
\[\partial_tu_i + L_iu_i + \boldsymbol{f}_0(\boldsymbol{u})_i + s_i(\boldsymbol{r}) + \partial_x\boldsymbol{f}_1(\boldsymbol{u})_i + \partial_y\boldsymbol{f}_2(\boldsymbol{u})_i = 0,\]
\[u_i(0,\boldsymbol{r})=g_i(\boldsymbol{r})\]
for $t\in[0,1]$, $\boldsymbol{r}=(x,y)\in[0,1]^2$.
Periodic boundary conditions along both axes are employed for simplicity.

For the reaction and convection flux functions, we take
\begin{equation}\label{eq:mv_interact_term}
	\boldsymbol{f}_l(\boldsymbol{u})_i = \sum_{0\le j\le d_u-1}a_{lij}u_j + \sum_{0\le j\le k\le d_u-1}b_{lijk}u_ju_k
\end{equation}
for $l=0,1,2$, $i=0,\dots,d_u-1$.
The coefficient tensors $\boldsymbol{a}_l,\boldsymbol{b}_l$ are sparse, allowing a total of at
most $3d_u$ non-zero entries.
Each non-zero entry is set to be one with probability $0.5$, and drawn from $U([-1,1])$ otherwise.

The initial value $g_i(\boldsymbol{r})$, the homogeneous spatial second-order term $L_iu_i$, and the source term $s_i(\boldsymbol{r})$ are generated in the same way as DCR equations, independently for different $i$'s.
Each homogeneous diffusion $L_iu_i$ is allowed to vanish (i.e., $a_i(\boldsymbol{r})=0$) if $b_{1iii}=b_{2iii}=0$ holds,
in which case the two flux functions affecting $u_i$ are linear with respect to $u_i$ itself.

We note that for the case of $d_u = 3$, the two-dimensional Maxwell's equations with homogeneous media and periodic boundaries are included within the distribution of this dataset,
including both Transverse Electric (TE) and Transverse Magnetic (TM) forms.
That is, for a specific choice of the random coefficients, the specific PDE form becomes the TE form of Maxwell's equations:
\[\begin{aligned}
	\partial_t E_x &&+ \partial_y((-1)H_z) &= 0, \\
	\partial_t E_y &+ \partial_x(H_z)&&=0, \\
	\partial_t H_z &+ \partial_x((-1)E_y) &+ \partial_y E_x&=0,
\end{aligned}\]
where $E_x$ and $E_y$ are the electric field components, and $H_z$ is the magnetic field component perpendicular to the 2D plane.
In this setting, we make the identification $u_0 = E_x$, $u_1 = E_y$, and $u_2 = H_z$.

Similarly, for the TM form, the equations take the specific form:
\[\begin{aligned}
	\partial_t H_x &&+ \partial_y(E_z) &= 0, \\
	\partial_t H_y &+ \partial_x((-1)E_z) &&= 0, \\
	\partial_t E_z &+ \partial_x(H_y) &+ \partial_y((-1)H_x) &= 0,
\end{aligned}\]
where $H_x$ and $H_y$ are the magnetic field components in the 2D plane, and $E_z$ is the electric field component perpendicular to the plane.
In this case, we make the identification $u_0 = H_x$, $u_1 = H_y$, and $u_2 = E_z$.

List of datasets:
\begin{itemize}
	\item \verb|mvdcr_2|: The vector $\boldsymbol{u}$ has dimension $d_u=2$.
	\item \verb|mvdcr_2_0|: Dimension $d_u=2$.
		The tensors $\boldsymbol{b}_1$ and $\boldsymbol{b}_2$ have no non-zero entry, leading to linear flux functions.
		This makes the dataset easier to learn in pretraining.
	\item \verb|mvdcr_3_1|: Dimension $d_u=3$.
		The tensor $\boldsymbol{b}_1$ has at most one non-zero entry, and same for $\boldsymbol{b}_2$,
		which increases the probability for Dedalus to successfully generate a solution.
	\item \verb|mvdcr_4_0|: Dimension $d_u=4$.
		The tensors $\boldsymbol{b}_1$ and $\boldsymbol{b}_2$ have no non-zero entry.
\end{itemize}

\subsubsection{Divergence-constrained DCR equation (DC-DCR)}
The PDE system includes $d_u=2$ equations
\[\begin{aligned}
	\partial_tu_0 + L_0u_0 + \boldsymbol{f}_0(\boldsymbol{u})_0 + s_0(\boldsymbol{r}) + \partial_x\boldsymbol{f}_1(\boldsymbol{u})_0 + \partial_y\boldsymbol{f}_2(\boldsymbol{u})_0 + (-c_0)p + \partial_xp &= 0,\\
	\partial_tu_1 + L_1u_1 + \boldsymbol{f}_0(\boldsymbol{u})_1 + s_1(\boldsymbol{r}) + \partial_x\boldsymbol{f}_1(\boldsymbol{u})_1 + \partial_y\boldsymbol{f}_2(\boldsymbol{u})_1 + (-c_1)p + \partial_yp &= 0,
\end{aligned}\]
and an additional equation
\[\partial_xu_0 + \partial_yu_1 + c_0u_0 + c_1u_1 + c_2 = 0\]
that imposes a constraint on the divergence $\nabla\cdot\boldsymbol{u}=\partial_xu_0+\partial_yu_1$.
The coefficients $c_0,c_1$, and $c_2$ are randomly and independently chosen to be zero, one, or a number from $U([-1,1])$ with equal probability.
If $c_0=c_1=0$ holds, we set $c_2=0$ to make it compatible with the periodic domain, leading to the divergence-free condition.
In this case, an additional pressure gauge condition $\int_\Omega p(t,\boldsymbol{r})\,\mathrm{d}\boldsymbol{r}=0$ is included during data preparation to ensure uniqueness of solution.
Such a pressure gauge is adopted as a convention for all similar PDEs (including INS), and not included in the input of the PDEformer-2 model.

The other settings are shared with MV-DCR equations, including the initial conditions on $u_0$ and $u_1$,
the functions $\boldsymbol{f}_l(\boldsymbol{u})$, the spatio-temporal domain, and the periodic boundaries.
Note that the initial condition of the ``pressure'' variable $p$ does not need to be specified.

We note that the two-dimensional incompressible Navier-Stokes (INS) equation in conservation form
are included within the distribution of this dataset.
That is, for a specific choice of the random coefficients, the specific PDE form becomes the INS equation:
\begin{equation}\label{eq:ins}\begin{aligned}
	\partial_tu_0-a\Delta u_0+\partial_x(u_0^2)+\partial_y(u_0u_1)+\partial_xp&=0,\\
	\partial_tu_1-a\Delta u_1+\partial_x(u_0u_1)+\partial_y(u_1^2)+\partial_yp&=0,\\
	\partial_xu_0+\partial_yu_1&=0.
\end{aligned}\end{equation}
If readers are more familiar with the convection form involving $\boldsymbol{u}\cdot\nabla\boldsymbol{u}$,
its mathematical equivalence%
\footnote{Here we only consider classical solutions of PDEs, in which case $\boldsymbol{u}(t,\boldsymbol{r})$ is sufficiently smooth.}
to the conservation form can be shown using the identity
$\nabla\cdot(\boldsymbol{uu}^\mathrm{T})=\boldsymbol{u}\cdot\nabla\boldsymbol{u}+(\nabla\cdot\boldsymbol{u})\boldsymbol{u}=\boldsymbol{u}\cdot\nabla\boldsymbol{u}$
when $\nabla\cdot\boldsymbol{u}=0$ holds.

List of datasets:
\begin{itemize}
	\item \verb|dcdcr_icA|: Basic DC-DCR dataset.
		The initial values $g_0(\boldsymbol{r}),g_1(\boldsymbol{r})$ are generated as independent GRFs~\eqref{eq:grf},
		and are not guaranteed to satisfy the divergence constraint.
		The \textbf{i}nitial \textbf{c}ondition is thus called \textbf{a}rbitrary (icA) in this sense.
		The Dedalus solver we utilize will make the solution comply with the divergence constraint automatically.
	\item \verb|dcdcr_icV|:
		The \textbf{i}nitial \textbf{c}ondition is \textbf{v}alid (icV), being compatible with the divergence constraint.
		To achieve this, we first draw a GRF $\psi(\boldsymbol{r})$ according to~\eqref{eq:grf},
		and set $\tilde g_0(\boldsymbol{r})=c_1\psi+\psi_y$, $\tilde g_1(\boldsymbol{r})=-(c_0\psi+\psi_x)$,
		in which the \verB{differentiation} is computed using finite-difference.
		We then adjust the mean and variance randomly to obtain the ultimate initial values,
		and reset $c_2$ to ensure \verB{compatibility} with them.
\end{itemize}

\subsubsection{Multi-variable wave equation (MV-Wave)}
Similar to MV-DCR,
the PDE system includes $d_u$ equations and initial conditions, in which the one indexed by $i\in\{0,\dots,d_u-1\}$ takes the form
\[\partial_{tt}u_i + \mu_i(\boldsymbol{r})\partial_tu_i + L_iu_i + \boldsymbol{f}_0(\boldsymbol{u})_i + s_i(\boldsymbol{r}) + \partial_x\boldsymbol{f}_1(\boldsymbol{u})_i + \partial_y\boldsymbol{f}_2(\boldsymbol{u})_i = 0,\]
\[u_i(0,\boldsymbol{r})=g_i(\boldsymbol{r}),\quad\partial_tu_i(0,\boldsymbol{r})=h_i(\boldsymbol{r})\]
for $t\in[0,1]$, $\boldsymbol{r}=(x,y)\in[0,1]^2$.
Each homogeneous spatial second-order term $L_iu_i$ is allowed to vanish (i.e., $a_i(\boldsymbol{r})=0$) if $a_{1ii}=a_{2ii}=b_{1iik}=b_{2iik}=b_{1iki}=b_{2iki}=0$ holds for all $k$'s,
in which case the two flux functions affecting $u_i$ do not depend on $u_i$ itself.
The coefficients are generated as in the Wave equation described in Supplementary Note~\ref{app:data_pretrain_wave}.

List of datasets:
\begin{itemize}
	\item \verb|mvwave_2|: The vector $\boldsymbol{u}$ has dimension $d_u=2$.
	\item \verb|mvwave_3|: Dimension $d_u=3$.
	\item \verb|mvwave_4|: Dimension $d_u=4$.
\end{itemize}

\subsubsection{Divergence-constrained wave equation (DC-Wave)}
Similar to DC-DCR,
the PDE system includes $d_u=2$ equations
\[\begin{aligned}
	\partial_{tt}u_0 + \mu_0(\boldsymbol{r})\partial_tu_0 + L_0u_0 + \boldsymbol{f}_0(\boldsymbol{u})_0 + s_0(\boldsymbol{r}) + \partial_x\boldsymbol{f}_1(\boldsymbol{u})_0 + \partial_y\boldsymbol{f}_2(\boldsymbol{u})_0 + (-c_0)p + \partial_xp &= 0,\\
	\partial_{tt}u_1 + \mu_1(\boldsymbol{r})\partial_tu_1 + L_1u_1 + \boldsymbol{f}_0(\boldsymbol{u})_1 + s_1(\boldsymbol{r}) + \partial_x\boldsymbol{f}_1(\boldsymbol{u})_1 + \partial_y\boldsymbol{f}_2(\boldsymbol{u})_1 + (-c_1)p + \partial_yp &= 0,
\end{aligned}\]
and an additional divergence constraint condition
\[\partial_xu_0 + \partial_yu_1 + c_0u_0 + c_1u_1 + c_2 = 0.\]
The other settings are shared with MV-Wave.
As the pressure value often exceeds the range $[-10,10]$, we do not record it in the final dataset as in DC-DCR, so as to make learning easier for our model.

List of datasets:
\begin{itemize}
	\item \verb|dcwave_icA|: Basic DC-Wave dataset like \verb|dcdcr_icA|.
		The initial values $g_0(\boldsymbol{r}),g_1(\boldsymbol{r}),h_0(\boldsymbol{r}),h_1(\boldsymbol{r})$ are generated as independent GRFs~\eqref{eq:grf},
		and are not guaranteed to satisfy the divergence constraint.
	\item \verb|dcwave_icV|:
		The initial values $g_0(\boldsymbol{r}),g_1(\boldsymbol{r})$ as well as $h_0(\boldsymbol{r}),h_1(\boldsymbol{r})$ are generated as in \verb|dcdcr_icV|
		to be compatible with the divergence constraint.
\end{itemize}

\subsubsection{Generalized shallow-water equation (G-SWE)}
The PDE takes the form
\[\begin{aligned}
	h_t + L_0h + \boldsymbol{f}_0([h;u;v])_0 + s_0(\boldsymbol{r}) &+ ((h+H(\boldsymbol{r}))u)_x + ((h+H(\boldsymbol{r}))v)_y &= 0,\\
	u_t + L_1u + \boldsymbol{f}_0([h;u;v])_1 + s_1(\boldsymbol{r}) &+ uu_x + vu_y + G_1h_x &= 0,\\
	v_t + L_2v + \boldsymbol{f}_0([h;u;v])_2 + s_2(\boldsymbol{r}) &+ uv_x + vv_y + G_2h_y &= 0,
\end{aligned}\]
\[h(0,\boldsymbol{r})=g_h(\boldsymbol{r}),\ u(0,\boldsymbol{r})=g_u(\boldsymbol{r}),\ v(0,\boldsymbol{r})=g_v(\boldsymbol{r})\]
for $t\in[0,1]$, $\boldsymbol{r}=(x,y)\in[0,1]^2$.
Periodic boundary conditions are employed for simplicity.

The initial water height $g_h(\boldsymbol{r})$ is taken to be a positive random field
like the inhomogeneous $a(\boldsymbol{r})$ of DCR equation, with lower-bound $0.1$ and upper-bound $4$.
The initial velocity $g_u(\boldsymbol{r}),g_v(\boldsymbol{r})$ are independent GRFs according to~\eqref{eq:grf}.
The homogeneous diffusion operator $L_i$ (allowed to vanish) and the source term $s_i(\boldsymbol{r})$ are generated in the same way as DCR equations, independently for different $i$'s.
The interaction term $\boldsymbol{f}_0$ is given in~\eqref{eq:mv_interact_term},
which also appears in the multi-variable equations.

The base water height $H(\boldsymbol{r})$ can be (1) zero, (2) one,
(3) a random positive scalar satisfying $\log_{10}H(\boldsymbol{r})\sim U([-3,\log_{10}4])$.
or (4) a positive random field with lower-bound $10^{-3}$ and upper-bound $4$.
The probability for the four cases is $5:1:1:5$.
For $i=1,2$, $G_i$ is a random positive scalar that equals one with probability $0.1$, and sampled according to $\log_{10}G_i\sim U([-1,0])$ otherwise.

We note that for a specific choice of the random coefficients, the specific PDE form becomes the regular shallow-water equations:
\[\begin{aligned}
	h_t + ((h+H(\boldsymbol{r}))u)_x + ((h+H(\boldsymbol{r}))v)_y &= 0,\\
	u_t - \nu\Delta u + ku + (-F)v + uu_x + vu_y + Gh_x &= 0,\\
	v_t - \nu\Delta v + kv + Fu + uv_x + vv_y + Gh_y &= 0,
\end{aligned}\]
with kinematic viscosity $\nu$, gravitational acceleration $G$, viscous drag coefficient $k$, and Coriolis force coefficient $F$.

The only dataset is named \verb|swe|.

\subsubsection{Steady-state elasticity equation}
The PDE takes the form
\[\begin{aligned}
	(\lambda^*(\boldsymbol{r})(u_x+v_y)+2\mu(\boldsymbol{r})u_x)_x+\left(\mu(\boldsymbol{r})(u_y+v_x)\right)_y+f_1(\boldsymbol{r})&=0,\\
	\left(\mu(\boldsymbol{r})(u_y+v_x)\right)_x+(\lambda^*(\boldsymbol{r})(u_x+v_y)+2\mu(\boldsymbol{r})v_y)_y+f_2(\boldsymbol{r})&=0
\end{aligned}\]
for $\boldsymbol{r}=(x,y)\in[0,1]^2$.
Boundary condition types are randomly chosen from Dirichlet and Neumann,
independently for the two variables $u,v$ and the four edges of the square domain.
The boundary values are randomly chosen from zero, a random scalar, and a one-dimensional random field.

The \verB{plane} stress assumption is employed to reduce the 3D elasticity equation to a 2D PDE. Here, $\lambda^*(\boldsymbol{r})$
and $\mu(\boldsymbol{r})$ are the Lam\'e parameters in the 2D plane stress case. These Lam\'e parameters are related to the
Young's modulus $E$ and the Poisson's ratio $\nu$ by
\[\begin{aligned}
	\lambda^*(\boldsymbol{r}) &= \frac{E(\boldsymbol{r})\nu(\boldsymbol{r})}{(1+\nu(\boldsymbol{r}))(1-\nu(\boldsymbol{r}))},\\
	\mu(\boldsymbol{r}) &= \frac{E(\boldsymbol{r})}{2(1+\nu(\boldsymbol{r}))}.
\end{aligned}\]
Each of Young's modulus $E(\boldsymbol{r})$ and Poisson's ratio $\nu(\boldsymbol{r})$ is randomly chosen from a random positive
and a 2D random field. Here, the Poisson's ratio $\nu(\boldsymbol{r})$ is constrained to be in the range $[0.01,1.0]$.

The external force $f_i(\boldsymbol{r})$ is randomly chosen from zero, a random constant, and a 2D random field.

The only dataset is named \verb|elasticsteady|.

\subsection{Specific PDEs}\label{app:dataset_forward}
PDEformer-2 is finetuned to solve several specific PDEs as described in `Adaptation to specific PDEs' in the main text.
A detailed description of the datasets will be presented in this subsection,
including three customized datasets generated by Dedalus~\cite{Dedalus} and two publicly available datasets.

\begin{table}[htpb]
	\centering
	\caption{Summary of the downstream specific PDE datasets,
		including (1) dataset names, (2) data source,
		(3) number of time steps used in training (i.e., $n_t$),
		(4) input coefficient scalars that may differ across samples (``/'' means no such scalars),
		(5) input coefficient fields,
		(6) output variables for the model to make predictions,
		and (7) boundary conditions.
	The number of testing samples is fixed to $100$ for all datasets.}
	\label{tab:finetune_datasets_summary}
	\begin{tabular}{cccccc}
		\toprule
		\textbf{Name} & Sine-Gordon & INS-Tracer & INS-Pipe & Wave-Gauss & Wave-C-Sines \\
		\midrule
		\textbf{Source} & By Dedalus & By Dedalus & By Dedalus & \PDEgym~\cite{Poseidon} & \cite{RIGNO} \\
		\textbf{T-steps} & 100 & 100 & 100 & 16 & 20 \\
		\textbf{Scalars} & $\nu$ & / & $\nu$ & / & / \\
		\textbf{Fields} & $g,h$ & $g_u,g_v,g_s$ & $g_u,g_v,s_0,s_1$ & $g,c^2,c$ & $g$ \\
		\textbf{Outputs} & $u$ & $u,v,s$ & $u,v$ & $u$ & $u$ \\
		\textbf{BCs} & Periodic & Periodic &
			\makecell{No-slip at\\top \& bottom} &
			\makecell{Absorbing at\\four square edges} &
			\makecell{Dirichlet at\\boundary circle} \\
		\bottomrule
	\end{tabular}
\end{table}

\subsubsection{Sine-Gordon equation}
The PDE takes the form
\[u_{tt}-a\Delta u+(-1)\sin(u)=0,\]
\[u(0,\boldsymbol{r})=g(\boldsymbol{r}),\quad u_t(0,\boldsymbol{r})=h(\boldsymbol{r})\]
for $t\in[0,1]$, $\boldsymbol{r}=(x,y)\in[0,1]^2$.
Periodic boundary conditions along both axes are employed.
The initial values $g(\boldsymbol{r}),h(\boldsymbol{r})$ are generated as independent GRFs given in~\eqref{eq:grf}.
The squared velocity $a=c^2$ satisfies $\log_{10}a\sim U([-2,\log_{10}4])$.
We note that this dataset lies within the distribution of the Wave dataset in pretraining (to be more specific, \verb|wave_base|).

\subsubsection{Incompressible Navier-Stokes equation with tracer (INS-Tracer)}\label{app:data_forward_ins_tracer}
The PDE takes the form
\[\begin{aligned}
	u_t-\nu\Delta u+(u^2)_x+(uv)_y+p_x&=0,\\
	v_t-\nu\Delta v+(uv)_x+(v^2)_y+p_y&=0,\\
	s_t-D\Delta s+(us)_x+(vs)_y&=0,\\
	u_x+v_y&=0,
\end{aligned}\]
\[u(0,\boldsymbol{r})=g_u(\boldsymbol{r}), v(0,\boldsymbol{r})=g_v(\boldsymbol{r}), s(0,\boldsymbol{r})=g_s(\boldsymbol{r})\]
for $t\in[0,1]$, $\boldsymbol{r}=(x,y)\in[0,1]^2$.
Periodic boundary conditions along both axes are employed.
The initial values $g_u(\boldsymbol{r}),g_v(\boldsymbol{r})$, and $g_s(\boldsymbol{r})$ are generated as independent GRFs given in~\eqref{eq:grf}.
The fluid viscosity is fixed as $\nu=10^{-3}$, and the particle diffusivity is fixed as $D=10^{-2}$.
Like DC-Wave equations, we do not record the pressure variable $p$ in the dataset.

Compared with the INS equation~\eqref{eq:ins}, an additional variable $s$ representing the tracer particle density is included.
This specific PDE therefore lies outside the distribution of the pretraining dataset.
However, the interaction between the tracer density $s$ and the fluid velocity $u,v$ is very similar to those appeared in the MV-DCR equations in Supplementary Note~\ref{app:data_pretrain_mvdcr},
and the PDE can thus be viewed as a new combination of physical mechanisms encountered during pretraining.

\paragraph{Input PDE forms in ablation}
In the ablation study shown in Fig.~\ref{figs:ablation}(c) in the main text, we use the INS-Tracer data samples to finetune PDEformer-2-fast, but feed PDEs of different forms into the model.
For the \emph{equivalent PDE} case, the model input is the convection form of the INS-Tracer equation:
\[\begin{aligned}
	u_t-\nu\Delta u+uu_x+vu_y+p_x&=0,\\
	v_t-\nu\Delta v+uv_x+vv_y+p_y&=0,\\
	s_t-D\Delta s+us_x+vs_y&=0,\\
	u_x+v_y&=0,
\end{aligned}\]
\[u(0,\boldsymbol{r})=g_u(\boldsymbol{r}), v(0,\boldsymbol{r})=g_v(\boldsymbol{r}), s(0,\boldsymbol{r})=g_s(\boldsymbol{r}).\]
For the \emph{incomplete PDE} case, the model input takes the form
\[\begin{aligned}
	u_t+(u^2)_x+(uv)_y&=0,\\
	v_t+(uv)_x+(v^2)_y&=0,\\
	s_t+(us)_x+(vs)_y&=0,
\end{aligned}\]
\[u(0,\boldsymbol{r})=g_u(\boldsymbol{r}), v(0,\boldsymbol{r})=g_v(\boldsymbol{r}), s(0,\boldsymbol{r})=g_s(\boldsymbol{r}),\]
in which the viscosity and particle diffusivity terms, the pressure variable, and the divergence-free condition are removed.
For the \emph{unknown PDE} case, the model input takes the form
\[\begin{aligned}
	u_t+f_1(u,v,s)&=0,\\
	v_t+f_2(u,v,s)&=0,\\
	s_t+f_3(u,v,s)&=0,
\end{aligned}\]
\[u(0,\boldsymbol{r})=g_u(\boldsymbol{r}), v(0,\boldsymbol{r})=g_v(\boldsymbol{r}), s(0,\boldsymbol{r})=g_s(\boldsymbol{r}),\]
in which the functions $f_1,f_2,f_3$ are assumed not to have known mathematical formulations, and the corresponding computational graph is constructed according to Supplementary Note~\ref{app:dag_unknown}.
\verC{We note that the related finetuning experiment also tests the model's ability to adapt to a node type genuinely unseen during pretraining, as it involves \texttt{AT} nodes that are introduced only at finetuning with randomly initialized embeddings.}
For the \emph{wrong PDE} case, the model input is the Lorenz system
\[\begin{aligned}
	u_t+\sigma v+(-\sigma)u&=0,\\
	v_t+\rho u+(-1)v+(-1)us&=0,\\
	s_t+uv+(-\beta)s&=0,
\end{aligned}\]
\[u(0,\boldsymbol{r})=g_u(\boldsymbol{r}), v(0,\boldsymbol{r})=g_v(\boldsymbol{r}),\ s(0,\boldsymbol{r})=g_s(\boldsymbol{r}),\]
with coefficients $\sigma=10$, $\rho=8/3$, $\beta=0.8$.
Involving no spatial \verB{differentiations}, this system actually degenerates to independent ODE systems for different spatial points.
Note that all these PDEs are merely inputs of the PDEformer-2 model, and we do not need to obtain their numerical solutions.

\subsubsection{Incompressible Navier-Stokes equation in a pipe (INS-Pipe)}
The PDE takes the form
\[\begin{aligned}
	u_t-\nu\Delta u+s_0(\boldsymbol{r})+(u^2)_x+(uv)_y+p_x&=0,\\
	v_t-\nu\Delta v+s_1(\boldsymbol{r})+(uv)_x+(v^2)_y+p_y&=0,\\
	u_x+v_y&=0,
\end{aligned}\]
\[u(0,\boldsymbol{r})=g_u(\boldsymbol{r}),\ v(0,\boldsymbol{r})=g_v(\boldsymbol{r}).\]
The initial values $g_u(\boldsymbol{r}),g_v(\boldsymbol{r})$ and the source terms $s_0(\boldsymbol{r}),s_1(\boldsymbol{r})$ are generated as independent GRFs given in~\eqref{eq:grf}.
The fluid viscosity $\nu$ is a random scalar satisfying $\log_{10}\nu\sim U([-3,-2])$.
We do not record the pressure variable $p$ in the dataset.

Different from the INS equation~\eqref{eq:ins},
the domain $\Omega=[0,1]^2$ is a periodic pipe along the $x$-axis,
and we impose no-slip boundary conditions $u|_\Gamma=v|_\Gamma=0$ on the boundary walls $\Gamma=[0,1]\times\{0,1\}$.
This specific PDE therefore lies outside the distribution of the pretraining dataset.
As homogeneous Dirichlet boundary conditions can appear in the single-variable DCR and wave equations,
this can be viewed as a new combination of the interior PDE and boundary conditions.

\subsubsection{Wave-Gauss}
The \emph{Wave-Gauss} dataset is part of the publicly available \PDEgym~\cite{Poseidon} dataset series,
and we give a brief description of it for convenience to the readers.
The PDE is a wave equation of the form
\[u_{tt}-c(\boldsymbol{r})^2\Delta u=0,\]
\[u(0,\boldsymbol{r})=g(\boldsymbol{r}),\quad u_t(0,\boldsymbol{r})=0\]
for $\boldsymbol{r}=(x,y)\in\Omega=[0,1]^2$.
The original time interval is different from $[0,1]$ used in the other datasets.
We rescale the time variable $t$ to resolve this, and the wave velocity values $c(\boldsymbol{r})$ are scaled accordingly.

The domain $\Omega$ has absorbing boundaries.
To inform PDEformer-2 of this, we add the Mur boundary condition
\[(u_t+c(\boldsymbol{r})\partial_nu)|_{\partial\Omega}=0,\]
into the PDE input for the model during finetuning.
We note that all evolutionary equations on square domains in pretraining are periodic in at least one axis, but Wave-Gauss absorbs waves in its four edges.

Both the initial value $g(\boldsymbol{r})$ and the wave velocity field $c(\boldsymbol{r})$ are generated as linear combinations of Gaussian densities of the form
$\exp(-\|\boldsymbol{r}-\boldsymbol{r}_i\|_2^2/2\sigma_i^2)$.
This is also very different from the GRF and its variants employed in the pretraining dataset.

The solutions are recorded on $16$ equispaced timesteps (including the initial value at $t=0$) within the temporal domain.
This differs from the $101$ timesteps used in the previous datasets we generated ourselves.
Following the convention of the authors, we train the models with normalized data, making the entire dataset have zero mean and unit variance.
The form of the PDE input does not need to be adapted accordingly, since the wave equation is linear, and the linearly transformed solutions still solve the PDE.

\subsubsection{Wave-C-Sines}
The \emph{Wave-C-Sines} dataset released by~\cite{RIGNO} considers the wave equation on a disk with radius 1.
The wave velocity is $2$, and the system evolves to the terminal time $0.2$.
After a spatio-temporal rescaling operation to fit our domain, the PDE takes the form
\[u_{tt}-0.01\Delta u=0,\]
\[u(0,\boldsymbol{r})=g(\boldsymbol{r}),\quad u_t(0,\boldsymbol{r})=0\]
for $t\in[0,1]$, $\boldsymbol{r}=(x,y)\in\Omega=D((0.5,0.5);0.5)$,
with Dirichlet boundary condition $u|_{\partial\Omega}=0$.
The initial condition $g(\boldsymbol{r})$ is generated as a superposition of sinusoidal terms, with coefficients decaying slowly as the frequency increases.
This makes it more oscillatory than the GRFs in our pretraining dataset.

The solutions are recorded on $16,431$ uniformly sampled points inside the disk, which is different from the polar-coordinate-based grids in the pretraining disk datasets illustrated in Fig.~\ref{fig:dedalus_grids_disk}.
To fit the CNN-based function encoder of PDEformer-2 as well as the branch network of DeepONet (CNN), we interpolate the initial condition to the $128\times 128$ uniform grid before feeding it into the network.
Solution interpolation is avoided, since to ensure a fair error metric against the ground-truth labels, we may need to interpolate the Cartesian predictions back to scattered points, which consumes computational resources.
Ignoring the initial frame, a total of $20$ equispaced time steps are predicted by the models.
As the solutions in the original dataset have a small magnitude, they are multiplied by $100$ to better fit the neural networks, and the linear PDE preserves its form as in Wave-Gauss.

\subsection{Inverse problems}
\subsubsection{Inverse problems on DCR equations}
\label{app:inv_dcr}
In inverse problems, including recovering scalar coefficients, system identifications, and recovering source fields, the PDEs have the same form as shown in Supplementary Note~\ref{app:data_pretrain_dcr} despite some differences in the coefficients.
The source term $s(\boldsymbol{r})$ is set to a GRF (scalar is no longer allowed),
the diffusivity coefficient $a(\boldsymbol{r})$ is set to a scalar with distribution $\log_{10}a(\boldsymbol{r})\sim U([-3,-2])$ (i.e., homogeneous case),
and the parameter $J$ is set to zero, meaning that no sinusoidal term is included in the PDE.
For convenience, we apply periodic boundary conditions along both axes.
For each PDE, we fix the PDE form and parameters, and solve it with 25 different initial conditions.
The observations of these solutions are used for inverse problems.
The initial conditions are generated as independent GRFs given in~\eqref{eq:grf}.

\subsubsection{Inverse problems on wave equations}
\label{app:inv_wave}
The PDEs follow the same form as shown in Supplementary Note~\ref{app:data_pretrain_wave} despite some differences in the coefficients.
The source term $s(\boldsymbol{r})$ is set to a GRF (scalars no longer allowed), the squared velocity field $a(\boldsymbol{r})$ is set as in the inhomogeneous case,
and the parameter $J$ is set to zero, meaning that there is no sinusoidal term in the PDE.
For convenience, we apply periodic boundary conditions along both axes.
For each equation, we fix the PDE structure, scalar coefficients, and the damping term $\mu(\boldsymbol{r})$, and solve it under 100 different initial conditions and source terms.
The resulting solution observations are then used to formulate inverse problems.
The initial conditions are generated as independent Gaussian random fields (GRFs) as described in~\eqref{eq:grf}.

\subsection{Classical solvers and data grids}
Wave-Gauss and Wave-C-Sines are publicly available datasets.
We generate the steady-state elasticity dataset using the FEniCSx~\cite{FEniCSx} package, a classical solver based on finite-element method.
All the remaining datasets mentioned in this section are generated by Dedalus V3~\cite{Dedalus} based on spectral methods.
Unless specified otherwise, the spatial solution domain is discretized into $128\times 128$ grid points.

\begin{figure}[htbp!]
	\centering
	\begin{subfigure}[p]{0.3\textwidth}
		\centering
		\includegraphics[width=\textwidth]{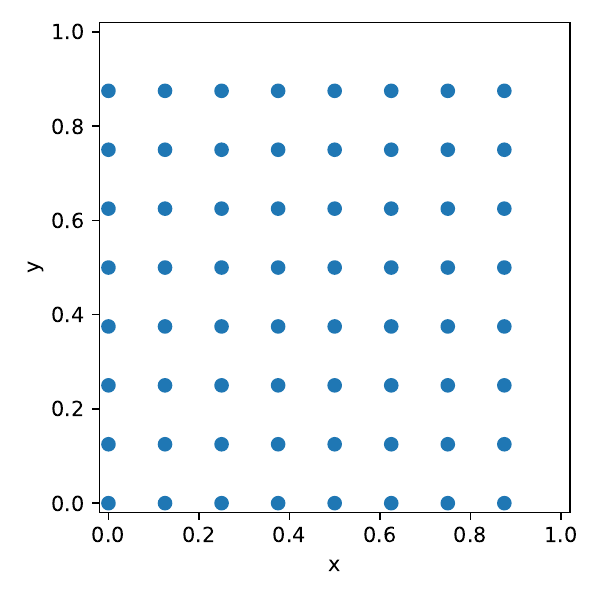}
		\caption{}
	\end{subfigure}%
	~ 
	\begin{subfigure}[p]{0.3\textwidth}
		\centering
		\includegraphics[width=\textwidth]{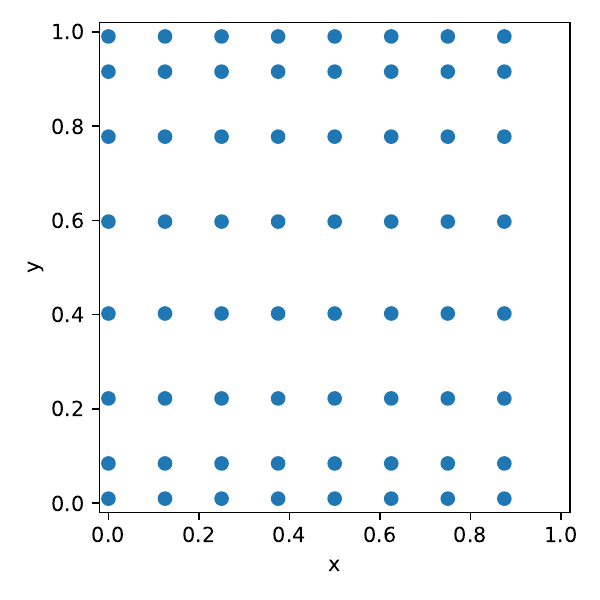}
		\caption{}
	\end{subfigure}%
	~ 
	\begin{subfigure}[p]{0.3\textwidth}
		\centering
		\includegraphics[width=\textwidth]{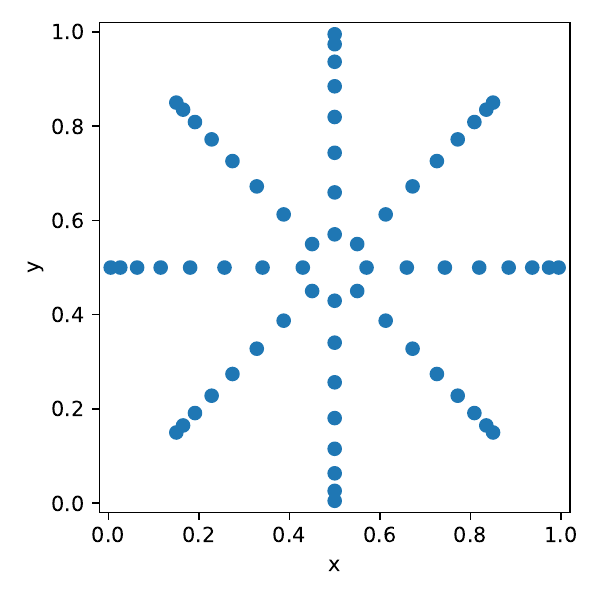}
		\caption{}%
		\label{fig:dedalus_grids_disk}
	\end{subfigure}%
	~ 
	\caption{Dedalus solution grids for different computational domains $\Omega$.
		For illustrative purposes, the resolution is $8\times 8$ instead of $128\times 128$ in the real datasets.
		(a) $\Omega=[0,1]^2$ is periodic along both axes.
		(b) $\Omega=[0,1]^2$ is periodic along the $x$-axis.
		(c) $\Omega=D((0.5,0.5);0.5)$ is a disk, in which the solution grids are based on polar coordinates.
	}%
	\label{figs:dedalus_grids}
\end{figure}
As illustrated in Fig.~\ref{figs:dedalus_grids},
Dedalus uses different solution grids for different computational domains $\Omega$.
Based on spectral methods, Dedalus employs a uniform grid along the periodic axis.
For the non-periodic axis, the grid points cluster quadratically towards the endpoints of the spatial interval.
The FEniCSx solver uses a uniform $128\times 128$ grid with triangular elements
for the steady-state elasticity equation.

Since PDEformer-2 employs an INR to decode mesh-free solutions,
we can utilize solution samples recorded on different grid points directly during training,
and no interpolation is required.
Indeed, the input coefficient fields need to be given on uniform grids to fit the CNN function encoder.
These fields, based on GRFs, are directly generated and stored on the uniform grids in our datasets,
and interpolation is only required in the data preparation stage to specify the discretized initial values.
No interpolation is applied in the training stage.

Snapshots of time-dependent PDE solutions are recorded with a time stepsize $\delta t_\text{data}=0.01$, yielding a total of $101$ time steps for each data sample.
For the three specific PDE datasets generated by Dedalus, we ignore the initial time-step ($t=0$) in the solution label during model training,
since the MindSpore implementation of FNO3D requires an even grid resolution along each axis.
For some datasets, Dedalus may need to proceed with a smaller time stepsize $\delta t_\text{solver}$ to ensure numerical stability.
Lacking an effective selection principle across diverse pretraining PDEs, we choose the numerical time stepsize by trial-and-error.
More specifically, we solve a PDE instance with a trial stepsize, and if Dedalus fails to solve, we make another trial with the stepsize divided by four, or directly terminate if the stepsize is already too small.

\verB{
\subsection{Data filtering}\label{sec:data_filtering}

As mentioned in `Pretrained results' in the main text,
during pretraining data generation, we observe that Dedalus sometimes produces solutions containing \texttt{NaN} values for certain PDEs.
Even when the solutions consist of normal floating-point numbers,
non-physical artifacts can still be identified through visual inspection.

One cause of these issues is that some randomly generated PDEs do not admit a valid solution
itself in the theoretical aspect.
For example, consider the PDE $u_t - a\Delta u +u^2 = 0$ with spatially constant initial value $u(0,x,y)=-2$.
The unique solution is $u(t,x,y) = 2/(2t-1)$,
which is only defined for $t\in[0,1/2)$ and blows up at $t=1/2$.
No solution exists over the full time interval $t\in[0,1]$.

Another cause, in the numerical aspect, is that traditional solvers often require problem-specific tuning.
A careful case-by-case selection of the solver parameters may lead to a valid solution without \texttt{NaN} values
and eliminate the non-physical patterns,
but is intractable in our pretraining data generation since the randomly sampled PDEs differ greatly from each other.

We attempted to filter anomalous data by computing the PDE residual via finite differences.
However, this na\"ive criterion proved ineffective:
we observed that the residual can also be very large on data containing clear shock waves,
exceeding the residual values on some samples that contain visible artifacts.
This makes it impossible to set a suitable threshold that filters out non-physical artifacts
without accidentally discarding valid samples with shocks.
We speculate that finite-difference-based residual computation is unsuitable for data
involving strong discontinuous shocks.
Designing a more reliable artifact screening criterion is left as a meaningful future direction.

As a practical measure, we ultimately filter out only those samples where
the Dedalus solution contains \texttt{NaN} values, or where any recorded variable $u$
satisfies $\|u\|_{L^\infty} > 10$.
All remaining data, including those that may contain non-physical artifacts,
are retained and used equally for pretraining.

\begin{table}[htpb]
	\centering
	\caption{Statistics of the data filtering process.
		Note that the three specific PDE datasets (INS-Tracer, INS-Pipe, Sine-Gordon) correspond to the downstream tasks described in Supplementary Note~\ref{app:dataset_forward},
	and the filtering conditions were rarely triggered for them.}
	\label{tab:datagen_statistics}
	\begin{tabular}{lrrrrr}
		\toprule
		\textbf{Dataset} & \textbf{Accept (\%)} & \textbf{\texttt{NaN} (\%)} & \textbf{$\|u\|_\infty>10$ (\%)} & \textbf{T/trial (s)} & \textbf{T/sample (s)} \\
		\midrule
		\multicolumn{6}{l}{\textit{DCR}} \\
		\texttt{dcr\_base} & 48.2 & 51.8 & 0.0 & 23.9 & 49.5 \\
		\texttt{dcr\_npX} & 26.9 & 73.0 & 0.1 & 22.9 & 85.1 \\
		\texttt{dcr\_npY} & 27.7 & 72.1 & 0.2 & 22.8 & 82.4 \\
		\texttt{dcr\_disk} & 34.9 & 64.9 & 0.2 & 153.3 & 439.7 \\
		\texttt{dcr\_sJ3} & 45.5 & 54.3 & 0.1 & 26.8 & 58.8 \\
		\texttt{dcr\_inhom} & 38.5 & 61.2 & 0.3 & 26.6 & 69.2 \\
		\midrule
		\multicolumn{6}{l}{\textit{Wave}} \\
		\texttt{wave\_base} & 63.4 & 32.3 & 4.3 & 7.6 & 12.0 \\
		\texttt{wave\_npX} & 32.1 & 65.1 & 2.9 & 8.4 & 26.2 \\
		\texttt{wave\_npY} & 31.9 & 65.1 & 3.0 & 8.2 & 25.7 \\
		\texttt{wave\_disk} & 45.1 & 51.4 & 3.6 & 41.3 & 91.6 \\
		\texttt{wave\_inhom} & 26.0 & 71.2 & 2.8 & 6.3 & 24.4 \\
		\midrule
		\multicolumn{6}{l}{\textit{MV-DCR}} \\
		\texttt{mvdcr\_2} & 25.6 & 72.1 & 2.3 & 37.8 & 147.4 \\
		\texttt{mvdcr\_2\_0} & 59.0 & 28.3 & 12.6 & 21.4 & 36.2 \\
		\texttt{mvdcr\_3\_0} & 40.0 & 41.6 & 18.4 & 34.5 & 86.2 \\
		\texttt{mvdcr\_3\_1} & 17.9 & 74.9 & 7.1 & 44.1 & 245.7 \\
		\texttt{mvdcr\_4\_0} & 30.1 & 52.0 & 17.8 & 45.5 & 151.0 \\
		\texttt{mvdcr\_5\_0} & 21.6 & 59.4 & 19.0 & 59.1 & 273.5 \\
		\midrule
		\multicolumn{6}{l}{\textit{DC-DCR}} \\
		\texttt{dcdcr\_icA} & 52.7 & 11.5 & 35.7 & 54.5 & 103.3 \\
		\texttt{dcdcr\_icV} & 68.2 & 24.8 & 7.0 & 71.1 & 104.2 \\
		\midrule
		\multicolumn{6}{l}{\textit{MV-Wave}} \\
		\texttt{mvwave\_2} & 82.9 & 7.1 & 10.0 & 9.7 & 11.7 \\
		\texttt{mvwave\_3} & 84.6 & 7.4 & 8.0 & 13.6 & 16.1 \\
		\texttt{mvwave\_4} & 84.8 & 5.2 & 10.0 & 17.4 & 20.6 \\
		\midrule
		\multicolumn{6}{l}{\textit{DC-Wave}} \\
		\texttt{dcwave\_icA} & 93.5 & 2.9 & 3.6 & 13.0 & 13.9 \\
		\texttt{dcwave\_icV} & 75.6 & 16.1 & 8.2 & 12.7 & 16.8 \\
		\midrule
		\multicolumn{6}{l}{\textit{G-SWE}} \\
		\texttt{swe} & 59.8 & 38.2 & 2.0 & 99.6 & 166.5 \\
		\midrule
		\multicolumn{6}{l}{\textit{Specific PDEs}} \\
		\texttt{Sine-Gordon} & 100.0 & 0.0 & 0.0 & 5.4 & 5.4 \\
		\texttt{INS-Tracer} & 99.8 & 0.0 & 0.2 & 61.6 & 61.7 \\
		\texttt{INS-Pipe} & 98.7 & 0.0 & 1.3 & 35.9 & 36.4 \\
		\midrule
		\textbf{Average} & 42.3 & 51.2 & 6.5 & 35.2 & 83.2 \\
		\bottomrule
	\end{tabular}
\end{table}

Table~\ref{tab:datagen_statistics} summarizes the outcome statistics of the data filtering process across all datasets involved in our data generation.
For each dataset, the first three columns report the rates of each possible outcome per solution attempt:
successfully generated samples (Accept),
samples whose Dedalus solution contains \texttt{NaN} values,
and samples rejected because $\|u\|_{L^\infty}>10$.
The last two columns report the average wall-clock time per individual solution attempt and the expected total wall-clock time needed to produce one accepted sample (including all preceding failed attempts), respectively.

It should be noted that the filtering operations and non-physical artifact issues
discussed above mainly concern the pretraining data.
For the three specific PDE datasets used in our comparison experiments
(Supplementary Note~\ref{app:dataset_forward}),
the filtering conditions were rarely triggered during data generation,
indicating that solutions containing \texttt{NaN} or having excessively large magnitudes are almost non-existent.
Furthermore, we did not observe any visually apparent artifacts in these datasets.
We therefore consider these data to be relatively reliable,
and the comparative conclusions drawn from different models are unlikely to be
significantly affected by data quality issues.

\subsection{Remark on impact of imperfect pretraining data}\label{sec:impact_imperfect}

We emphasize that the core contribution of this paper lies in the model we develop, rather than in the construction of the pretraining dataset.
Although our pretraining data contain certain imperfections, their impact on the main contributions of this work is relatively limited.

As discussed above, the pretraining data generation process introduces two primary types of imperfection.
First, \emph{selection bias} arises because the Dedalus solver occasionally fails to produce valid solutions (outputs containing \texttt{NaN}) for certain PDEs that nevertheless admit theoretical solutions.
These samples are necessarily discarded, causing the pretraining dataset to underrepresent PDE instances that are numerically challenging for traditional solvers.
Second, \emph{label noise} originates from non-physical artifacts appearing in some Dedalus-generated solutions.
These artifacts were not manually corrected or automatically filtered out, and may still be present in the pretraining data.

In the comparative experiments presented in ``Adaptation to specific PDEs'' in the main text, all models (PDEformer-2 and the baselines) are trained or finetuned on the same specific PDE datasets and evaluated on the same test sets.
As noted in the previous subsection, the specific PDE datasets exhibit relatively reliable data quality, with \texttt{NaN}-containing or excessively large-magnitude solutions being almost non-existent and no visually apparent artifacts observed.
Therefore, the comparative conclusions derived from these experiments are relatively trustworthy.
We acknowledge that further improving the pretraining data quality by reducing selection bias and label noise could potentially lead to stronger model performance.
Nevertheless, even under the current imperfect pretraining data conditions, PDEformer-2 has already demonstrated meaningful advantages over existing baseline models.

It should be noted that the selection bias discussed here primarily concerns large-scale, diverse pretraining datasets that encompass a wide variety of PDE forms, rather than training or test sets constructed for individual specific PDEs.
Selection bias in a test set would undermine the reliability of evaluation conclusions, while such bias in a task-specific training set may prevent the model from acquiring complete knowledge of that particular PDE.
Although the selection bias induced by numerical solver difficulty does affect the pretraining dataset, its practical impact is far less significant than the inherent limitations introduced by the pretraining dataset design.
Specifically, despite encompassing eight diverse PDE families, the current pretraining data still exclude many important equation forms, such as compressible Euler equations.
Furthermore, the uniformly distributed scalar coefficients in $[-1,1]$,
the Gaussian random fields for spatially varying coefficients,
and the regular-shaped computational domains
are far from representing the diverse physical scenarios and parameter distributions encountered in real-world applications.
Therefore, mitigating this specific type of solver-induced selection bias is not treated as a priority in our current research agenda for future directions.

We also note that the presence of selection bias and label noise in the pretraining dataset serves, to some extent, as a testbed for evaluating model robustness.
We anticipate that future PDE foundation models trained at even larger scales will need to integrate data from a broader range of sources.
Different datasets may employ different numerical solvers, adopt different PDE formulations, and consequently introduce different types of label noise, or experience numerical solver failures due to various causes that lead to selection bias.
If a model cannot adequately cope with such imperfections in the pretraining data, practitioners may need to design tailored strategies for each individual dataset and solver to screen or repair noisy labels, and to recover samples discarded due to numerical solver failures.
This would incur substantial additional human effort.
Conversely, a model that exhibits inherent robustness to label noise and solver-induced selection bias can substantially alleviate such costs.

By pretraining on this imperfect dataset, we are able to evaluate the model's ability to handle both types of data imperfection.
As demonstrated in ``Robustness to non-physical artifacts'' in the main text and in ``Diagnosis of large solution discrepancy'' (Supplementary Note~\ref{app:training}),
the inductive bias of the model architecture appears to enable it to partly overcome the negative effects of label noise in the pretraining data.
Regarding the selection bias caused by numerical solver difficulty, although occasional success cases have been observed,
as shown in ``Performance on failure cases of the teacher solver'' in the main text and ``More about failure cases of the teacher solver'' (Supplementary Note~\ref{app:training}),
the overall average performance remains unsatisfactory (Supplementary Fig.~\ref{fig:dedalus_failed_histogram}).
We look forward to more advanced model designs in future work that exhibit stronger robustness against both types of imperfections in pretraining data.
}

\section{Training settings and further results}\label{app:training}
\subsection{Pretraining}\label{app:pretraining}
\paragraph{Dynamic data buffer}
Pretraining is executed using the deep learning framework MindSpore%
\footnote{\url{https://www.mindspore.cn/}},
on a machine with 8 Ascend 910A NPUs, 192 CPU cores, 768 GB of RAM, and a local disk with about 3 TB storage.
The entire pretraining dataset stored on the remote large disk has a size of about 40 TB, which does not fit into the capacity of our local disk.
In order to resolve this, we maintain a data buffer in the local disk, which contains about 50 data files (with 1k samples in each file).
Each pretraining ``epoch'' only traverses the data samples stored in this buffer.
In parallel to training, new data files are downloaded from the remote disk, and will replace the existing files in the local buffer after downloading is completed.

\paragraph{Continual learning}
The large diverse dataset took a relatively long time for us to prepare (including PDE design, coding, and solver execution).
We decide not to postpone our pretraining until the entire dataset is ready.
Instead, we train the model with all the data available at the moment, and prepare new data in the mean time.
When new data is prepared, we load the previous model weights, and improve it using the extended dataset.
This procedure is repeated for over ten times, with each incremental training run for about one or two days.
\verB{We believe that such capability of being continuously improved with more data could be an important feature of a PDE foundation model in practice.}

\paragraph{More settings}
Each incremental training stage during pretraining runs from 200 to 500 ``epochs'' with batch size 80 (10 per device).
The Adam optimizer~\cite{Adam} without weight decay is used, clipping the gradient to a maximal norm of 1.
The maximum learning rate is selected as $10^{-4}$ in the early stages and $6\times 10^{-5}$ for the rest of stages.
Learning rate is linear warmed-up in the first 10 ``epochs'', and decays following a cosine schedule afterwards.
Before training starts, we preprocess data by constructing the computational graph for the specific PDE forms of the data samples,
and storing the information in \verB{separate} files.

Let the reference solution (of a single variable) $u^\text{ref}$ be recorded on the spatio-temporal points $\{(t_i,\boldsymbol{r}_i)\mid i\in\mathcal{I}\}$.
During evaluation, the model accuracy is measured by the normalized root mean square error (nRMSE)
\begin{equation}\label{eq:nrmse}
	\mathrm{nRMSE}(\hat u,u^\text{ref})=\frac{
		\sqrt{\frac{1}{|\mathcal{I}|}\sum_{i\in\mathcal{I}}(\hat u(t_i,\boldsymbol{r}_i)-u^\text{ref}|_i)^2}}
		{\epsilon+
		\sqrt{\frac{1}{|\mathcal{I}|}\sum_{i\in\mathcal{I}}(u^\text{ref}|_i)^2}}
,\end{equation}
in which the small constant $\epsilon=10^{-6}$ is included to prevent division by zero.
Two-dimensional PDE solutions often involve a lot of discretization points (typically $101\times 128\times 128$ for our pretraining dataset),
leading to a considerable computational cost for our INR decoder that deals with each grid point independently.
Therefore, to reduce the computational cost, the training loss is taken to be a Monte-Carlo estimate of~\eqref{eq:nrmse}.
Specifically, we randomly sample $K=8192$ independent indices $i_k\sim U(\mathcal{I}), k=1,\dots,K$, and take
\begin{equation}\label{eq:loss}
	\mathrm{nRMSE}(\hat u,u^\text{ref};\{i_k\})=\frac{
		\sqrt{\frac{1}{K}\sum_{k=1}^{K}(\hat u(t_{i_k},\boldsymbol{r}_{i_k})-u^\text{ref}|_{i_k})^2}}
		{\epsilon+
		\sqrt{\frac{1}{K}\sum_{k=1}^{K}(u^\text{ref}|_{i_k})^2}}
\end{equation}
as the training loss.
Following our previous PDEformer-1~\cite{PDEformer1}, the small constant is increased to $\epsilon=0.05$ to \verB{stabilize} training.
For PDEs with $n_v$ variables, the loss is computed for one variable at a time, and a data sample is used $n_v$ times in every training ``epoch''.

\paragraph{Detailed dataset accuracies}
\begin{figure}[htbp!]
	\centering
	\begin{subfigure}[p]{0.45\textwidth}
		\centering
		\includegraphics[width=\textwidth]{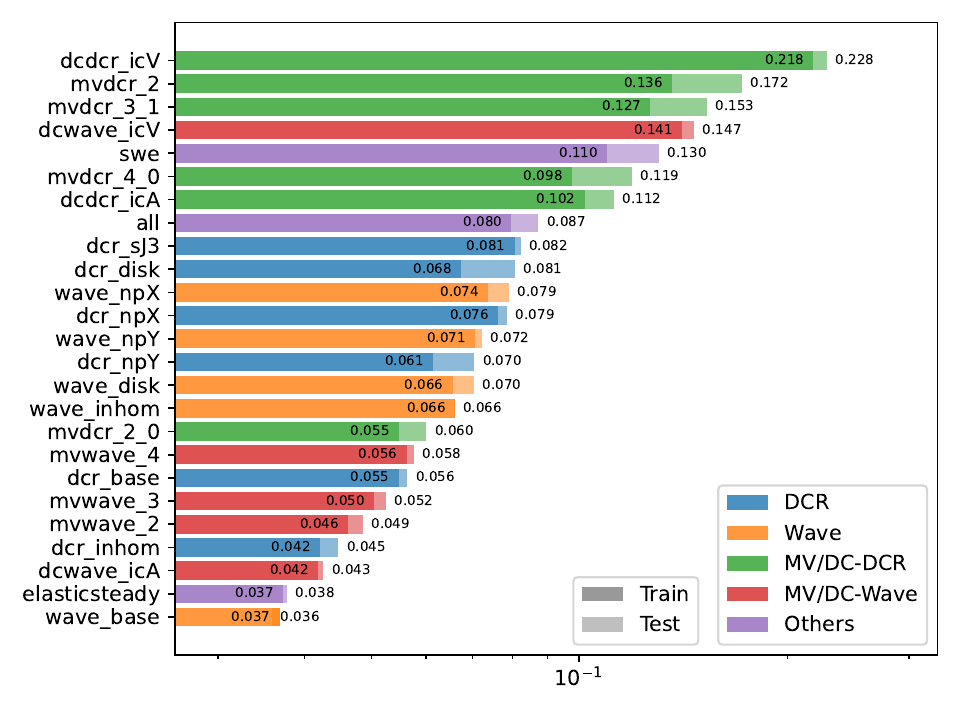}
		\caption{}
	\end{subfigure}%
	~ 
	\begin{subfigure}[p]{0.45\textwidth}
		\centering
		\includegraphics[width=\textwidth]{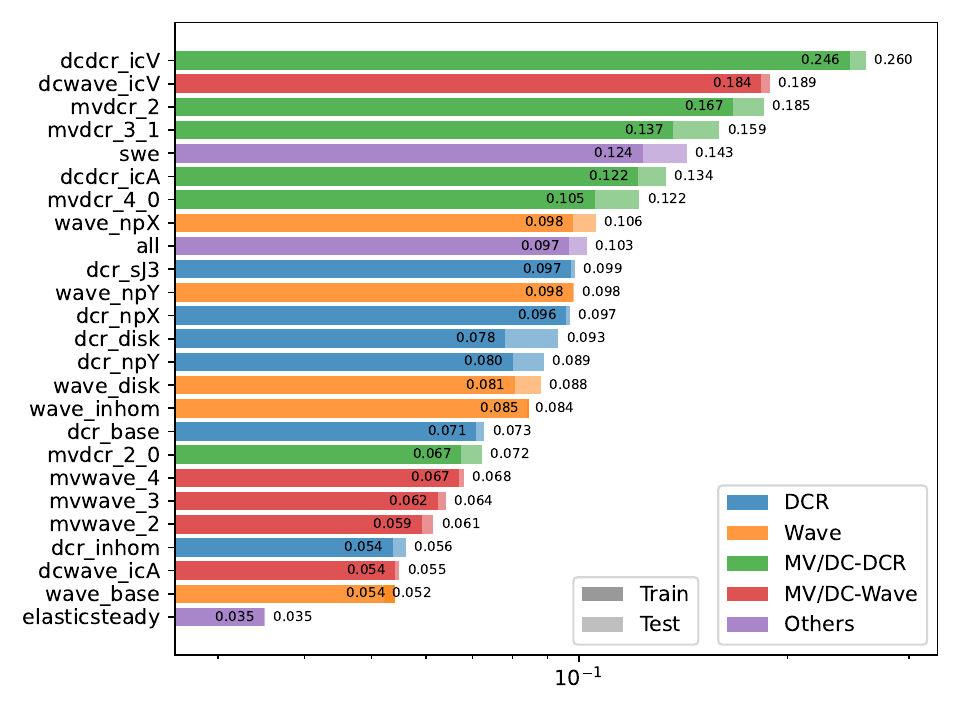}
		\caption{}
	\end{subfigure}%
	\caption{Detailed model accuracy (nRMSE) on the pretraining dataset,
		including (a) PDEformer-2-base and (b) PDEformer-2-fast models.
	}%
	\label{figs:pretrain_err_bars}
\end{figure}
Fig.~\ref{figs:pretrain_err_bars} shows the detailed model accuracy (measured by~\eqref{eq:nrmse}) post pretraining.
Each test dataset contains 1k samples.
Although most of the training datasets contain more than 100k samples, we only use the first 1k samples to compute the mean nRMSE for each dataset,
since the entire dataset does not fit into the capacity of our local disk.

\verB{
\paragraph{Estimation of equivalent pretraining ``epochs''}
The continual learning procedure described above, involving a lot of incremental training stages, makes it impractical to maintain a detailed record of the entire pretraining history.
To obtain a more principled estimate of the required pretraining resources,
we conduct a controlled from-scratch pretraining experiment using the PDEformer-2-base architecture.
The model is trained on the same complete pretraining dataset, starting from randomly initialized weights.
The maximum learning rate is set to $10^{-4}$, with a 10-``epoch'' linear warmup and exponential decay towards $0.04\times 10^{-4}$ over the originally planned 1000 ``epochs''.
Due to an unexpected interruption, training was terminated before reaching 1000 ``epochs''; the data shown in Fig.~\ref{fig:pretrain_scaling} covers approximately the first 500 ``epochs''.
Each training ``epoch'' takes approximately 7 minutes and 40 seconds on the same hardware as the main pretraining.

During training, we monitor the model performance on three representative test datasets:
\texttt{dcr\_base}, \texttt{wave\_base}, and \texttt{mvdcr\_2\_0}.
Evaluation is performed at ``epochs'' $\{N_k\}$ where $\sqrt{N_k}$ are uniformly spaced between $0$ and $\sqrt{1000}$ (with $K=11$ evaluation points, i.e., $N_k=(k\sqrt{1000}\,/\,10)^2$ for $k=0,\dots,10$).
This schedule ensures the measured points are well distributed across the log-scale axes of Fig.~\ref{fig:pretrain_scaling},
in contrast to periodic evaluation which would cause clustering at later ``epochs''.

\begin{figure}[htbp]
	\centering
	\includegraphics[width=0.6\linewidth]{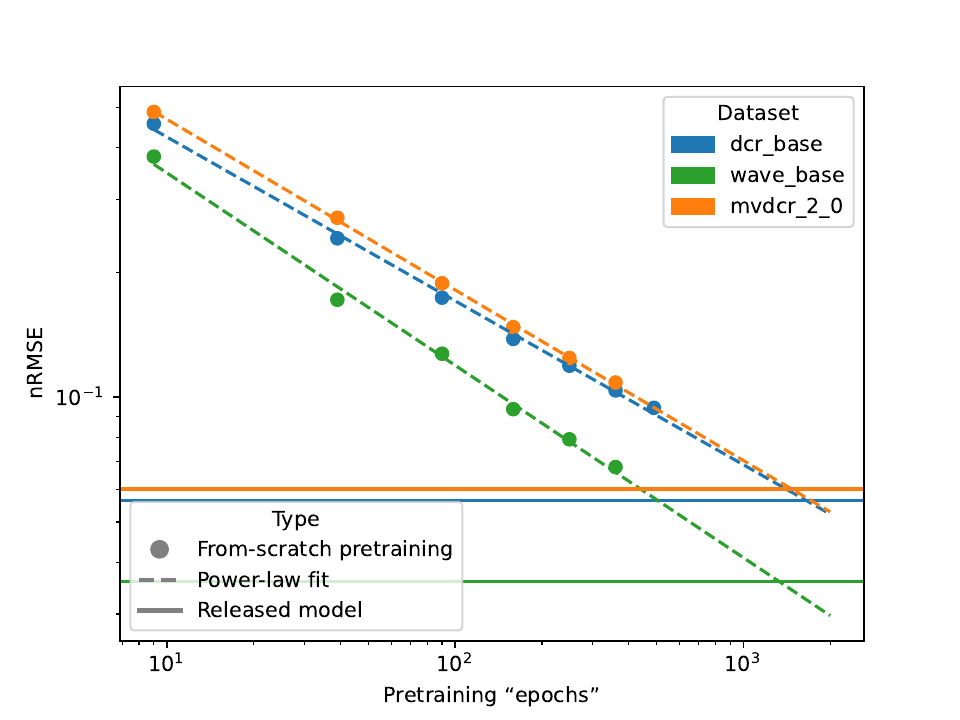}
	\caption{Estimating the equivalent pretraining ``epochs'' of PDEformer-2-base.
		The nRMSE on three representative test datasets is monitored at selected intervals during from-scratch pretraining.
		Recorded values (scatter points) follow a power-law trend with respect to the number of training ``epochs''.
		Power-law fits (dashed lines) are extrapolated to estimate the ``epochs'' required to match the error level of the released model checkpoint (solid horizontal lines).
	}
	\label{fig:pretrain_scaling}
\end{figure}

As shown in Fig.~\ref{fig:pretrain_scaling}, despite the early termination, the recorded errors already exhibit a clear power-law relationship with the number of training ``epochs'' in the log-log scale.
We fit a power-law curve $\mathrm{nRMSE} \propto N^{-\alpha}$ (where $N$ is the number of ``epochs'') to the measured points of each dataset,
and extrapolate the fitted lines to estimate the total training ``epochs'' required to reach the error level of our released model checkpoint.
The extrapolated intersection points are approximately 1649, 1331, and 1470 ``epochs'' for \texttt{dcr\_base}, \texttt{wave\_base}, and \texttt{mvdcr\_2\_0}, respectively.
Taking a conservative overall estimate of 1500 ``epochs'', and multiplying by the ``per-epoch'' time of 7 minutes and 40 seconds,
the estimated total from-scratch pretraining time is approximately 192 hours, or about 8 days.
}

\paragraph{Diagnosis of large solution discrepancy}
\begin{figure}[htbp]
	\centering
	\includegraphics[width=0.9\linewidth]{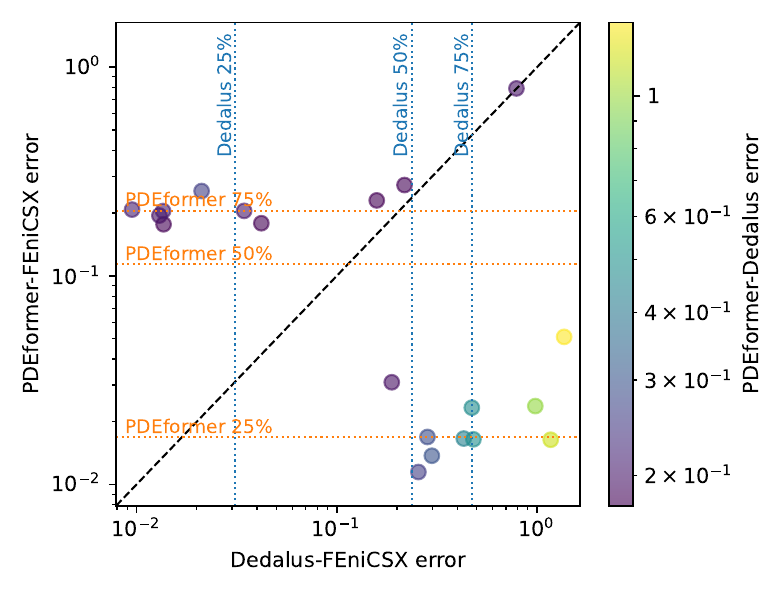}
	\caption{Scatter plot of nRMSE compared with FEniCSx solutions on the 20 selected cases where Dedalus and PDEformer-2-base predictions deviate significantly.
		Point color stands for the nRMSE between the Dedalus solution and the PDEformer-2 prediction.}
	\label{fig:nonphysical_pattern_scatter}
\end{figure}
When PDEformer-2 predictions deviate significantly from Dedalus solutions on in-distribution PDEs,
we may identify two possible causes:
(1) PDEformer-2 does not generalize well to this case,
(2) Dedalus itself produces non-physical patterns as mentioned in `Pretraining results' in the main text.
To make an investigation of this, we focus specifically on the DCR (more accurately, \verb|dcr_base|) dataset,
and select the 20 samples exhibiting the largest nRMSE (PDEformer-2-base vs. Dedalus) from the test set (with 1000 samples in total).
We regenerate the solutions using the FEniCSx solver, and employ them as the new reference solutions to compute the nRMSE.
According to Fig.~\ref{fig:nonphysical_pattern_scatter}, for the 10 samples with the largest solution discrepancies
(with point color ranging from yellow to blue), the error is more likely attributed to the unreliable solutions produced by Dedalus (case (2)),
since most of the points are located in the lower-right corner of the figure.
As for the remaining 10 samples with relatively smaller discrepancies (points in purple),
the error tends to originate from the inefficient generalization of PDEformer-2,
as most of the points lie in the upper-left part.
Our results indicate that, although a neural solver typically exhibits a certain level of generalization error for most cases,
the worst-case performance could potentially be superior to its teacher solver,
since the model inductive bias can prevent it from producing highly non-physical predictions.
Such extreme-case reliability might be crucial in some real scenarios.
Indeed, we acknowledge that the above investigation is still limited.
More sufficient and \verB{rigorous} validation of the model would be necessary for future deployment in real applications.

\paragraph{More about failure cases of the teacher solver}
\begin{figure}[htbp]
	\centering
	\includegraphics[width=0.9\linewidth]{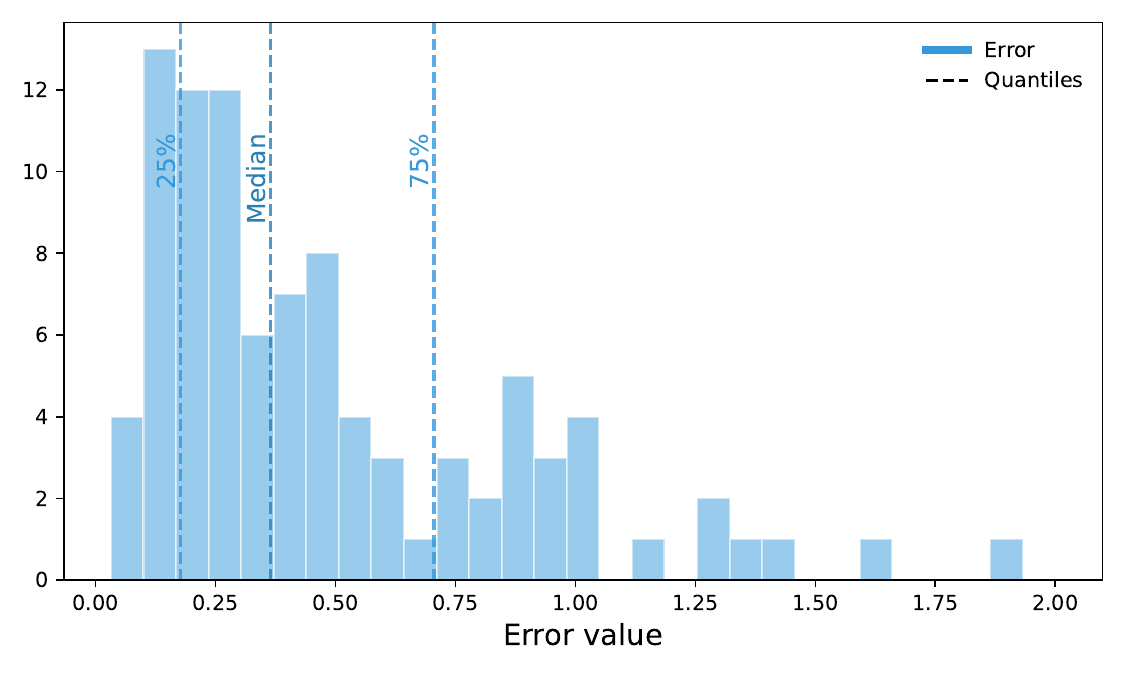}
	\caption{Histogram of PDEformer-2-base prediction nRMSE compared with FEniCSx solutions on \verB{94 filtered} failure cases of Dedalus.
	}
	\label{fig:dedalus_failed_histogram}
\end{figure}
\begin{figure}[tbp]
	\centering
	\includegraphics[width=0.9\linewidth]{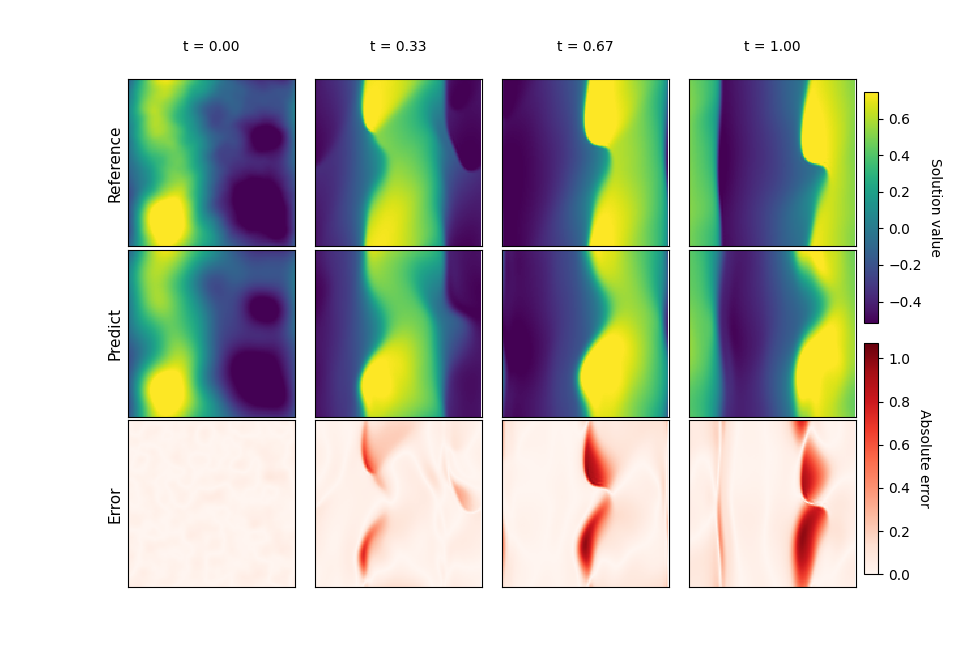}
	\[u_t-\tilde{a}\nabla\cdot(\tilde{a}\nabla u)+s(r)+c_{01}u+c_{02}u^2+u^3
	+(c_{11}u+c_{12}u^2+c_{13}u^3)_x+(u^2+c_{23}u^3)_y=0\]
	\caption{Prediction of PDEformer-2-base when Dedalus fails, showing a representative case.
		We employ FEniCSx to prepare the reference solution.
		While our model successfully captures the overall physical phenomena, its imperfect handling of shock evolution leads to observable prediction errors.
	}
	\label{fig:dedalus_failed_sample_avg}
\end{figure}
We have mentioned in `Pretraining results' in the main text that Dedalus fails to solve some PDEs during pretraining data generation.
To make a further investigation of these cases, we collect a total of 100 such failure cases when generating the \verb|dcr_base| dataset,
\verB{
and regenerate their solutions using the FEniCSx solver.
Although FEniCSx successfully produced numerical solutions for all 100 cases,
we observed that some of these re-generated solutions exhibit anomalously large amplitudes, with the maximal $L^\infty$ norm reaching the order of $10^6$.
We speculate that such cases may correspond to PDEs that do not admit a valid solution in the theoretical sense
as discussed in Supplementary Note~\ref{sec:data_filtering},
and the finite numerical values were obtained only due to discretization errors and possible numerical safeguard mechanisms in the solver,
which prevented the appearance of \texttt{NaN} values.
Taking this into account, we further filter the FEniCSx re-generated solutions by excluding those with $\|u\|_{L^\infty}>10$,
resulting in the removal of 6 cases.
We then compare PDEformer-2-base predictions against the FEniCSx reference solutions on the remaining 94 cases.
}
As can be seen in the histogram in Fig.~\ref{fig:dedalus_failed_histogram},
PDEformer-2 can generalize to some of the cases, but the overall performance would seem unsatisfactory.
A representative case is shown in Supplementary Fig.~\ref{fig:dedalus_failed_sample_avg}.

\subsection{Baseline experiments on specific PDEs}\label{app:adaptation}
If not specified otherwise, all experiments on specific PDEs uses MindSpore, with a maximum learning rate of $10^{-4}$ for the Adam optimizer.
We use 8 NPUs with batch size 40 (5 per device) for $80$ or more training samples, and a single NPU otherwise,
with batch size 5 for 20 samples, 2 for 4 samples.
For the baseline models, all solution components are predicted concurrently, so the evaluation nRMSE~\eqref{eq:nrmse} takes into account the entire solution.
This is different from PDEformer-2, which computes~\eqref{eq:nrmse} independently for different solution components, and then takes the average to obtain the final prediction error.
The training loss~\eqref{eq:loss} is computed in a similar way, with the grid point sampling applied only to PDEformer-2 and DeepONets.
The remaining optimization settings are shared with pretraining.

For the two baseline foundation models, Poseidon and MPP, we use the official PyTorch implementation adapted for GPUs.
Finetuning starts from the official pretrained model weights.
Note that Poseidon and MPP do not predict the entire spatio-temporal solution directly in a single inference.
For Poseidon, we feed the initial value into the model, and let it predict all the system states afterwards, so a trajectory with $n_t+1$ time snapshots will be used $n_t$ times in each epoch.
\verB{Its finetuning is executed on 8 GPUs with batch size $192$ (24 per device) for Poseidon-B,
and $64$ (8 per device) for Poseidon-L.
To mitigate the difference introduced by a smaller batch size, we set the gradient accumulation steps to $3$ for Poseidon-L,
since the effective global batch size is the product of per-device batch size, the number of devices, and the number of gradient accumulation steps.}
For MPP, we randomly select two consequent snapshots from each trajectory to train the model,
with the former snapshot serving as model input and the latter as the output label.
Each trajectory is used $2.5$ times in average in every ``epoch''.
The finetuning uses a single GPU with batch size $10$.
As less hardware resources are used to finetune MPP, we divide its time consumption by 8 before preparing Fig.~\ref{figs:ablation}(b) in the main text. 
The remaining configurations for finetuning are the same as in the original implementation.

Training or finetuning runs for 1k epochs, except for the Poseidon (200 epochs) and MPP (500 ``epochs'')
results in Fig.~\ref{figs:ablation}(a) in the main text, since their special training strategy leads to a longer epoch time.
To avoid showing results with overfitting,
\verB{
we compute the error value~\eqref{eq:nrmse} on the validation set with $100$ samples for every $50$ epochs, and record the best model checkpoint.
The reported performance, using this selected checkpoint, is then computed on an independent test set with $100$ samples.
}

Geo-FNO runs on a single GPU using the official PyTorch implementation.
For each model, the software and hardware platforms used remain consistent in all the experiments involved, including few-shot adaptation, adaptation speed, and prediction speed.
When we assess the prediction speed of PDEformer-2, \emph{solely preparing an INR} only requires producing the embedding vectors $\mu_j^\ell$ mentioned in Supplementary Note~\ref{app:arch_graph_encode}
(i.e., forward propagation of the graph Transformer), without making any spatio-temporal queries for the INR decoder that follows.

\verB{
\subsubsection{Controlled comparison across OOD scenarios}
\label{app:ood_stress}

To investigate how different types of distribution shift affect PDEformer-2's generalization, we conduct a controlled comparison across several out-of-distribution (OOD) scenarios.
The incompressible Navier-Stokes (INS) equation with periodic boundaries is selected as the base PDE system.
All experiments use 80 training samples and compare three configurations:
PDEformer-2-fast zero-shot (direct inference without any finetuning),
PDEformer-2-fast finetuned, and an FNO trained from-scratch on each dataset.
Both PDEformer-2-fast finetuned and FNO are trained for 1,000 epochs with validation-based model selection.
The FNO serves as a difficulty-calibrated baseline: as a specialized neural operator without any pretraining, its performance reflects only the intrinsic difficulty of fitting the target data, unaffected by distribution shift relative to any pretraining distribution.
PDEformer-2's zero-shot performance reflects the raw transferability of the pretrained model, while its finetuned performance additionally benefits from task-specific adaptation.
Comparing the three helps separate data difficulty, pretraining transferability, and the gains from finetuning.

\paragraph{Stress axis 1: GRF smoothness $\alpha$.}
The initial conditions of the PDEs are generated as Gaussian random fields (GRFs) with covariance operator $50^2(-\Delta + 5^2 I)^{-\alpha}$, where $\alpha$ controls the spatial smoothness.
Larger $\alpha$ produces smoother fields with long-range correlations, while smaller $\alpha$ yields rougher fields with finer oscillations.
We test $\alpha \in \{3.5, 3.0, 2.8, 2.5, 2.3, 2.2\}$, with $\alpha=3.0$ (aligned with~\eqref{eq:grf}) as the in-distribution setting.
Fig.~\ref{fig:grf_by_alpha} visualizes representative GRF samples at each $\alpha$, allowing qualitative assessment of how the spatial oscillation pattern evolves as $\alpha$ varies.

\begin{figure}[htbp]
	\centering
	\begin{subfigure}[b]{\linewidth}
		\centering
		\includegraphics[width=0.8\linewidth]{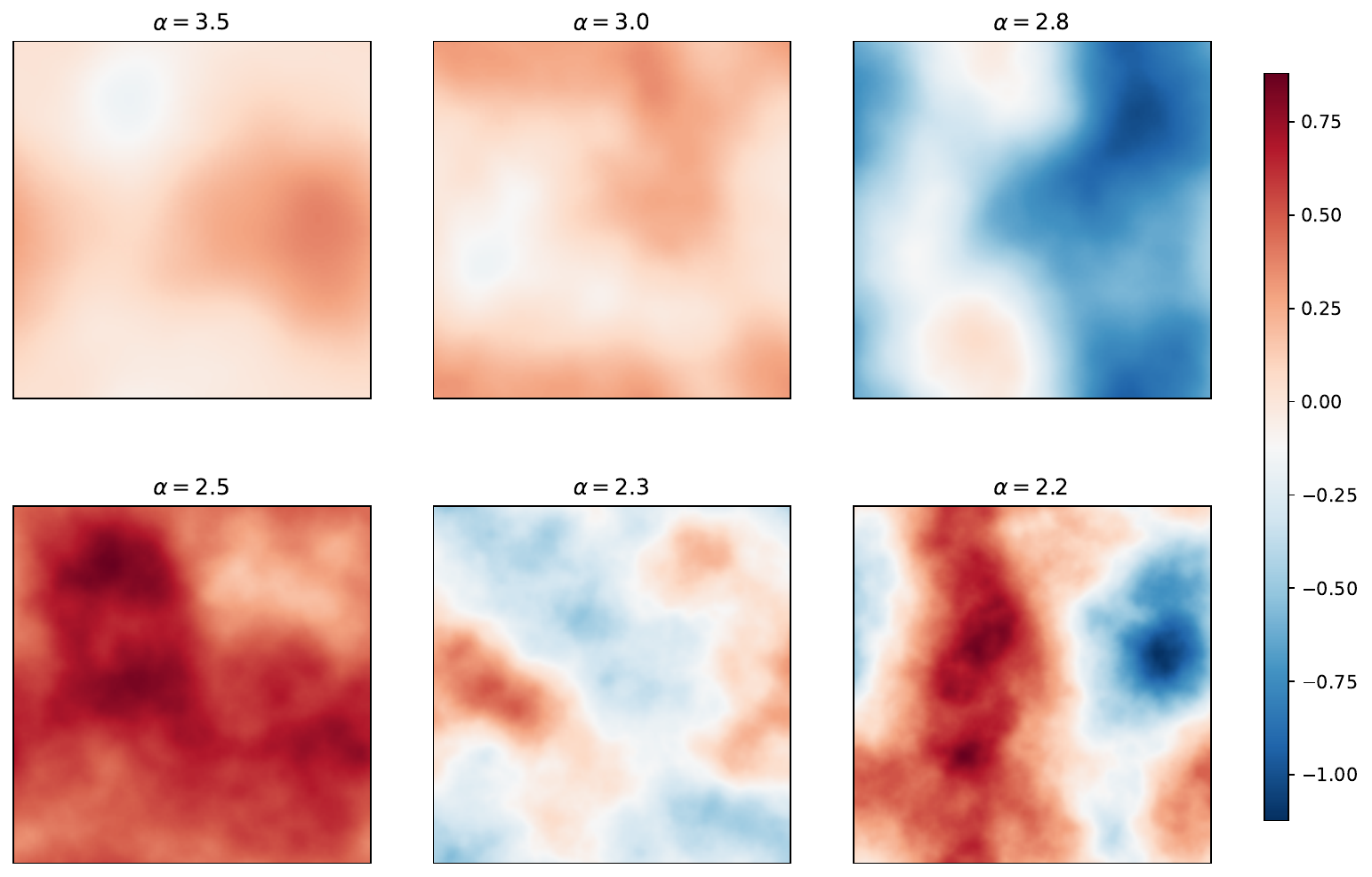}
		\caption{Raw fields (as generated, without preprocessing).}
		\label{fig:grf_raw}
	\end{subfigure}
	\\[6pt]
	\begin{subfigure}[b]{\linewidth}
		\centering
		\includegraphics[width=0.8\linewidth]{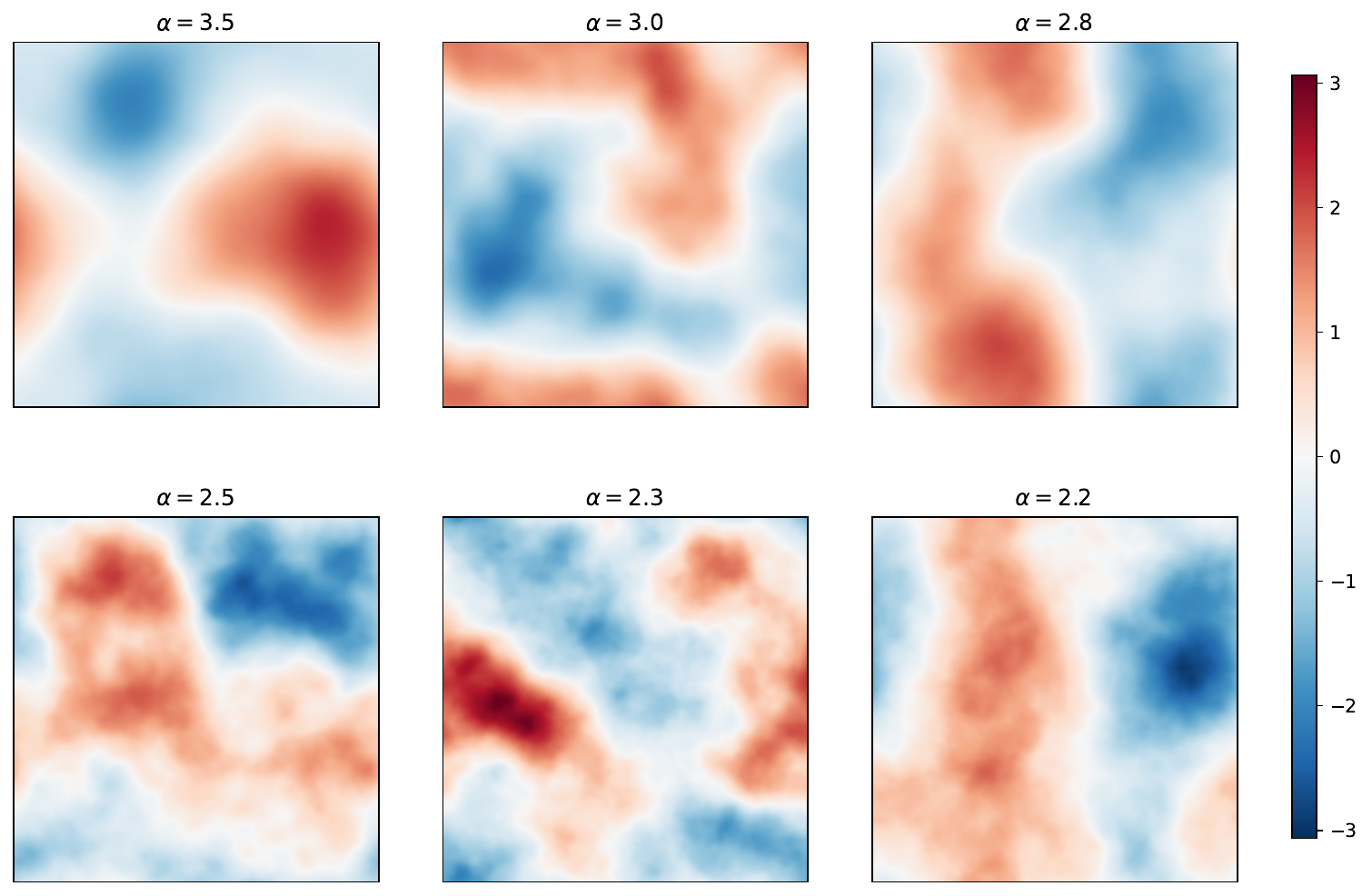}
		\caption{Normalized fields (zero mean, unit variance, shared colorbar).}
		\label{fig:grf_norm}
	\end{subfigure}
	\caption{Representative Gaussian random field samples for different smoothness parameters $\alpha$.
		(a) Raw fields show the original amplitude variation across $\alpha$.
		(b) Normalized fields use a common color scale, facilitating comparison of spatial oscillation patterns.
		As $\alpha$ decreases, the fields exhibit progressively finer spatial structures.
	}
	\label{fig:grf_by_alpha}
\end{figure}

\paragraph{Stress axis 2: Physical composition.}
We also include OOD scenarios arising from new combinations of known physical mechanisms:
coupling tracer transport with the incompressible Navier-Stokes equations (INS-Tracer),
and applying no-slip boundary conditions on the top and bottom pipe walls to the same system (INS-Pipe).
These correspond to the two ``new combination'' datasets discussed in the main text.

\paragraph{Results.}
Fig.~\ref{fig:finetune_by_ood_type} presents the test nRMSE across all eight scenarios.

\begin{figure}[htbp]
	\centering
	\includegraphics[width=\linewidth]{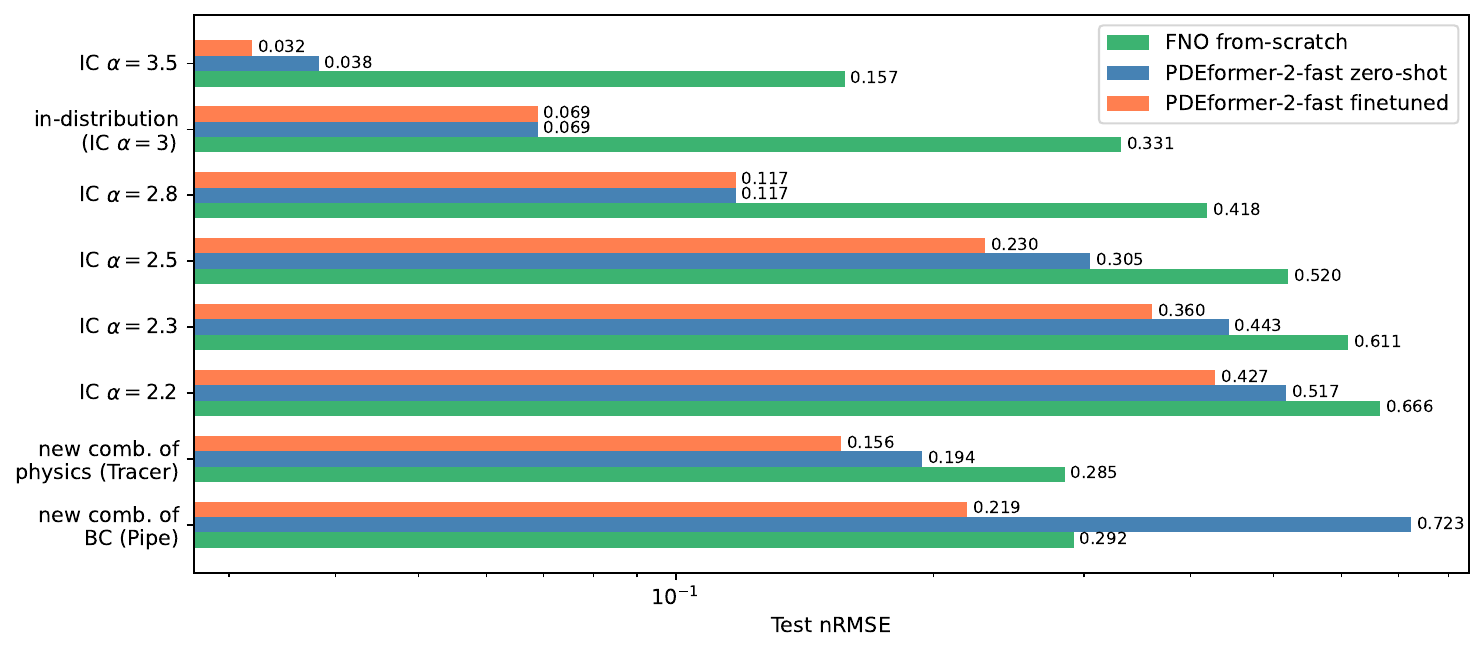}
	\caption{Test nRMSE of FNO (trained from-scratch), PDEformer-2-fast zero-shot, and PDEformer-2-fast (finetuned) under different OOD scenarios.
		The INS equation with periodic boundaries is selected as the base PDE system.
		Training and finetuning use 80 samples, while zero-shot uses no training data.
		FNO reflects data difficulty only, zero-shot reflects the raw transferability of the pretrained model, and finetuning additionally benefits from task-specific adaptation.
	}
	\label{fig:finetune_by_ood_type}
\end{figure}

Several observations emerge from the results.
First, the FNO error grows monotonically from $0.157$ at $\alpha=3.5$ (smoother than in-distribution) to $0.666$ at $\alpha=2.2$ (highly oscillatory), confirming that rougher GRFs are intrinsically more challenging to fit.
This trend aligns with the visual progression in Fig.~\ref{fig:grf_by_alpha}: lower $\alpha$ introduces finer spatial structures that demand higher model capacity to resolve.

Second, PDEformer-2-fast zero-shot already outperforms FNO in seven out of eight scenarios tested, suggesting that the pretrained model carries a certain level of transferable knowledge even without any task-specific adaptation.
The sole exception is INS-Pipe, where zero-shot nRMSE ($0.723$) substantially exceeds FNO ($0.292$).
Note that this scenario involves a new combination of the INS equation with no-slip boundary conditions, a configuration absent from the pretraining data.
The poor zero-shot result suggests that such novel combinations can pose a significant challenge for direct transfer, even when each individual component has been seen during pretraining.

Third, after finetuning, PDEformer-2-fast consistently outperforms FNO across all scenarios, but its relative advantage narrows as $\alpha$ deviates further from the in-distribution value of $3.0$.
The finetuned-to-FNO error ratio increases from approximately $0.20$ at $\alpha=3.5$ to $0.64$ at $\alpha=2.2$, suggesting that the benefits of pretraining diminish under stronger distribution shift.
This is consistent with the intuition noted in the main text that ``transfer efficiency may also depend on similarity to pretraining data.''

Fourth, for the two ``new combination'' scenarios, PDEformer-2 maintains a meaningful advantage over FNO after finetuning (error ratios of $0.55$ for INS-Tracer and $0.75$ for INS-Pipe), indicating that the model can effectively recombine known physical mechanisms even when they appear in unfamiliar configurations.
Notably, INS-Pipe exhibits the most dramatic improvement from zero-shot to finetuned, with nRMSE dropping from $0.723$ to $0.219$, which echoes the main-text observation that finetuning ``efficiently restores superiority'' for this scenario.

Overall, this comparison reveals that PDEformer-2's generalization under distribution shift is governed by both the intrinsic difficulty of the target data and the degree of similarity to the pretraining distribution.
The use of FNO as a difficulty-calibrated baseline
enables a more nuanced understanding than reporting absolute performance numbers in isolation.

\subsection{Prediction speed}\label{app:prediction_speed}
\paragraph{Software and hardware for prediction speed}
Unless otherwise noted, all prediction speed measurements in this section are conducted using MindSpore on a single Ascend 910B NPU,
processing one PDE sample at a time (i.e., a single sample per batch, while the number of time steps or spatio-temporal query points within that sample may vary as needed).
Note that this hardware differs from the pretraining setup (Supplementary Note~\ref{app:pretraining}), which conducts distributed data-parallel training on a machine with 8 Ascend 910A NPUs.
For the baseline foundation models Poseidon and MPP, which only support PyTorch, their inference times are originally measured on an NVIDIA V100 GPU.
To enable a fairer comparison across these distinct platforms, we calibrate the platform difference using DeepONet as a reference model and apply a rescaling procedure (described at the end of this section) to estimate equivalent times on our MindSpore+Ascend 910B setting.
}

\paragraph{Traditional solver resolution}
When we reduce the grid resolution for Dedalus, the error is computed against the resolution-128 results, consistent with the neural solvers.
If we use solutions with a higher resolution for training and testing, we speculate that the computed errors would become larger for low-resolution traditional solvers, and do not change significantly for neural solvers.

\paragraph{Solely preparing an INR}
Solely preparing an INR corresponds to the case in which
we only need to make queries at a very limited number of spatio-temporal locations (such as the inverse problem observations),
and the majority of the prediction time can be attributed to preparing an implicit neural representation for the solution.

\verB{
\paragraph{Inference time versus query count}
\begin{figure}[htbp]
	\centering
	\includegraphics[width=\linewidth]{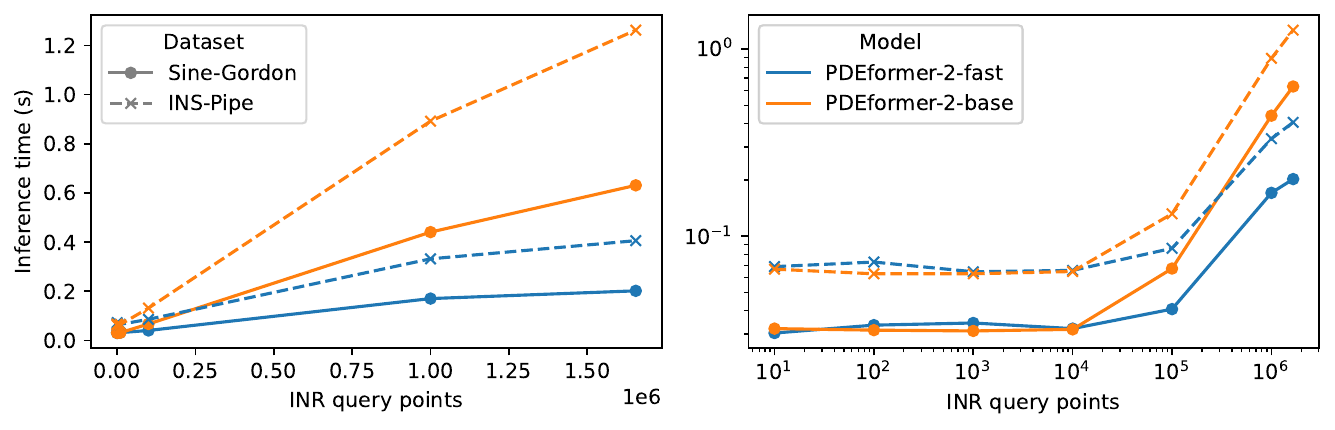}
	\caption{Inference time of PDEformer-2-base and PDEformer-2-fast versus the number of INR query points,
		evaluated on the Sine-Gordon (solid lines, circles) and INS-Pipe (dashed lines, crosses) datasets.
		The left panel uses linear axes, while the right panel adopts log-log scaling.
		For small query counts (up to $\sim 10^4$), the inference time is dominated by the graph Transformer and remains nearly constant.
		As the number of queries increases, the INR decoder cost becomes the primary factor, and the speed advantage of the fast variant over the base variant grows accordingly.
		The rightmost data point on each curve corresponds to the full-solution with $101\times 128^2$ resolution.
	}
	\label{fig:inference_time_vs_query}
\end{figure}
Table~\ref{tab:inference_time} in the main text reports inference time under three representative
query modes: full spatio-temporal output, final time step only, and solely preparing an INR.
To complement these discrete modes, Fig.~\ref{fig:inference_time_vs_query} provides a continuous view of how inference time varies with the number of INR query points, spanning from ten queries to the full $101\times 128^2$ solution.
The nearly flat region at low query counts confirms that the time for solely preparing an INR
(approximately $0.02$~s in Table~\ref{tab:inference_time}) accurately reflects the query-free cost.
Beyond approximately $10^4$ query points, the INR decoder queries start to contribute noticeably
to the total cost, and the per-sample time at full resolution is consistent with full
spatio-temporal output in Table~\ref{tab:inference_time}.
Fig.~\ref{fig:inference_time_vs_query} thus bridges the discrete query modes in the main text, providing a more complete characterization of PDEformer-2's inference behavior.

\paragraph{Cross-platform fair comparison of prediction speed}
As noted at the beginning of this section, the baseline foundation models Poseidon and MPP are implemented in PyTorch and their inference times are originally measured on an NVIDIA V100 GPU,
while our PDEformer-2 and other baselines (DeepONet, FNO, U-Net) use MindSpore on an Ascend 910B NPU.
To enable a fairer comparison across these distinct platforms, we calibrate the platform difference using DeepONet as a reference model,
since DeepONet can be executed on both platforms with identical network architecture and inference logic,
and both its CNN-branch and MLP-branch variants have similar computational characteristics.
We acknowledge that the speed comparison derived from DeepONet may not perfectly transfer to Transformer-based architectures with different computational profiles.
The rescaled estimates to be presented below should therefore be interpreted as approximate.

\begin{table}[htbp]
	\centering
	\caption{Inference time (seconds) of DeepONet under different software and hardware configurations, used as calibration for cross-platform rescaling.
		All measurements process a single PDE sample at a time (sample batch size~1).
		``OOM'' indicates out-of-memory on the GPU, preventing full-resolution inference.
		The platform ratio $0.024/0.040=0.6$ (MindSpore+Ascend~910B over PyTorch+V100 at $1/8$ resolution) is adopted as the rescaling factor.
	}
	\label{tab:deeponet_calib}
	\begin{tabular}{lllcc}
		\toprule
		Framework & Hardware & Query mode & CNN (s) & MLP (s) \\
		\midrule
		MindSpore & Ascend~910B & Full ($101\times 128^2$) & 0.086 & 0.085 \\
		MindSpore & Ascend~910B & $1/8$ resolution & 0.024 & 0.024 \\
		PyTorch & V100 & $1/8$ resolution & 0.040 & 0.040 \\
		PyTorch & V100 & Full ($101\times 128^2$) & OOM & OOM \\
		\bottomrule
	\end{tabular}
\end{table}
\begin{table}[htbp]
	\centering
	\caption{Original (PyTorch+V100) and platform-rescaled (estimated MindSpore+Ascend~910B) per-sample inference time of Poseidon and MPP,
		measured on the Sine-Gordon (SG) and INS-Pipe datasets.
	Rescaled values are obtained by multiplying the original V100 times by $0.6$ (Table~\ref{tab:deeponet_calib}).}
	\label{tab:platform_rescale}
	\begin{tabular}{llrrrr}
		\toprule
		Model & Inference strategy & \multicolumn{2}{c}{Original (s)} & \multicolumn{2}{c}{Rescaled (s)} \\
		& & SG & INS-Pipe & SG & INS-Pipe \\
		\midrule
		Poseidon-B & all timesteps (T-bsz$=100$) & 0.402 & 0.403 & 0.241 & 0.242 \\
		Poseidon-B & per-step (T-bsz$=1$) & 11.86 & 11.95 & 7.12 & 7.17 \\
		Poseidon-B & final time step only & 0.120 & 0.120 & 0.072 & 0.072 \\
		Poseidon-L & all timesteps (T-bsz$=100$) & 0.843 & 0.844 & 0.506 & 0.506 \\
		Poseidon-L & per-step (T-bsz$=1$) & 11.99 & 12.15 & 7.19 & 7.29 \\
		Poseidon-L & final time step only & 0.124 & 0.120 & 0.074 & 0.072 \\
		MPP-B & autoregressive rollout & 5.612 & 4.666 & 3.367 & 2.800 \\
		MPP-L & autoregressive rollout & 18.50 & 18.35 & 11.10 & 11.01 \\
		\bottomrule
	\end{tabular}
\end{table}

Table~\ref{tab:platform_rescale} lists the original PyTorch+V100 inference times of Poseidon and MPP together with their platform-rescaled estimates (multiplied by $0.6$).
All measurements process a single PDE sample at a time (sample batch size~1), and the number of time steps packed into one forward pass (time batch size) may vary across strategies.
For Poseidon, we consider three inference modes:
(i) feeding all $n_t=100$ time steps as a single batch (T-bsz$=100$), which maximizes hardware utilization while respecting the single-sample routine;
(ii) processing each time step independently (T-bsz$=1$);
and (iii) predicting only the final time step, corresponding to the ``Final time step only''
entry in Table~\ref{tab:inference_time} in the main text.
Strategy (i) yields faster per-sample inference than strategy (ii) that requires repeated model inference to obtain the full predicted trajectory,
and is adopted for the ``Full spatio-temporal output'' entry in Table~\ref{tab:inference_time}
in the main text.
For MPP, inference is inherently autoregressive: the model predicts one time step at a time and requires sequential rollout to obtain a full trajectory, as explained in Supplementary Note~\ref{app:model_baseline}.

For the adaptation speed comparison in Fig.~\ref{figs:ablation}(b) in the main text,
the finetuning time curves of Poseidon and MPP are similarly rescaled using the factor $0.6$,
assuming that the per-iteration time ratio between the two platforms is comparable for finetuning (which involves both forward and backward passes) and for inference.
}

\subsection{Inverse problems}
\subsubsection{Recovering scalar coefficients}
In this task, we consider diffusion-convection-reaction equations mentioned in Supplementary Note~\ref{app:inv_dcr}.

\begin{figure}[htbp]
	\centering
	\includegraphics[width=\linewidth]{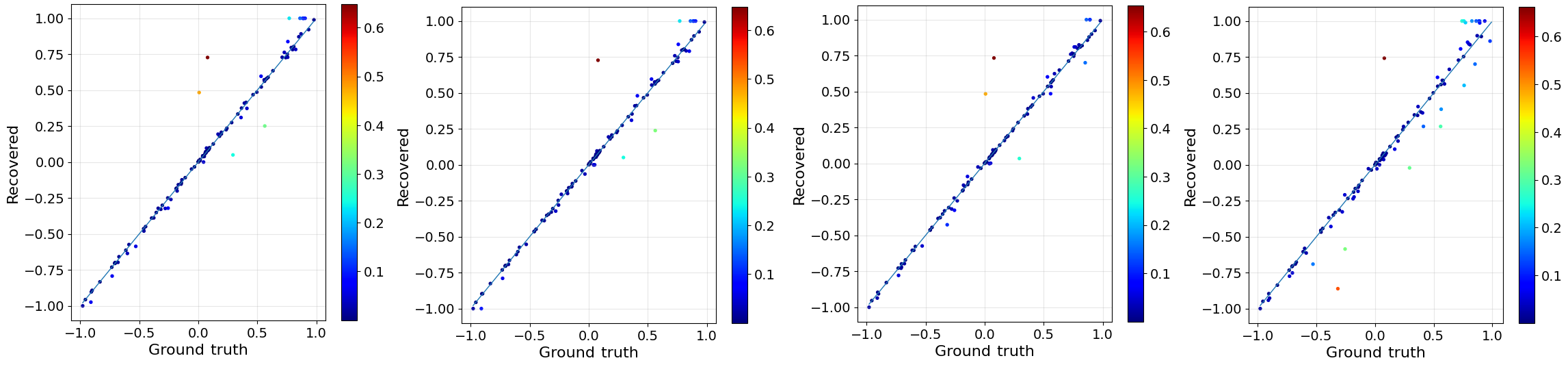}
	\caption{Results of the scalar coefficient recovery problem under varying noise levels added to the initial conditions.
	The noise levels are set as 1\%, 3\%, 10\%, and 30\%, respectively.}
	\label{fig:inv-coef-ic}
\end{figure}

\begin{figure}[htbp]
	\centering
	\includegraphics[width=\linewidth]{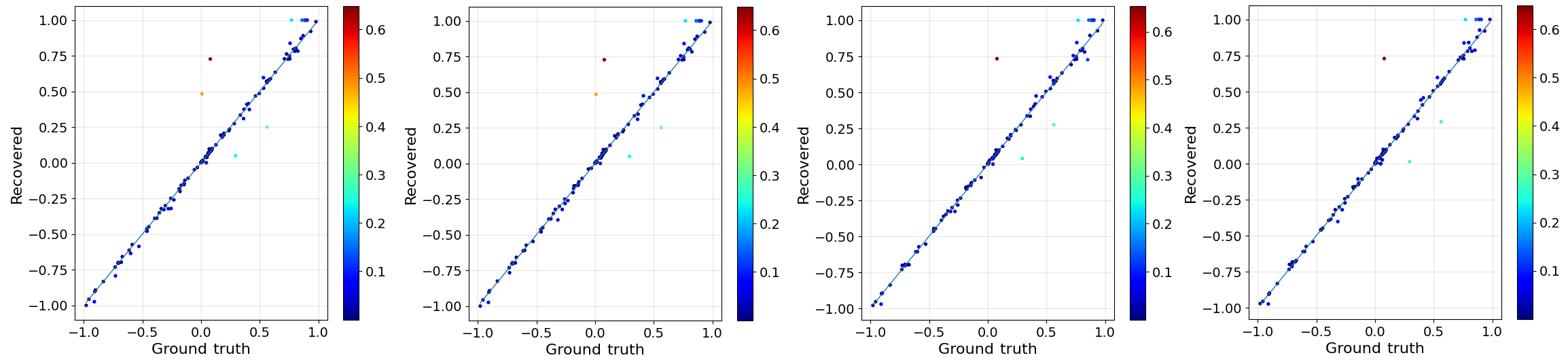}
	\caption{Results of the scalar coefficient recovery problem under varying noise levels added to the observations.
	The noise levels are set as 1\%, 3\%, 10\%, and 30\%, respectively.}
	\label{fig:inv-coef-noise}
\end{figure}

\begin{figure}[htbp]
	\centering
	\includegraphics[width=\linewidth]{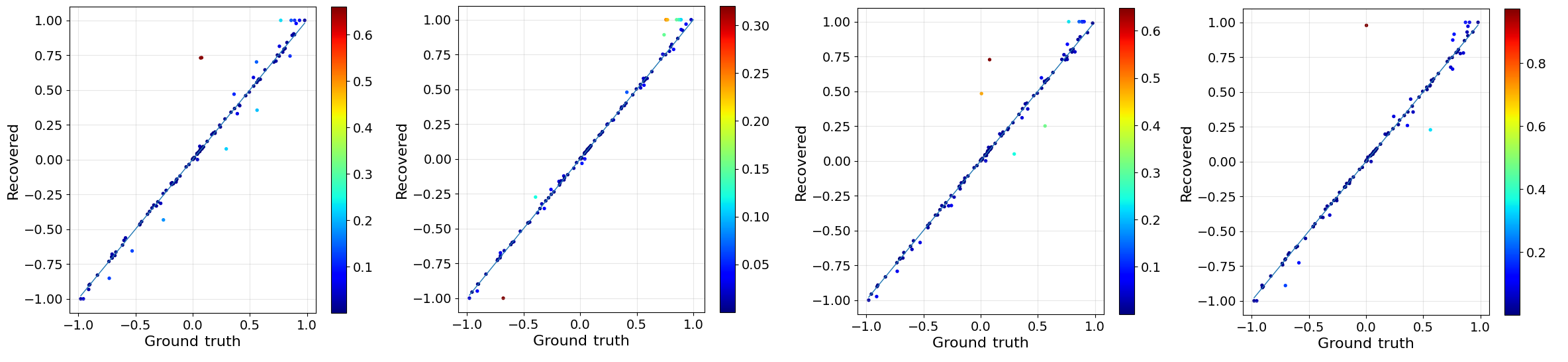}
	\caption{Results of the scalar coefficient recovery problem with varying numbers of observation points.
	The numbers of spatio-temporal points are set as $32 \times 5$, $64 \times 10$, $128 \times 20$, and $256 \times 40$, respectively.}
	\label{fig:inv-coef-sample-points}
\end{figure}

\begin{figure}[htbp]
	\centering
	\includegraphics[width=\linewidth]{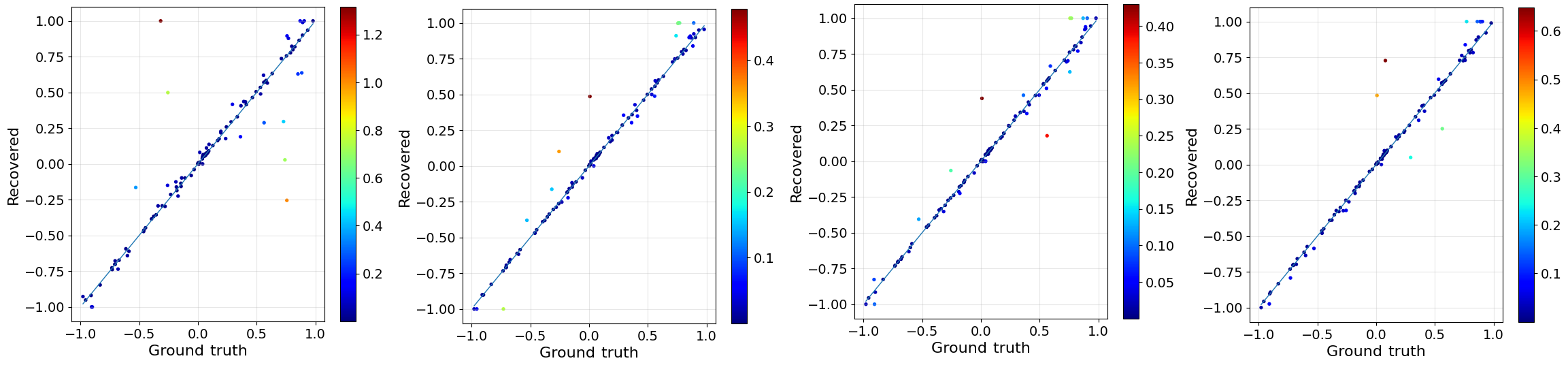}
	\caption{Results of the scalar coefficient recovery problem with varying numbers of observable solutions.
	The numbers of solutions are 1, 5, 10, and 25, respectively.}
	\label{fig:inv-coef-sample}
\end{figure}

We use the PSO algorithm with a swarm size of 100 and a maximum generation of 200 to solve the corresponding optimization problem.
Experiments are conducted on a single NPU (same for all experiments on inverse problems).
We also conduct experiments with different noise levels, the number of spatio-temporal locations, and the number of initial conditions for the PDE, with the results given in Fig.~\ref{fig:inv-coef-ic}--\ref{fig:inv-coef-sample}.
From these experiments, we can conclude that PDEformer-2 is robust and adaptable to noise and various settings in inverse problems.

\label{subsec:inv-coef}

\subsubsection{System identification}
\label{subsec:sys-id}
As a complement to the main text, we present in Table~\ref{tab:sys_id} the complete results of the system identification experiments.
All coefficients were initialized to zero, and the identification was performed using the Particle Swarm Optimization (PSO) algorithm with a swarm size of 150 and a maximum of 150 iterations.
The table compares the ground-truth (GT) coefficients with the recovered values (Rec.), demonstrating the algorithm's ability to approximate a wide range of coefficient configurations across different PDE instances.
\begin{table}
\caption{Results of system identification, including the ground-truth (GT) and recovered coefficients.
We mark the zero coefficients that are mistakenly recovered as non-zero (with absolute value over the threshold $0.1$) in bold.}
\label{tab:sys_id}
\begin{minipage}{\textwidth}
\footnotesize
\begin{tabular}{lccccccccccc}
	\toprule
	& $a$ & $c_{01}$ & $c_{02}$ & $c_{03}$ & $c_{11}$ & $c_{12}$ & $c_{13}$ & $c_{21}$ & $c_{22}$ & $c_{23}$ \\
	\midrule
	GT & 0.007 & -0.376 & 0.000 & \textbf{0.000} & 1.000 & 1.000 & 1.000 & -0.732 & 0.000 & 0.000 \\
	Rec. & 0.007 & -0.536 & 0.087 & \textbf{0.484} & 1.000 & 1.000 & 1.000 & -0.726 & 0.012 & 0.020 \\
	\midrule
	GT & 0.009 & -0.187 & 1.000 & 1.000 & -0.110 & 0.000 & 0.237 & 1.000 & 1.000 & 0.995 \\
	Rec. & 0.009 & -0.172 & 1.000 & 1.000 & -0.103 & 0.018 & 0.245 & 1.000 & 1.000 & 1.000 \\
	\midrule
	GT & 0.005 & 0.000 & 0.884 & 0.000 & 0.181 & 0.000 & 1.000 & 0.008 & 1.000 & 0.879 \\
	Rec. & 0.007 & 0.089 & 1.000 & 0.086 & 0.185 & -0.026 & 1.000 & 0.016 & 1.000 & 0.898 \\
	\midrule
	GT & 0.002 & 1.000 & 0.000 & 0.412 & 0.128 & 0.000 & \textbf{0.000} & 0.061 & 1.000 & -0.005 \\
	Rec. & 0.004 & 1.000 & -0.009 & 1.000 & 0.254 & 0.099 & \textbf{-1.000} & -0.029 & 1.000 & 0.683 \\
	\midrule
	GT & 0.069 & 1.000 & 0.729 & -0.317 & 0.000 & 0.196 & 1.000 & -0.118 & 1.000 & 0.000 \\
	Rec. & 0.083 & 1.000 & 0.760 & -1.000 & -0.019 & 0.213 & 1.000 & -0.129 & 1.000 & 0.060 \\
	\midrule
	GT & 0.001 & \textbf{0.000} & 0.825 & \textbf{0.000} & 0.000 & 1.000 & 1.000 & 1.000 & 1.000 & 0.000 \\
	Rec. & 0.001 & \textbf{-0.174} & 1.000 & \textbf{0.694} & 0.006 & 1.000 & 1.000 & 1.000 & 1.000 & 0.001 \\
	\midrule
	GT & 0.002 & 1.000 & 0.000 & 0.980 & 1.000 & 1.000 & \textbf{0.000} & 0.000 & 0.000 & 1.000 \\
	Rec. & 0.429 & 1.000 & 0.090 & 0.999 & 1.000 & 1.000 & \textbf{1.000} & 0.044 & 0.061 & 0.830 \\
	\midrule
	GT & 0.038 & 0.592 & 0.000 & \textbf{0.000} & 1.000 & 1.000 & 0.000 & 0.000 & -0.169 & 1.000 \\
	Rec. & 0.020 & 1.000 & -0.092 & \textbf{-1.000} & 1.000 & 1.000 & -0.033 & 0.016 & -0.183 & 0.937 \\
	\midrule
	GT & 0.003 & 1.000 & 1.000 & 0.000 & 1.000 & 0.000 & 0.788 & -0.708 & 1.000 & -0.145 \\
	Rec. & 0.003 & 1.000 & 1.000 & 0.093 & 1.000 & -0.013 & 0.791 & -0.699 & 1.000 & -0.248 \\
	\midrule
	GT & 0.083 & 1.000 & 1.000 & 0.759 & 0.000 & 0.000 & -0.978 & 1.000 & -0.183 & 1.000 \\
	Rec. & 0.081 & 1.000 & 1.000 & 1.000 & -0.002 & 0.001 & -0.980 & 1.000 & -0.193 & 1.000 \\
	\bottomrule
	\end{tabular}
\end{minipage}
\end{table}
\subsubsection{Recovering source fields}
\label{subsec:inv-field}
In this task, we still consider equations mentioned in Supplementary Note~\ref{app:inv_dcr}, with the goal of recovering the source term $s(\boldsymbol{r})$.
For the optimization problem~\eqref{eq:inverse-opt}, we set $\lambda = 0.01$, choose the regularization term as $\mathcal{R}(s) = \|\nabla s\|_2$, and use the Adam optimizer~\cite{Adam} with a learning rate of $0.1$ for 150 iterations.
The initial guess of $s(\boldsymbol{r})$ is taken to be a zero field.
To fit into the CNN function encoder of PDEformer-2, we discretize $s(\boldsymbol{r})$ on a grid with resolution $128\times 128$,
and treat the grid values as independent variables during the optimization process.
Spatial derivatives in the regularization term is computed using finite difference.

\subsubsection{Recovering wave velocity fields}
\label{subsec:fwi}
In this task, we consider wave equations mentioned in Supplementary Note~\ref{app:inv_wave}
Still considering the optimization problem~\eqref{eq:inverse-opt}, we take $\lambda = 0.001, \mathcal{R}(a) = \|\nabla a\|_2$ and use the Adam optimizer with a learning rate of 0.02 for $20000$ iterations.
In the experiments, it is found that the optimization may converge to local minima in some cases.
We resolve this issue by running independent optimizations with different initializations, and take the best result as the final prediction.
Specifically, for $k=1,2,\dots,9$, the velocity field is initialized as $a_k^\text{init}(\boldsymbol{r})=0.4k$ (a constant without spatial dependency), and the optimization terminates at $\hat{a}_k$.
The final prediction is taken as $\hat{a}=\hat{a}_K$ with $K=\arg\min_{k=1,\dots,9}\mathcal{L}(a_k)$.

\verB{
\subsection{Model variant with a mesh-free function encoder}\label{app:wdfe_variant}
The main models presented in the main text, namely PDEformer-2-base and PDEformer-2-fast,
both achieve mesh-free output by employing an INR decoder.
However, since both models adopt a CNN function encoder, their input is constrained to a fixed $128\times 128$ grid.
To further enhance the model's flexibility with respect to discretization,
we test a variant that replaces the CNN function encoder with one capable of accepting arbitrary scattered point inputs.
Although this variant achieves comparable inference speed to the CNN-based models,
its na\"ive mesh-free function encoder substantially increases the training cost, making large-scale pretraining and routine experimentation much less efficient.
Consequently, we do not introduce it as a primary model in the main text.
In the following, we present the detailed results of this variant.

\subsubsection{Weighted DeepSet function encoder}
Based on our prior 1D model PDEformer-1~\cite{PDEformer1}, the function $s(\boldsymbol{r})$ is now represented numerically as a series of scattered points $\{(\boldsymbol{r}_j,s(\boldsymbol{r}_j))\mid j\in \mathcal{J}\}$,
and the function encoder performs as
\[\xi_i=\text{Function-Encoder}(\{(\boldsymbol{r}_j,s(\boldsymbol{r}_j))\mid j\in \mathcal{J}\})_k.\]
Inspired by VIDON~\cite{VIDON}, we employ a weighted variant of the DeepSet~\cite{DeepSet} architecture, in the form
\[\begin{split}
	&\qquad\text{Function-Encoder}(\{(\boldsymbol{r}_j,s(\boldsymbol{r}_j))\mid j\in \mathcal{J}\})
	\\&= \psi^3\left(\sum_{j\in \mathcal{J}}\mathrm{softmax}_j\Bigl(\psi^2(\boldsymbol{r}_j,s(\boldsymbol{r}_j))\Bigr)\odot\psi^1(\boldsymbol{r}_j,s(\boldsymbol{r}_j))\right)
	\\&= \psi^3\left(\frac{\sum_{j\in \mathcal{J}}\exp\Bigl(\psi^2(\boldsymbol{r}_j,s(\boldsymbol{r}_j))\Bigr)\odot\psi^1(\boldsymbol{r}_j,s(\boldsymbol{r}_j))}{\sum_{j\in \mathcal{J}}\exp\Bigl(\psi^2(\boldsymbol{r}_j,s(\boldsymbol{r}_j))\Bigr)}\right)
,\end{split}\]
where $\psi^1,\psi^2:\R^2\to\R^{512}$ are represented by Poly-INRs (Fig.~\ref{fig:INR-decoder} without external modulations),
$\psi^3:\R^{512}\to\R^{N\times d_e}$ by a MLP, all with two hidden layers of width $d_\psi=512$.
The exponential and division operations are applied elementwise,
and $\odot$ denotes elementwise multiplication (i.e., Hadamard product) of two vectors.
We take $N=4$ as in the main PDEformer-2 variants with a CNN function encoder.

\subsubsection{Pretraining}
To avoid the prohibitive cost of from-scratch pretraining,
we initialize the model from the pretrained PDEformer-2-base checkpoint,
and replace the CNN function encoder with a randomly initialized weighted DeepSet function encoder.
The total number of model parameters increases marginally from 82.65M to 84.70M.

The pretraining follows the same dataset and training protocol as the two main models,
and proceeds in two stages.
In the first stage, we freeze the backbone (i.e., the graph Transformer and the INR decoder)
and only update the parameters of the newly introduced function encoder,
using a learning rate of $10^{-4}$ for $20$ ``epochs''.
In the second stage, the entire network is jointly finetuned with a learning rate of $2\times 10^{-5}$ for $22$ ``epochs''.
Each training ``epoch'' takes approximately 1 hour and 40 minutes, substantially longer than
the 7 minutes and 40 seconds required by PDEformer-2-base,
due to the higher training cost of the new mesh-free function encoder.
With such a prohibitively high ``per-epoch'' time, we could only afford a limited number of training ``epochs'',
and the pretraining has therefore not yet reached full convergence.
A longer training would likely yield further accuracy improvements,
but for all downstream comparisons presented below, we use the checkpoint obtained at this early stage.

\begin{figure}[htbp]
	\centering
	\includegraphics[width=0.5\textwidth]{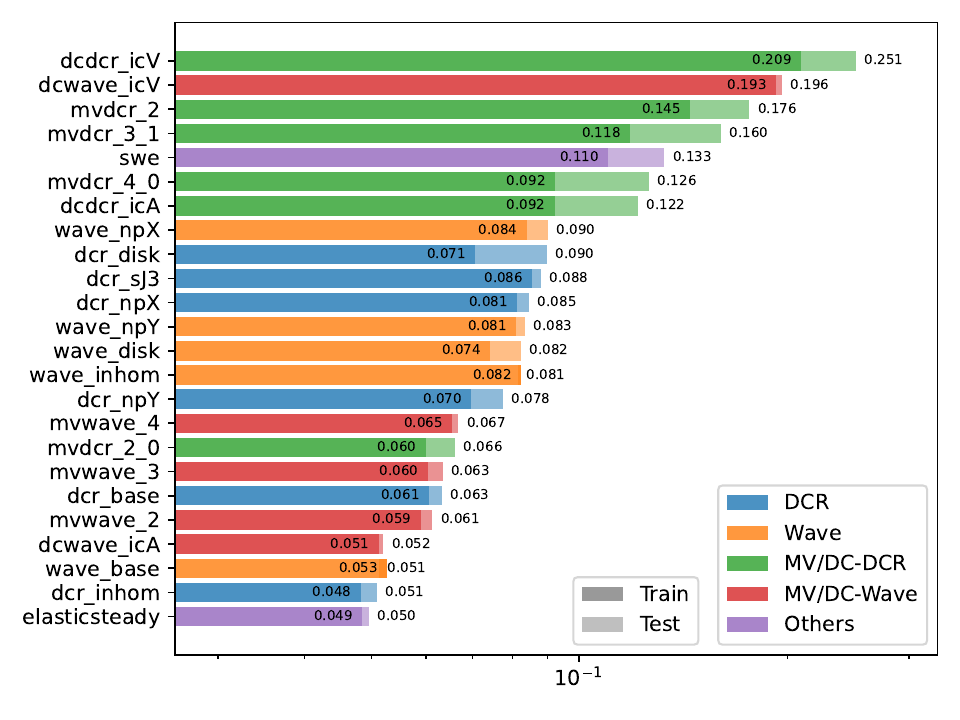}
	\caption{Detailed model accuracy (nRMSE) on the pretraining dataset for the PDEformer-2 variant with a mesh-free weighted DeepSet function encoder.}
	\label{figs:pretrain_err_bars_meshfree}
\end{figure}
After pretraining, we evaluate the model accuracy on the pretraining dataset.
The results are shown in Fig.~\ref{figs:pretrain_err_bars_meshfree}.
Compared with the original models in Fig.~\ref{figs:pretrain_err_bars},
the accuracy of this end-to-end mesh-free variant appears slightly worse than PDEformer-2-base, while slightly better than PDEformer-2-fast.

\subsubsection{Inference time}
Despite the substantially higher training cost, the inference speed of this variant remains comparable to that of PDEformer-2-base.
Specifically, on the Sine-Gordon dataset, the inference time per sample is 680ms for this mesh-free variant,
versus 631ms for PDEformer-2-base.
The fact that the overhead mainly manifests during training, rather than inference,
suggests that the bottleneck is associated with the backward pass of the na\"ively implemented mesh-free function encoder.
Possible causes include the backward computation of our na\"ive implementation interacting unfavorably with the low-level optimizations of the deep learning framework (MindSpore) and the accelerator hardware (Ascend 910A),
or the increased memory footprint during backpropagation that leads to less favorable hardware utilization.
However, we currently lack the relevant expertise to conduct a rigorous attribution study, and leave this investigation for future work.

\subsubsection{Results on specific PDEs}
Since this variant only differs from PDEformer-2-base in its function encoder,
we focus on a single representative case from the specific-PDE experiments in the main text,
namely Wave-C-Sines from~\cite{RIGNO}.
Defined on a disk domain with irregular discretization,
this PDE is the most relevant scenario for evaluating the mesh-free input capability of the new function encoder.

\begin{figure}[htbp]
	\centering
	\includegraphics[width=0.5\textwidth]{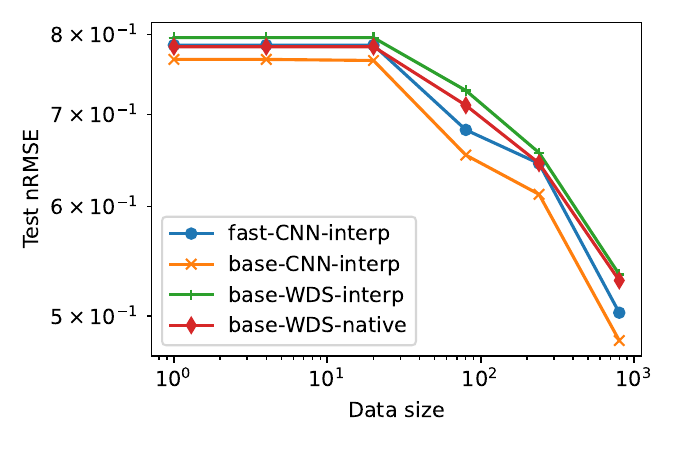}
	\caption{Test nRMSE versus finetuning data size on Wave-C-Sines.
		Model labels follow the convention ``\{INR scale\}-\{function encoder\}-\{IC input mode\}''.
		INR scale: ``base'' (PDEformer-2-base decoder) or ``fast'' (PDEformer-2-fast decoder).
		Function encoder: ``CNN'' or ``WDS'' (weighted DeepSet).
		IC input mode: ``interp'' (initial conditions interpolated to a uniform $128\times 128$ grid before feeding into the network)
		or ``native'' (scattered IC points fed directly without interpolation).
		Thus ``fast-CNN-interp'' corresponds to PDEformer-2-fast in the main text,
		and the two WDS configurations represent the mesh-free variant.
	}
	\label{fig:meshfree_wavec}
\end{figure}

We compare four configurations in Fig.~\ref{fig:meshfree_wavec}:
the two main models, PDEformer-2-base (base-CNN-interp) and PDEformer-2-fast (fast-CNN-interp);
and the mesh-free variant under two input modes, namely with ICs interpolated to a uniform grid (base-WDS-interp)
and with scattered IC points fed directly (base-WDS-native).
Across all training data sizes, the mesh-free variant slightly underperforms PDEformer-2-fast,
and lags further behind PDEformer-2-base.
Notably, feeding the initial conditions directly as scattered points (native)
yields slightly better accuracy than first interpolating them to a uniform grid (interp).

We note that both input modes introduce unfavorable factors for prediction accuracy.
For the native scattered-input mode, the pretraining data (shared with the two main models)
is entirely defined on a uniform $128\times 128$ grid.
Even for data samples on the disk domain, the initial conditions are still specified on the uniform grid,
with only the interior points ultimately affecting the solution during data generation.
Consequently, when the model encounters irregular scattered point inputs during zero-shot inference or finetuning,
this constitutes an out-of-distribution scenario.
For the interpolated mode, the interpolation operation inevitably introduces some degree of information loss.

The current results suggest that the information loss caused by interpolation has a relatively larger negative impact,
while the mismatch in the spatial distribution of input points appears to have a smaller effect.
This may indicate that the weighted DeepSet function encoder architecture
possesses a favorable inductive bias and is relatively insensitive to the discretization of input functions.
We acknowledge, however, that these findings are preliminary,
and more extensive experiments are needed to thoroughly investigate this property.
}

\section{Usage examples for customized PDEs}\label{app:code_example}
In this section, we provide some examples about how to use PDEformer-2 to make predictions.
The following code snippets demonstrate the inference process for solving several customized PDEs using the PDEformer-2 model.
More detailed instructions, the corresponding prediction results, and further examples can be found in the Jupyter Notebook in our code repository: \url{https://github.com/functoreality/pdeformer-2/blob/main/PDEformer_inference.ipynb}.

\subsection{Solving 2D nonlinear conservation law}
Here is an example of using PDEformer-2 to solve a two-dimensional nonlinear conservation law on $[0,1]^2$ with periodic boundaries:
\[u_{t} + (u^2)_x + (-0.3u)_y = 0,\]
\[u(0,x,y) = \exp\left(-2\Bigl(x-\frac{1}{2}\Bigr)^2-8\Bigl(y-\frac{1}{2}\Bigr)^2\right).\]
\begin{lstlisting}[language=Python, caption=Using PDEformer-2 to solve a 2D nonlinear conservation law.]
# define the PDE to be solved
pde = PDENodesCollector()
u = pde.new_uf()
u_ic = np.exp(-4 * (x_fenc - 0.5)**2 - 8 * (y_fenc - 0.5)**2)
pde.set_ic(u, u_ic, x=x_fenc, y=y_fenc)
pde.sum_eq0(pde.dt(u), pde.dx(pde.square(u)), pde.dy(-0.3 * u))

# inference with resolution 32, and plot the results
pde_dag = pde.gen_dag(config)
x_plot, y_plot = np.meshgrid(np.linspace(0, 1, 32), np.linspace(0, 1, 32), indexing="ij")
u_pred = infer_plot_2d(model, pde_dag, x_plot, y_plot)
\end{lstlisting}

\subsection{Convection-diffusion equation}
The following example demonstrates how to solve the convection-diffusion equation:
\[u_{t} + \Bigl(\frac{1}{2}u\Bigr)_x + f(x,y) - 3\times 10^{-3}\Delta u = 0,\]
\[u(0,x,y) = \sin(2\pi x)\cos(4\pi y),\]
\[f(x,y) = \exp\left(-32\Bigl(x-\frac{1}{2}\Bigr)^2-32\Bigl(y-\frac{1}{2}\Bigr)^2\right).\]
\begin{lstlisting}[language=Python, caption=Solving the convection-diffusion equation.]
u_ic = np.sin(2 * np.pi * x_fenc) * np.cos(4 * np.pi * y_fenc)
f_arr = np.exp(-32 * (x_fenc - 0.5)**2 - 32* (y_fenc - 0.5)**2)
# plot source term
plt.figure(figsize=(3, 3))
plt.pcolormesh(x_fenc, y_fenc, f_arr, cmap="turbo")

# define PDE
pde = PDENodesCollector()
u = pde.new_uf()
pde.set_ic(u, u_ic, x=x_fenc, y=y_fenc)
f = pde.new_coef_field(f_arr, x=x_fenc, y=y_fenc)
pde.sum_eq0(u.dt, (0.5 * u).dx, f, -(3e-3 * (u.dx.dx + u.dy.dy)))

# inference and plot
pde_dag = pde.gen_dag(config)
pde_dag.plot(hide='aux')  # plot the computational graph
u_pred = infer_plot_2d(model, pde_dag, x_plot, y_plot)
\end{lstlisting}

\subsection{Conservative variant of 2D viscous Burgers' equation}
In this part, we illustrate how to solve the conservative variant of the two-dimensional viscous Burgers' equation:
\[u_t + (u^2)_x + (uv)_y - 10^{-3}\Delta u = 0,\]
\[v_t + (uv)_x + (v^2)_y - 10^{-3}\Delta v = 0,\]
where $u(0, x, y)$ and $v(0, x, y)$ are independent gaussian random fields.
\begin{lstlisting}[language=Python, caption=Solving the conservative variant of 2D viscous Burgers' equation.]
u_ic, v_ic = sample_grf(2)
pde = PDENodesCollector()
# add unknown fields
u = pde.new_uf()
v = pde.new_uf()
# set initial conditions
pde.set_ic(u, u_ic, x=x_fenc, y=y_fenc)
pde.set_ic(v, v_ic, x=x_fenc, y=y_fenc)
# define the PDEs
pde.sum_eq0(u.dt, (u * u).dx, (u * v).dy, -(1e-3 * (u.dx.dx + u.dy.dy)))
pde.sum_eq0(v.dt, (u * v).dx, (v * v).dy, -(1e-3 * (v.dx.dx + v.dy.dy)))
pde_dag = pde.gen_dag(config)
# inference and plot
u_pred = infer_plot_2d(model, pde_dag, x_plot, y_plot, snap_t=np.linspace(0, 1, 4))
\end{lstlisting}

\subsection{Wave equation on a disk domain with Mur boundary}
Finally, we present an example of solving the wave equation, concentrating our focus to a disk region $\Omega$ with radius $r_0=0.4$ and center $(x_0,y_0)=(0.5,0.6)$.
The wave can propagate outside the disk, so the absorbing boundary condition needs to be specified at the boundary circle, for which we select the Mur boundary condition:
\[u_{tt} - 0.7^2\Delta u = 0,\]
\[u(0,x,y) = \sin(2\pi x)\cos(4\pi y), \quad u_t(0,x,y) = 0,\]
\[\left(u_t + c \frac{\partial}{\partial n}u \right)\bigg|_{\partial \Omega} = 0.\]
\begin{lstlisting}[language=Python, caption=Solving the wave equation on a disk with absorbing Mur boundary condition.]
# define disk-shaped domain with polar coordinates
x0, y0, r0 = 0.5, 0.6, 0.4
disk_sdf = np.sqrt((x_fenc - x0)**2 + (y_fenc - y0)**2) - r0
r = r0 * np.linspace(0, 1, 16)**0.75
theta = np.linspace(-np.pi, np.pi, 64)
r, theta = np.meshgrid(r, theta, indexing="ij")
x_disk = x0 + r * np.cos(theta)
y_disk = y0 + r * np.sin(theta)

# define PDE and inference
pde = PDENodesCollector(dim=2)
domain = pde.new_domain(disk_sdf, x=x_fenc, y=y_fenc)
u = pde.new_uf(domain)
pde.set_ic(u, u_ic, x=x_fenc, y=y_fenc)
pde.set_ic(u.dt, 0, x=x_fenc, y=y_fenc)
c = 0.7
pde.sum_eq0(u.dt.dt, -(c**2 * (u.dx.dx + u.dy.dy)))
boundary = pde.new_domain(np.abs(disk_sdf), x=x_fenc, y=y_fenc)
pde.bc_sum_eq0(boundary, [u.dt] + pde.dn_sum_list(u, domain, coef=c))
pde_dag = pde.gen_dag(config)
_ = infer_plot_2d(model, pde_dag, x_disk, y_disk)
\end{lstlisting}

\section{Visualizations}\label{app:visual}
In this section, we present additional visualizations of PDEformer-2-base predictions across different datasets, including those given in Supplementary Note~\ref{app:dataset} and not shown in the main text.
These include predicted solutions, ground-truth, and the corresponding error maps.
For time-dependent equations, snapshots are shown at $t = 0$, $0.33$, $0.67$, and $1.0$ to capture the temporal evolution.
These visualizations highlight PDEformer-2's ability to accurately model the underlying dynamics of the equations.

\subsection{DCR equation}
\begin{figure}
	\centering
	\includegraphics[width=\textwidth]{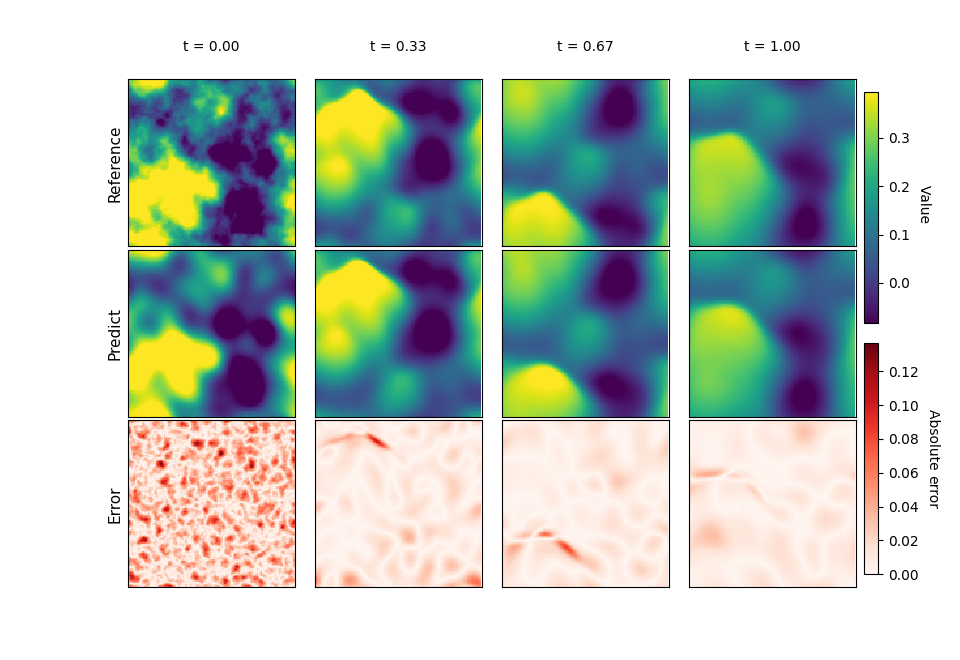}
	\begin{equation*}
		u_t - \tilde{a}\nabla\cdot(\tilde{a}\nabla u) + u^3 + (c_{21}u + u^3)_y = 0
	\end{equation*}
	\caption{Visualization of the DCR Equation on a square domain with periodic boundary conditions along both axes.}
	\label{fig:dcr_0}
\end{figure}

\begin{figure}
	\centering
	\includegraphics[width=\textwidth]{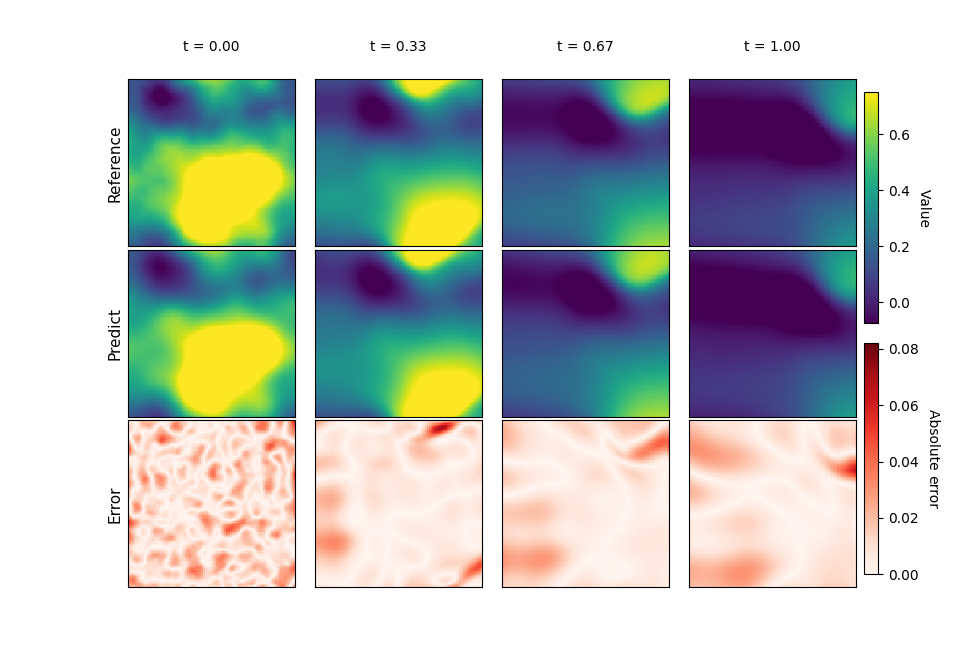}
	\begin{equation*}
		\begin{cases}
			u_t - a\Delta u + s + c_{01}u + c_{02}u^2 + (c_{11}u)_x + (c_{21}u + c_{23}u^3)_y = 0, \\
			\partial_n u \big|_{\Gamma_R} = 0, \\
			(\partial_n u + a_L(r)u)\big|_{\Gamma_L} = 0.
		\end{cases}
	\end{equation*}
	\caption{Visualization of the DCR equation on a square domain periodic along $y$-axis.}
	\label{fig:dcr_npX}
\end{figure}

\begin{figure}
	\centering
	\includegraphics[width=\textwidth]{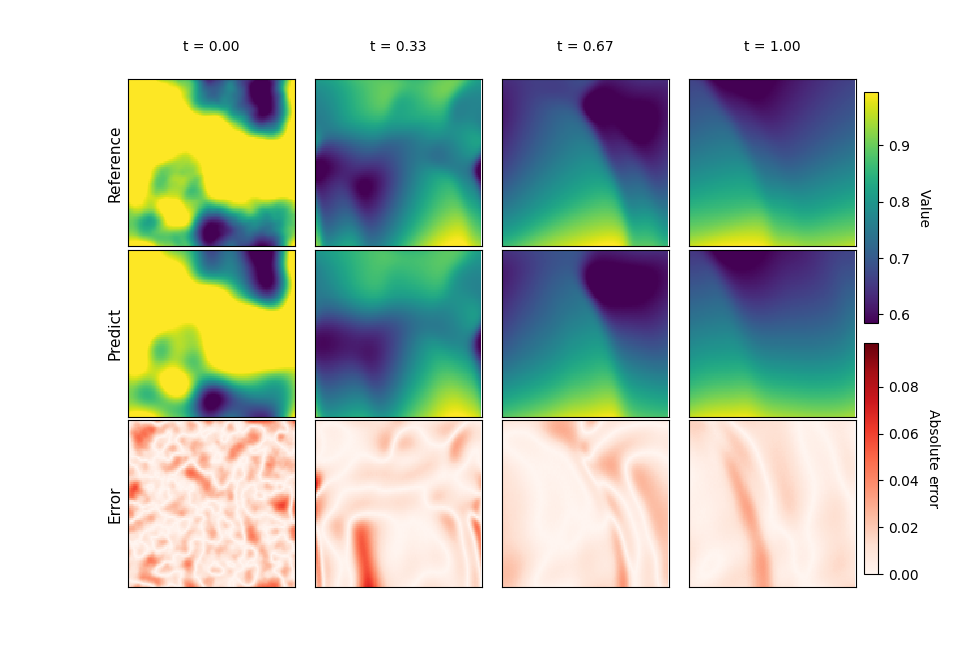}
	\begin{equation*}
		\begin{cases}
			u_t - \tilde{a}\nabla\cdot(\tilde{a}\nabla u) + c_{02}u^2 + u^3 + (u + c_{13}u^3)_x + (c_{21}u + u^2 + c_{23}u^3)_y = 0, \\
			\partial_n u \big|_{\Gamma_U} = 0, \\
			(\partial_n u + a_D(r)u + c_D)\big|_{\Gamma_D} = 0
		\end{cases}
	\end{equation*}
	\caption{Visualization of the DCR equation on a square domain periodic along $x$-axis.}
	\label{fig:dcr_npY}
\end{figure}

\begin{figure}
	\centering
	\includegraphics[width=\textwidth]{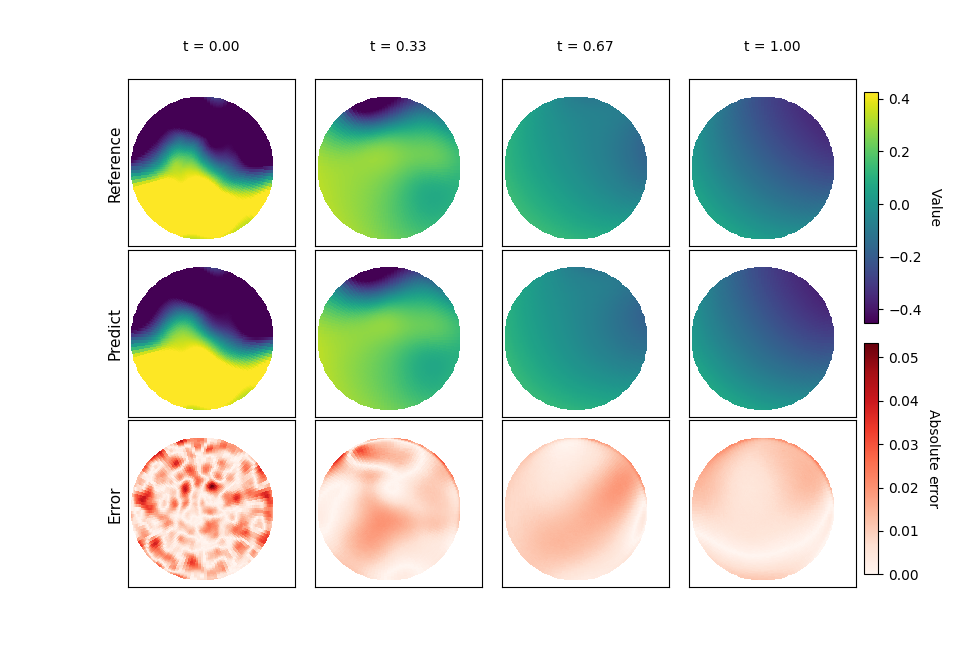}
	\begin{equation*}
		\begin{cases}
			u_t - \tilde{a}\nabla\cdot(\tilde{a}\nabla u) + 1 + c_{01}u + u^2 + u^3 + (u + c_{13}u^3)_x + u_y = 0, \\
			(u + \partial_n u + c_O)\big|_{\Gamma_O} = 0
		\end{cases}
	\end{equation*}
	\caption{Visualization of the DCR equation on a disk domain.}
	\label{fig:dcr_disk}
\end{figure}

We consider the DCR equation on different domains, including a square domain periodic along both axes, periodic only along $y$-axis, periodic only along $x$-axis, as well as a disk-shaped domain.
The visualizations for these cases are shown in Fig.~\ref{fig:dcr_0}, \ref{fig:dcr_npX}, \ref{fig:dcr_npY}, and~\ref{fig:dcr_disk} respectively.

\subsection{Wave equation}
\begin{figure}
	\centering
	\includegraphics[width=\textwidth]{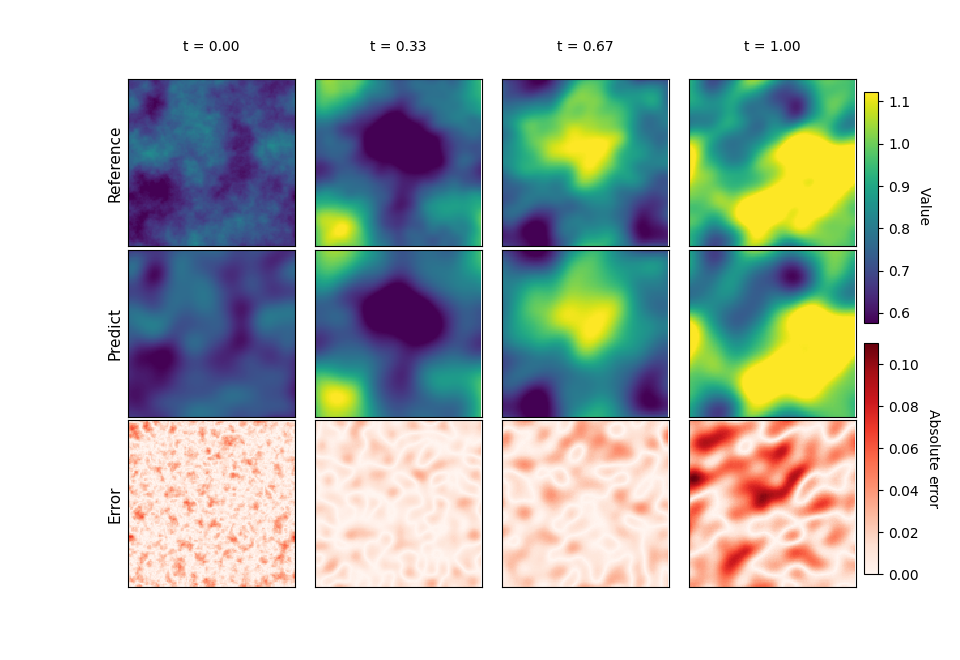}
	\begin{equation*}
		u_{tt}-\Delta u+c_{02}u^2+(c_{11}u+c_{12}u^2+c_{13}u^3)_x+(c_{22}u^2+c_{23}u^3)_y=0
	\end{equation*}
	\caption{Visualization of the wave equation on a square domain with periodic boundary conditions along both axes.}
	\label{fig:wave_0}
\end{figure}

\begin{figure}
	\centering
	\includegraphics[width=\textwidth]{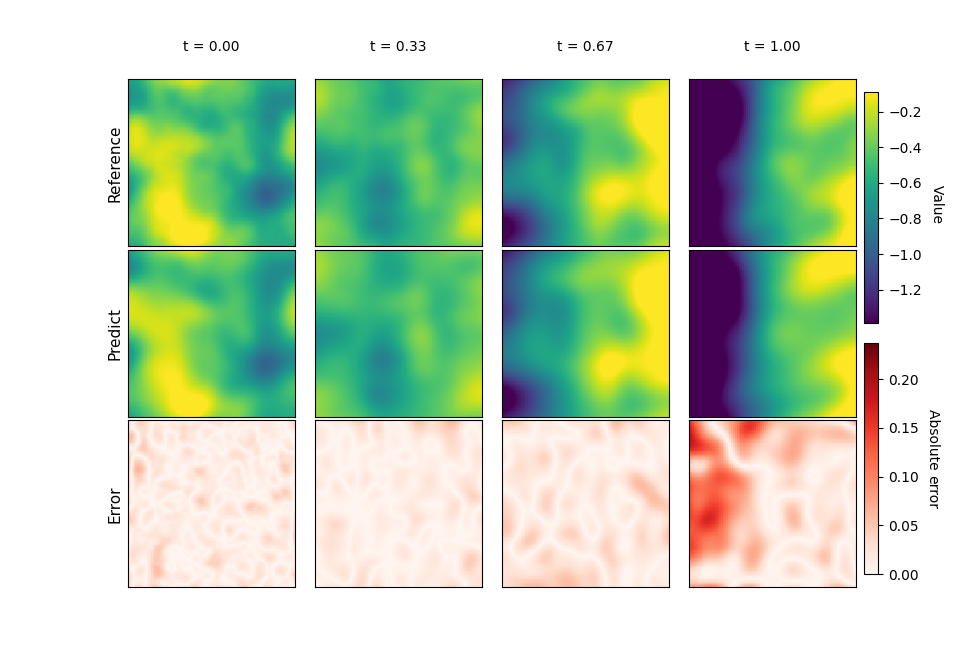}
	\begin{equation*}
		\begin{cases}
			u_{tt}-a\Delta u+c_{001}u+c_{002}u^2+c_{003}u^3+(u^3)_x\\
			+(c_{201}u+u^2+c_{210}\cos(c_{211}u+c_{212}u^2))_y=0 \\
			(\partial_nu+u)|_{\Gamma_R}=0,\\
			(u+b_L\partial_nu+c_L)|_{\Gamma_L}=0
		\end{cases}
	\end{equation*}
	\caption{Visualization of the wave equation on a square domain periodic along $y$-axis.}
	\label{fig:wave_npX}
\end{figure}

\begin{figure}
	\centering
	\includegraphics[width=\textwidth]{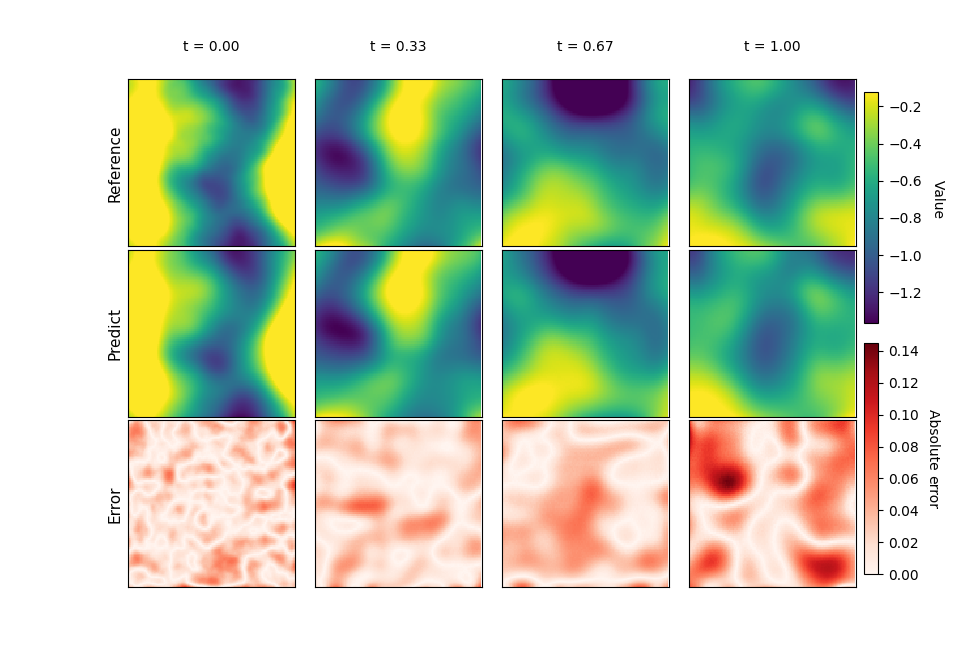}
	\begin{equation*}
		\begin{cases}
			u_{tt}+\mu(r)u_t-a\Delta u+s(r)+c_{001}u+u^2+(u+c_{102}u^2+c_{103}u^3)_x\\
			+(u^2+\sin(c_{211}u+c_{212}u^2)+\cos(u^2))_y=0 \\
			(\partial_nu+a_Uu+c_U)|_{\Gamma_U}=0,\\
			(\partial_tu+a_Du+b_D\partial_nu)|_{\Gamma_D}=0
		\end{cases}
	\end{equation*}
	\caption{Visualization of the wave equation on a square domain periodic along $x$-axis.}
	\label{fig:wave_npY}
\end{figure}

\begin{figure}
	\centering
	\includegraphics[width=\textwidth]{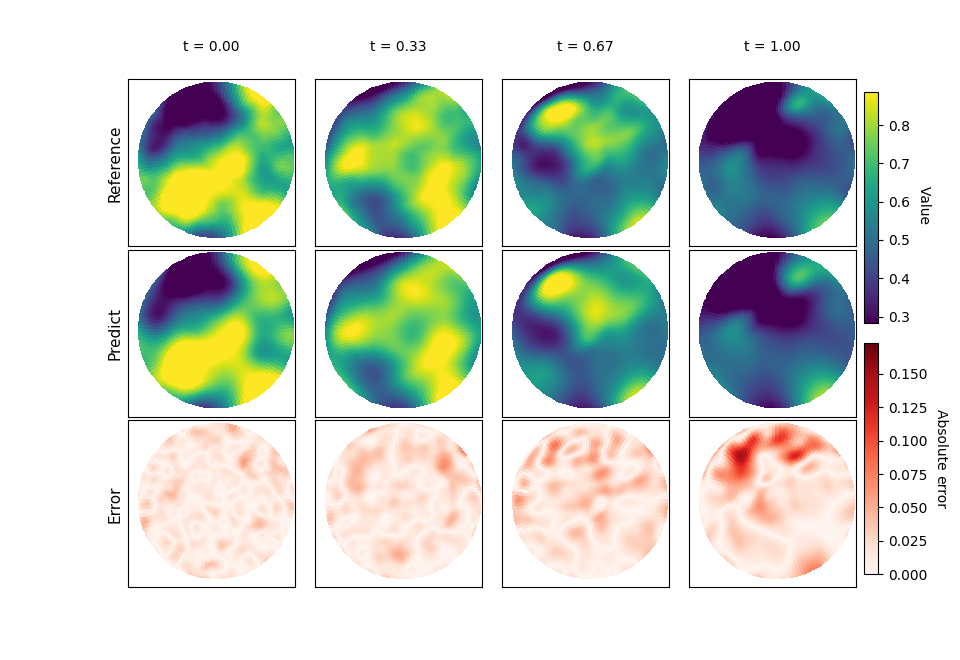}
	\begin{equation*}
		\begin{cases}
			u_{tt}+u_t-\tilde{a}\nabla\cdot(\tilde{a}\nabla u)+u^3+(c_{11}u+c_{12}u^2+c_{13}u^3)_x+(c_{22}u^2+u^3)_y=0\\
			(\partial_tu+a_Ou)|_{\Gamma_O}=0
		\end{cases}
	\end{equation*}
	\caption{Visualization of the wave equation on a disk domain.}
	\label{fig:wave_disk}
\end{figure}

We consider the wave equation on a square domain with periodic boundary conditions along both axes, periodic along $y$-axis, periodic along $x$-axis, and a disk domain with randomly generated boundary conditions.
The visualizations for these cases are shown in Fig.~\ref{fig:wave_0}, \ref{fig:wave_npX}, \ref{fig:wave_npY}, and~\ref{fig:wave_disk} respectively.
\subsection{Multi-variable DCR equation}
\begin{figure}
	\centering
	\includegraphics[width=0.8\textwidth]{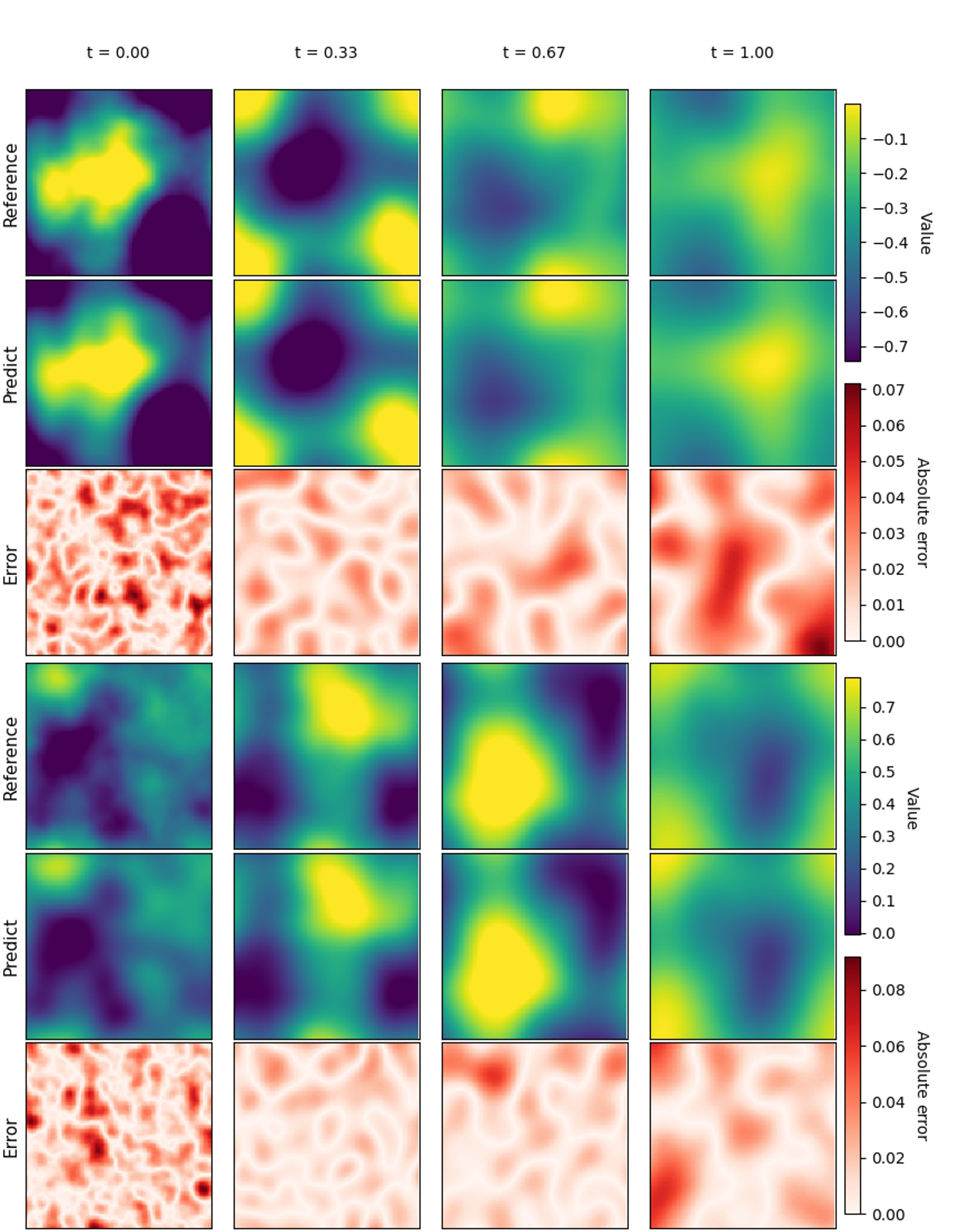}
	\begin{equation*}
		\begin{cases}
			\partial_tu_0+c_{001}u_1+u_0+\partial_x(u_0+u_1)+\partial_y(u_1)-a_0\Delta u_0=0, \\
			\partial_tu_1+u_1+u_0+\partial_x(c_{111}u_1+u_0)+\partial_y(u_0)+s_1(r)-\nabla\cdot(a_1\nabla u_1)=0
		\end{cases}
	\end{equation*}
	\caption{Visualization of the multi-variable DCR equation on a square domain that is periodic along both axes, including all solution components $u_0,u_1$.}
	\label{fig:mcdcr}
\end{figure}

We consider the multi-variable DCR equation on a square domain that is periodic along both axes.
Visualization of the variables $u_0,u_1$ is shown in Fig.~\ref{fig:mcdcr}.
\subsection{Divergence-constrained DCR equation}
\begin{figure}
	\centering
	\includegraphics[width=0.6\textwidth]{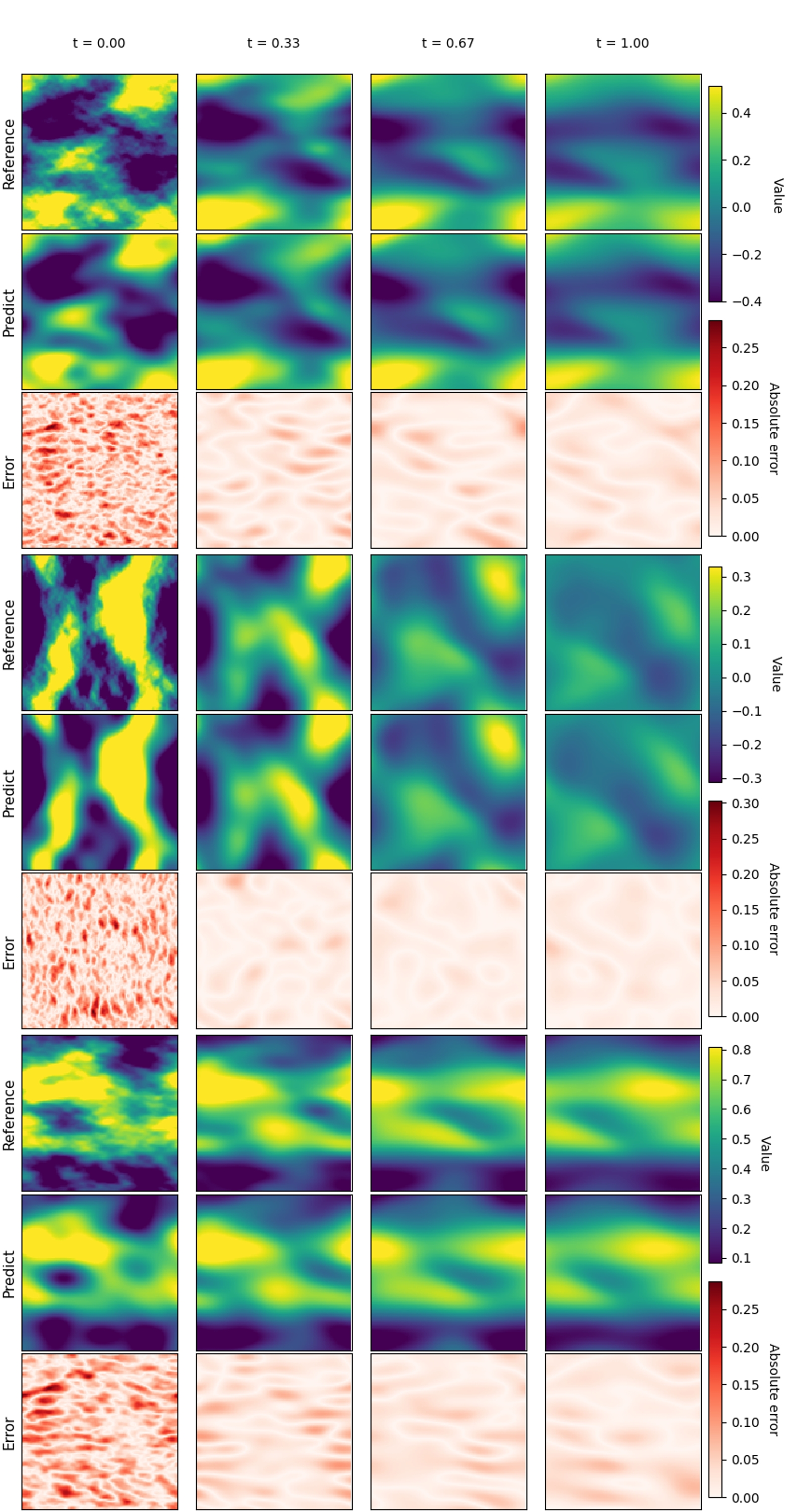}
	\begin{equation*}
		\begin{cases}
			\partial_tu_0+u_1+\partial_x(c_{100}u_0+c_{101}u_1)+\partial_y(c_{201}u_1)\\
			-\tilde{a}_0\nabla\cdot(\tilde{a}_0\nabla u_0)+\partial_xp=0, \\
			\partial_tu_1+u_1+c_{010}u_0+\partial_x(c_{111}u_1+c_{110}u_0)+\partial_y(u_0+b_{2101}u_0u_1)\\
			+s_1(r)-\nabla\cdot(a_1\nabla u_1)+\partial_yp+(-c_1)p=0, \\
			\partial_xu_0+\partial_yu_1+c_1u_1=0
		\end{cases}
	\end{equation*}
	\caption{Visualization of the divergence-constrained DCR equation on a square domain periodic along both axes and initial conditions satisfying the divergence constraint.}
	\label{fig:dcdcr_icV}
\end{figure}

\begin{figure}
	\centering
	\includegraphics[width=0.6\textwidth]{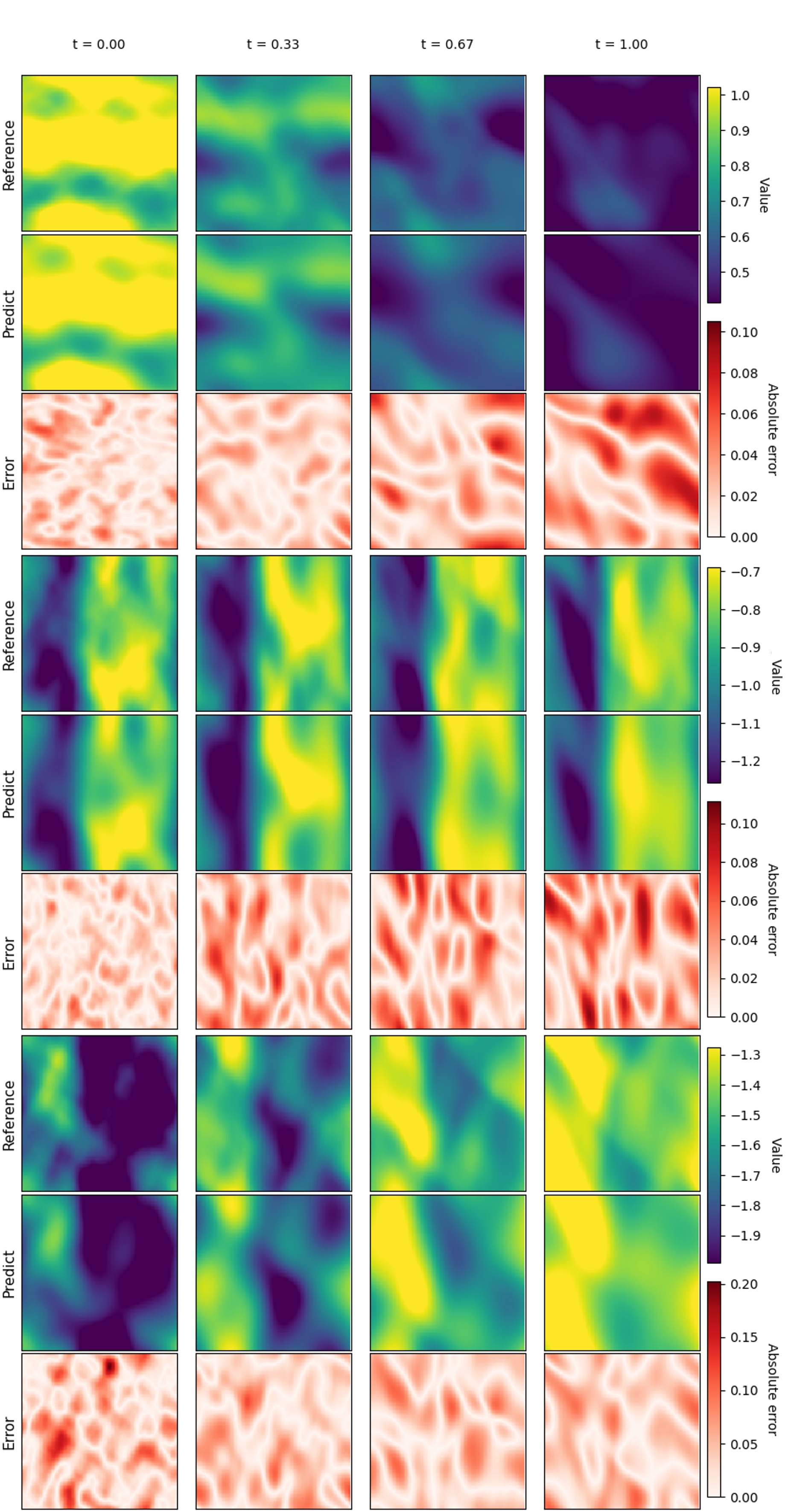}
	\begin{equation*}
		\begin{cases}
			\partial_tu_0+u_0+c_{001}u_1+\partial_x(u_0+u_0u_1+u_0^2)+\partial_y(c_{200}u_0+u_1^2)\\
			-\tilde{a}_0\nabla\cdot(\tilde{a}_0\nabla u_0)+\partial_xp=0, \\
			\partial_tu_1+c_{010}u_0+c_{011}u_1+\partial_x(u_0^2)+\partial_y(u_0+u_0u_1)+s_1+\partial_yp-p=0, \\
			\partial_xu_0+\partial_yu_1+u_1+c_2=0
		\end{cases}
	\end{equation*}
	\caption{Visualization of the divergence-constrained DCR equation on a square domain periodic along both axes and initial conditions violating the divergence constraint.}
	\label{fig:dcdcr_icA}
\end{figure}

We consider the divergence-constrained DCR equation on a square domain periodic along both axes, including those with initial conditions satisfying and violating the divergence constraint.
Numerical results for valid (satisfying the divergence constraint) initial conditions are shown in Fig.~\ref{fig:dcdcr_icV},
while cases with arbitrary (violating the divergence constraint) initial conditions are visualized in Fig.~\ref{fig:dcdcr_icA}.
\subsection{Multi-variable wave equation}
\begin{figure}
	\centering
	\includegraphics[width=0.8\textwidth]{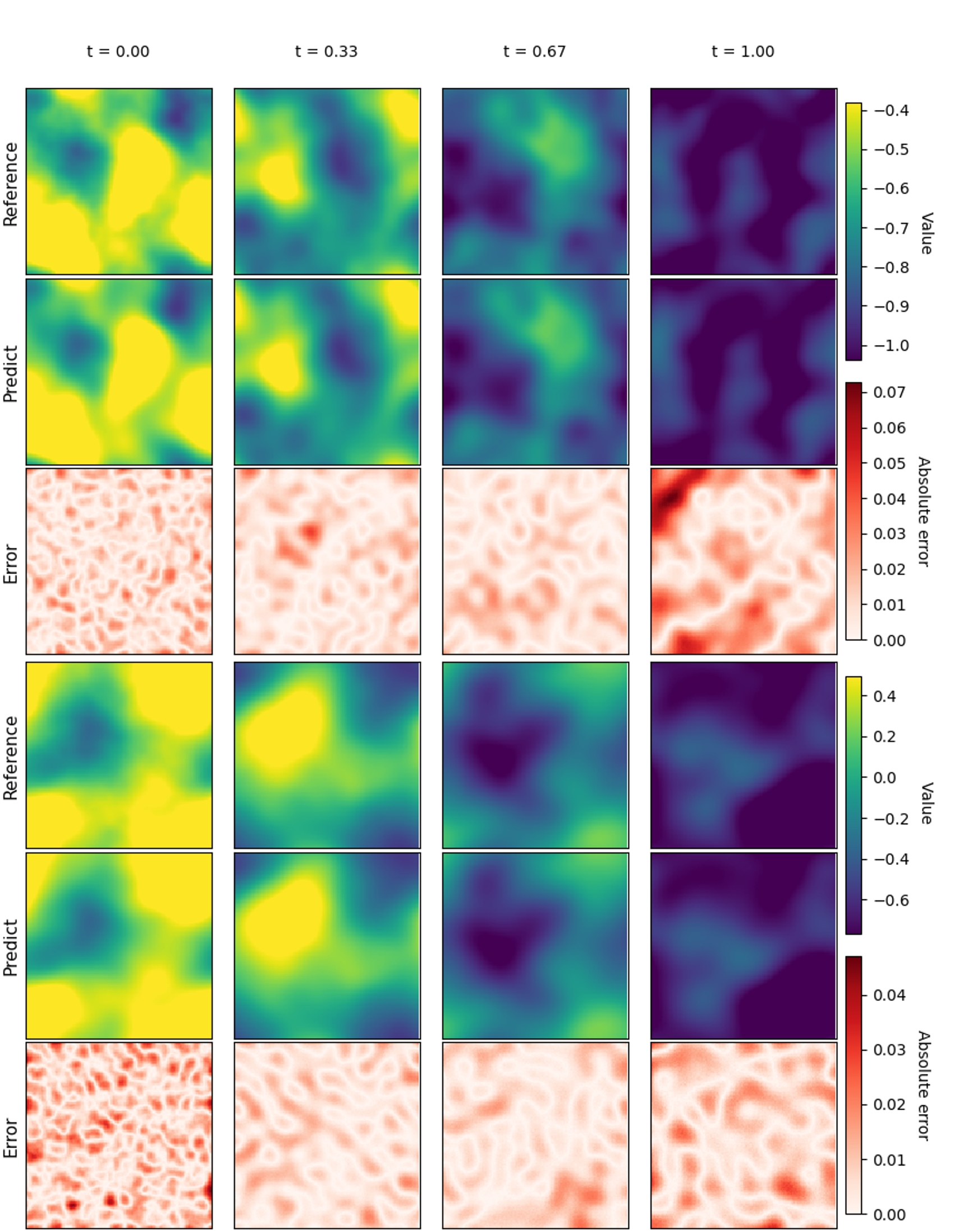}
	\begin{equation*}
		\begin{cases}
			\partial_{tt}u_0+\mu_0\partial_tu_0+u_1+c_{000}u_0+b_{0001}u_0u_1+\partial_x(u_0+u_1)+\partial_y(c_{200}u_0+u_0u_1)\\
			+s_0-a_0\Delta u_0=0, \\
			\partial_{tt}u_1+\mu_1\partial_tu_1+u_1+c_{010}u_0+\partial_y(b_{2100}u_0^2)-\tilde{a}_1\nabla\cdot(\tilde{a}_1\nabla u_1)=0
		\end{cases}
	\end{equation*}
	\caption{Visualization of the multi-variable wave equation on a square domain periodic along both axes, including all solution components $u_0,u_1$.}
	\label{fig:mcwave}
\end{figure}

We consider the multi-variable wave equation on a square domain periodic along both axes and only two variables for convenience.
Visualization of the solution components $u_0,u_1$ is shown in Fig.~\ref{fig:mcwave}.
\subsection{Divergence-constrained wave equation}
\begin{figure}
	\centering
	\includegraphics[width=0.6\textwidth]{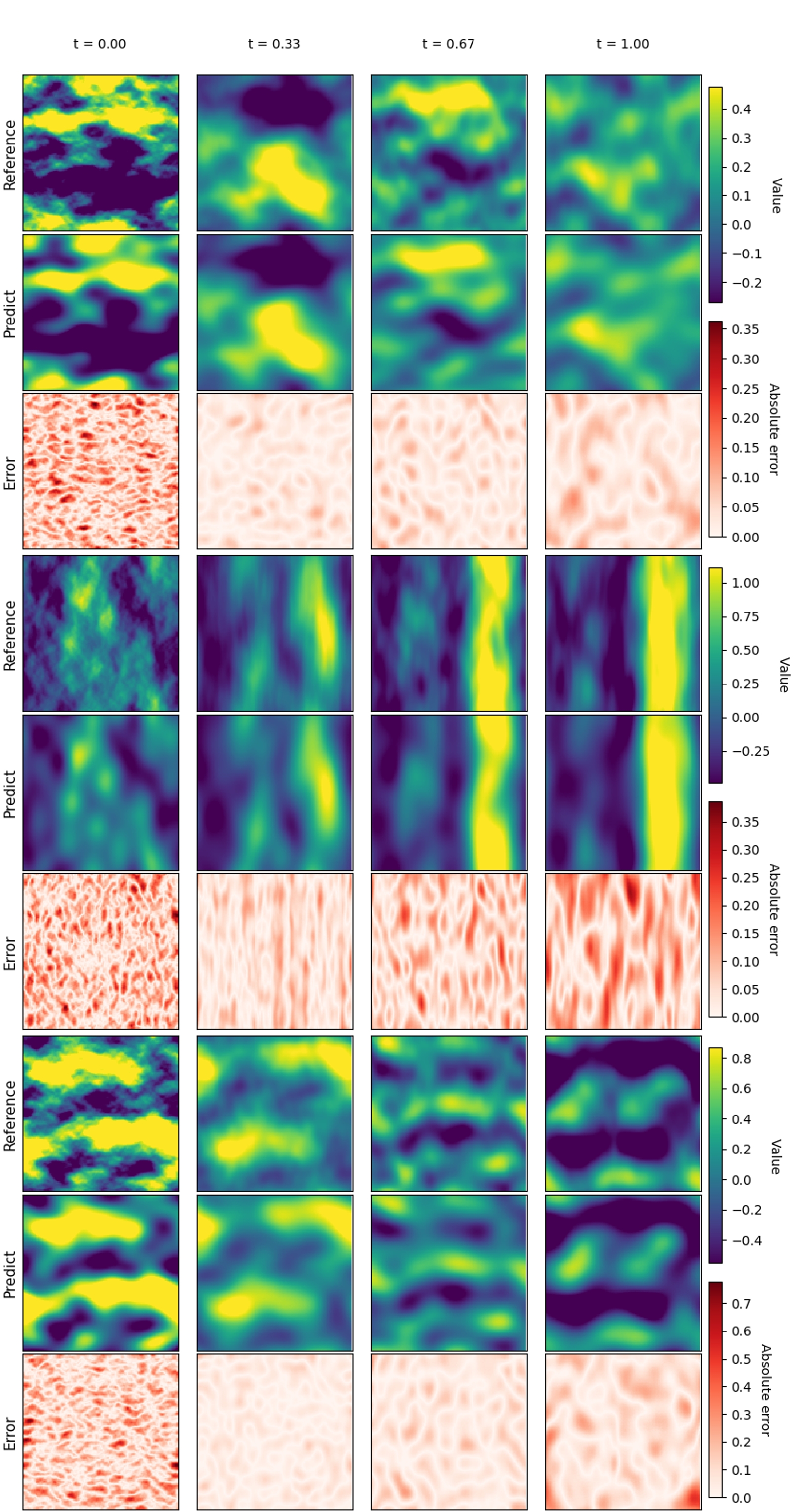}
	\begin{equation*}
		\begin{cases}
			\partial_{tt}u_0+\mu_0(r)\partial_tu_0+u_0+\partial_x(u_0+u_1^2+u_0^2)+\partial_y(u_0+u_1)\\
			-\tilde{a}_0\nabla\cdot(\tilde{a}_0\nabla u_0)+\partial_xp+(-c_0)p=0, \\
			\partial_{tt}u_1+\mu_1\partial_tu_1+u_1^2+\partial_x(b_{1101}u_0u_1+b_{1100}u_0^2)+\partial_y(u_0+c_{211}u_1)\\
			-a_1\Delta u_1+\partial_yp+(-c_1)p=0, \\
			\partial_xu_0+\partial_yu_1+c_0u_0+c_1u_1=0
		\end{cases}
	\end{equation*}
	\caption{Visualization of the divergence-constrained wave equation on a square domain periodic along both axes and initial conditions satisfying the divergence constraint.}
	\label{fig:dcwave_icV}
\end{figure}

\begin{figure}
	\centering
	\includegraphics[width=0.6\textwidth]{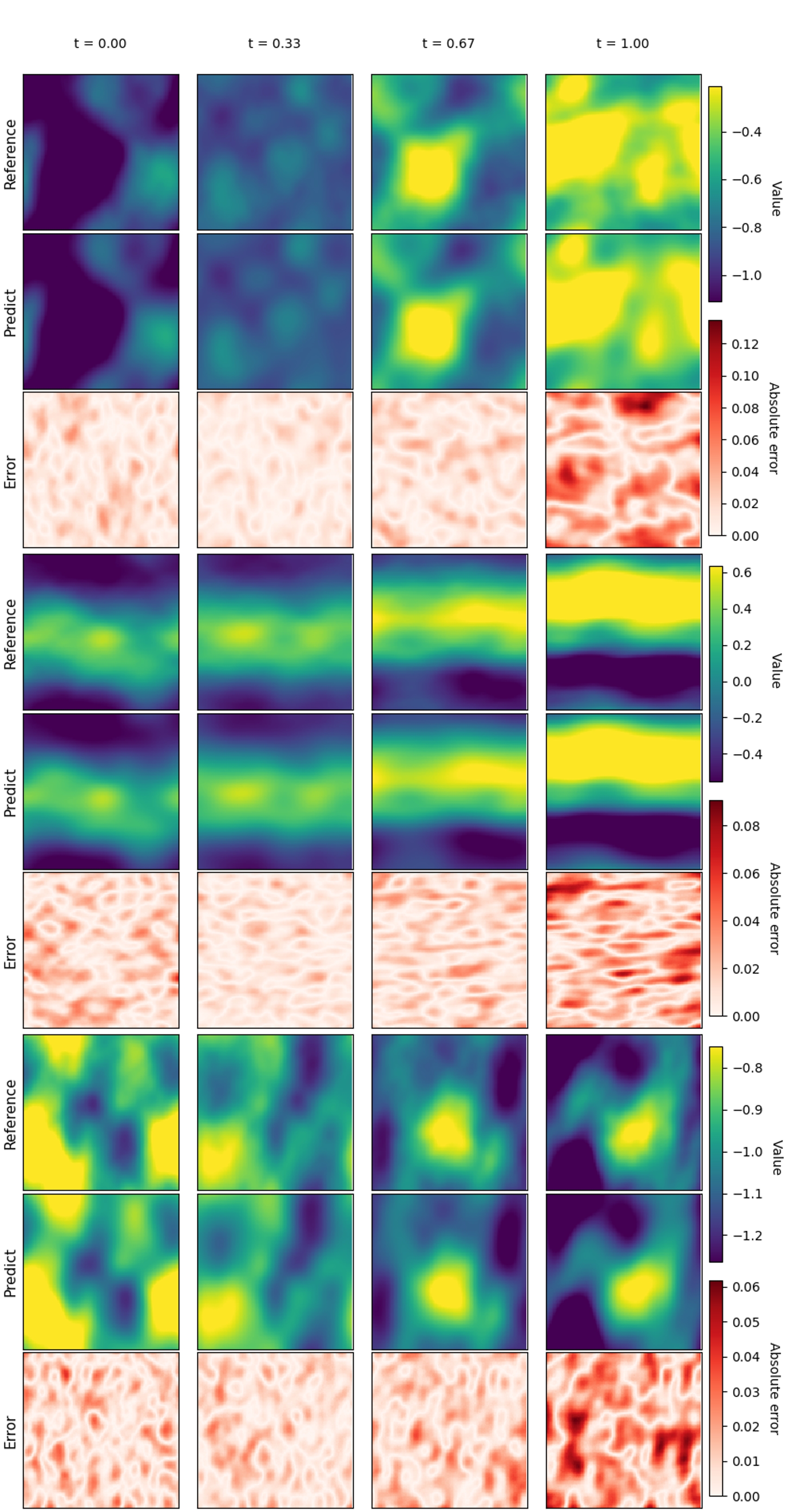}
	\begin{equation*}
		\begin{cases}
			\partial_{tt}u_0+\mu_0(r)\partial_tu_0+u_0+b_{0001}u_0u_1+b_{0011}u_1^2+\partial_x(c_{100}u_0+u_1)\\
			+\partial_y(u_0+c_{201}u_1+b_{2001}u_0u_1)
			+s_0-\nabla\cdot(a_0\nabla u_0)+\partial_xp+(-c_0)p=0\\
			\partial_{tt}u_1+c_{010}u_0+\partial_x(c_{111}u_1+c_{110}u_0)+\partial_y(c_{210}u_0+b_{2101}u_0u_1)\\
			-\nabla\cdot(a_1\nabla u_1)+\partial_yp+(-c_1)p=0\\
			\partial_xu_0+\partial_yu_1+c_0u_0+c_1u_1+1=0
		\end{cases}
	\end{equation*}
	\caption{Visualization of the divergence-constrained wave equation on a square domain periodic along both axes and initial conditions violating the divergence constraint.}
	\label{fig:dcwave_icA}
\end{figure}

We consider the divergence-constrained wave equation on a square domain periodic along both axes, including those with initial conditions satisfying and violating the divergence constraint.
Numerical results for valid initial conditions are shown in Fig.~\ref{fig:dcwave_icV}, while cases with arbitrary initial conditions are visualized in Fig.~\ref{fig:dcwave_icA}.
\subsection{Generalized shallow water equation}
\begin{figure}[htbp]
	\centering
	\includegraphics[width=0.6\textwidth]{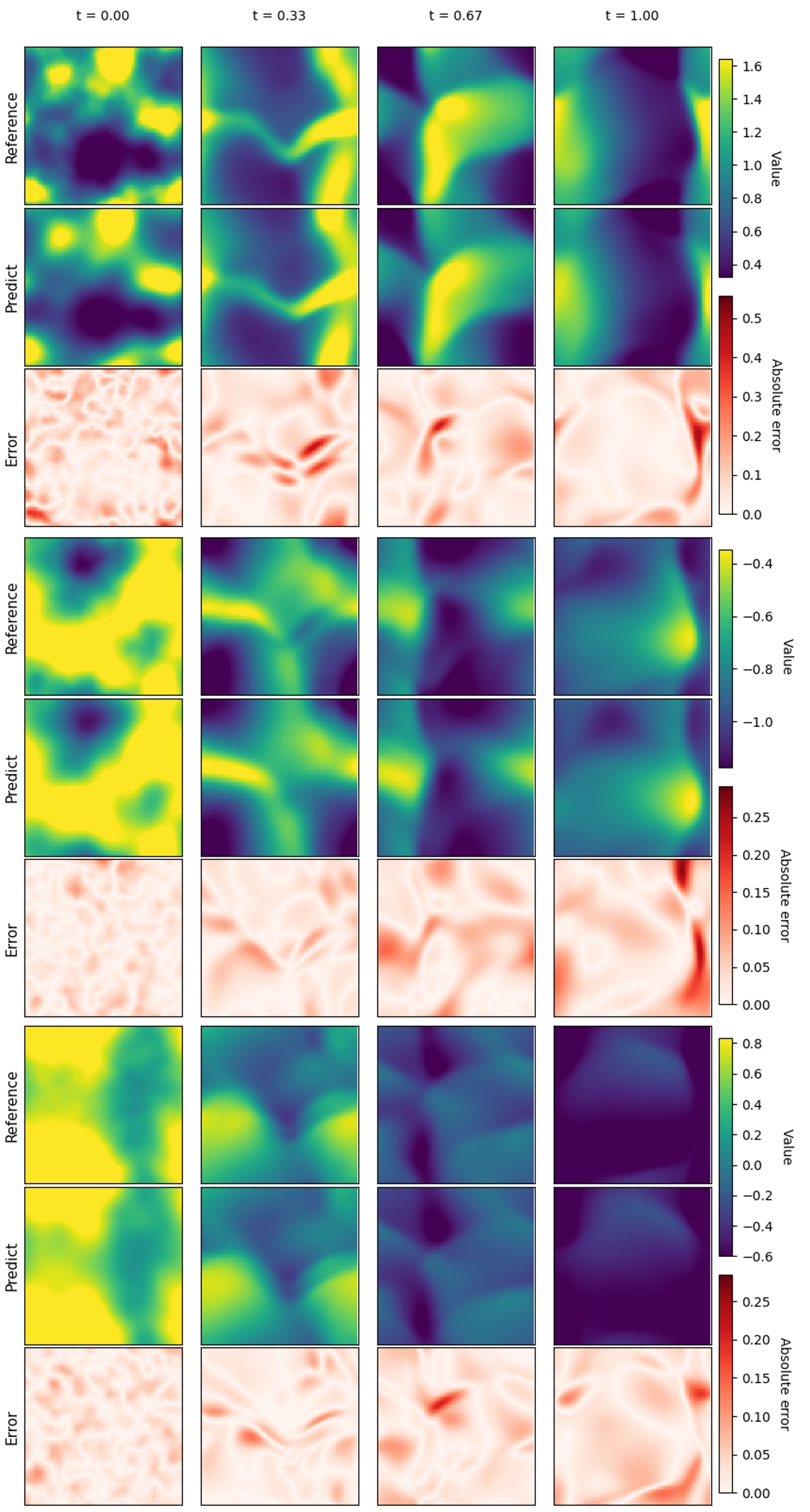}
	\begin{equation*}
		\begin{cases}
			h_t+((h+H)u)_x+((h+H)v)_y-a^h\Delta h+s^h(r)+c_{02}v+u=0, \\
			u_t+uu_x+vu_y+g_1h_x-a^u\Delta u+s^u(r)+v+c_{10}h=0, \\
			v_t+uv_x+vv_y+g_2h_y-a^v\Delta v+1+c_{20}h+c_{22}v+b_{201}hu=0
		\end{cases}
	\end{equation*}
	\caption{Visualization of the generalized shallow water equation on a square domain periodic along both axes, including all solution components $h,u,v$.}
	\label{fig:swe}
\end{figure}

We consider the generalized shallow water equation on a square domain periodic along both axes.
The visualization of the solution components $h,u,v$ is shown in Fig.~\ref{fig:swe}.

\subsection{Sine-Gordon equation}
\begin{figure}
	\centering
	\includegraphics[width=\textwidth]{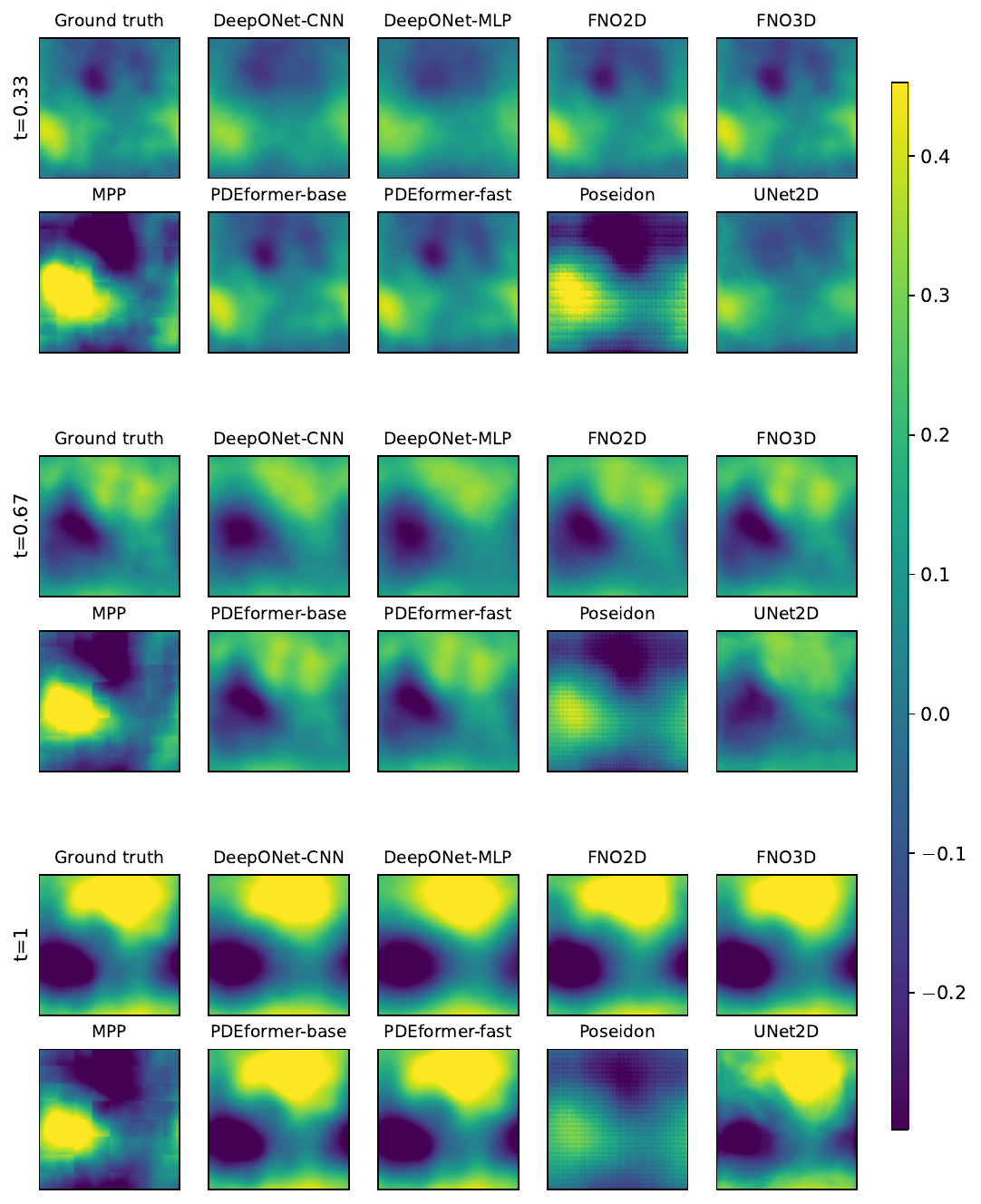}
	\begin{equation*}
		u_{tt}-a\Delta u+(-1)\sin(u)=0
	\end{equation*}
	\caption{Visualization of the Sine-Gordon equation on a square domain periodic along both axes.}
	\label{fig:sine-gordon}
\end{figure}

We consider the sine-Gordon equation on a square domain periodic along both axes, and train the models using 900 data samples.
Visualization of the solution is shown in Fig.~\ref{fig:sine-gordon}.
\subsection{INS-Tracer}
\begin{figure}
	\centering
	\includegraphics[width=\textwidth]{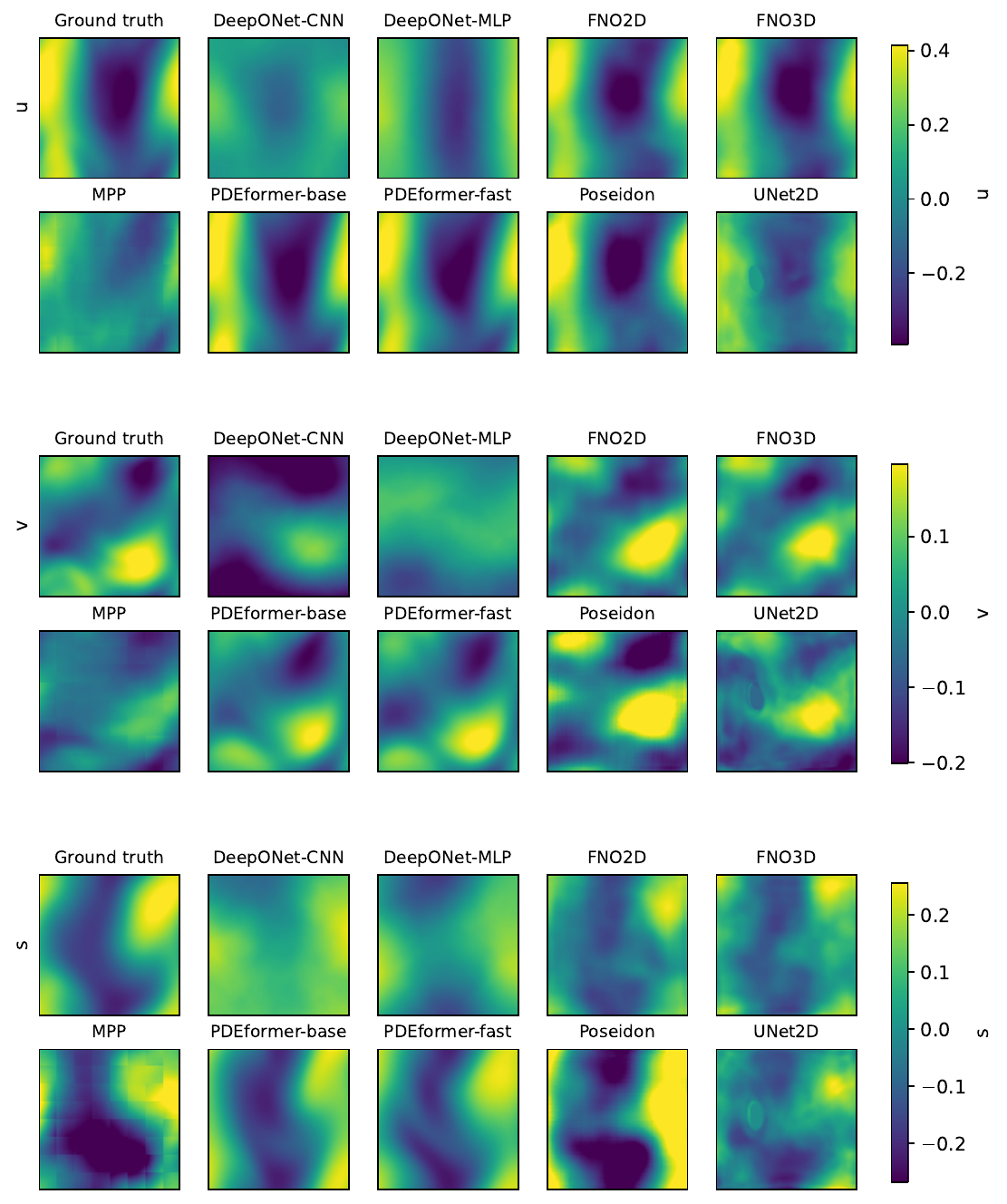}
	\begin{equation*}
		\begin{cases}
			u_t-\nu\Delta u+(u^2)_x+(uv)_y+p_x&=0,\\
			v_t-\nu\Delta v+(uv)_x+(v^2)_y+p_y&=0,\\
			s_t-D\Delta s+(us)_x+(vs)_y&=0,\\
			u_x+v_y&=0
		\end{cases}
	\end{equation*}
	\caption{Visualization of the incompressible Navier-Stokes equation with tracer particle (INS-Tracer) on a square domain periodic along both axes, showing the solution components $u,v,s$ at $t=1$.}
	\label{fig:ins-tracer}
\end{figure}

We consider the Incompressible Navier-Stokes Equation with Tracer on a square domain periodic along both axes, and train the models using 80 data samples.
Visualization of the variables $u,v,p$ on the last timestep $(t=1)$ is shown in Fig.~\ref{fig:ins-tracer}.
\subsection{INS-Pipe}
\begin{figure}
	\centering
	\includegraphics[width=\textwidth]{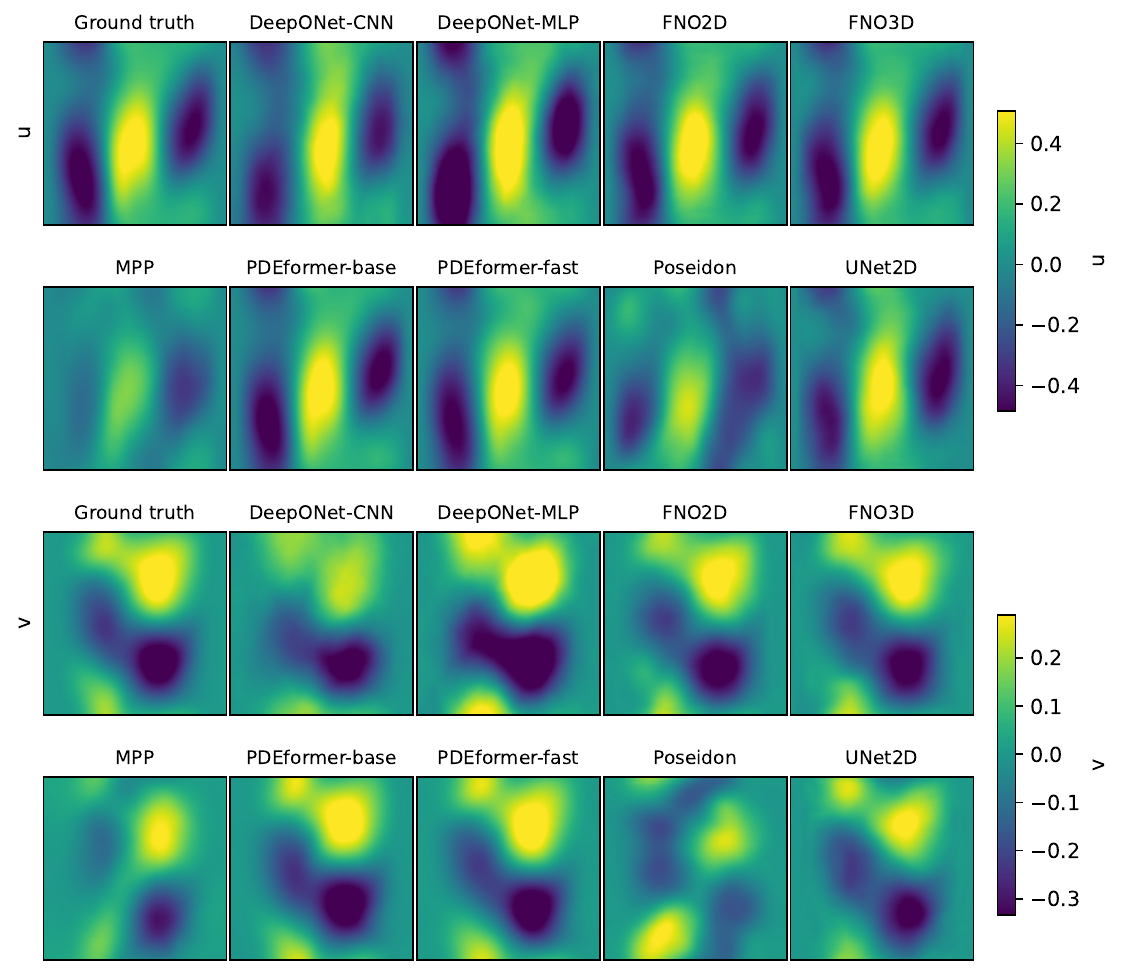}
	\begin{equation*}
		\begin{cases}
			u_t-\nu\Delta u+s_0(\boldsymbol{r})+(u^2)_x+(uv)_y+p_x&=0,\\
			v_t-\nu\Delta v+s_1(\boldsymbol{r})+(uv)_x+(v^2)_y+p_y&=0,\\
			u_x+v_y&=0,
		\end{cases}
	\end{equation*}
	\caption{Visualization of the incompressible Navier-Stokes equation in a pipe (INS-Pipe), showing the solution components $u,v$ at $t=1$.}
	\label{fig:ins-pipe}
\end{figure}

We consider the Incompressible Navier-Stokes Equation in a pipe on a square domain periodic along $x$-axis, and train the models using 900 data samples.
Visualization of the variables $u,v$ at the last timestep $(t=1)$ is shown in Fig.~\ref{fig:ins-pipe}.
\subsection{Wave-Gauss}
\begin{figure}
	\centering
	\includegraphics[width=\textwidth]{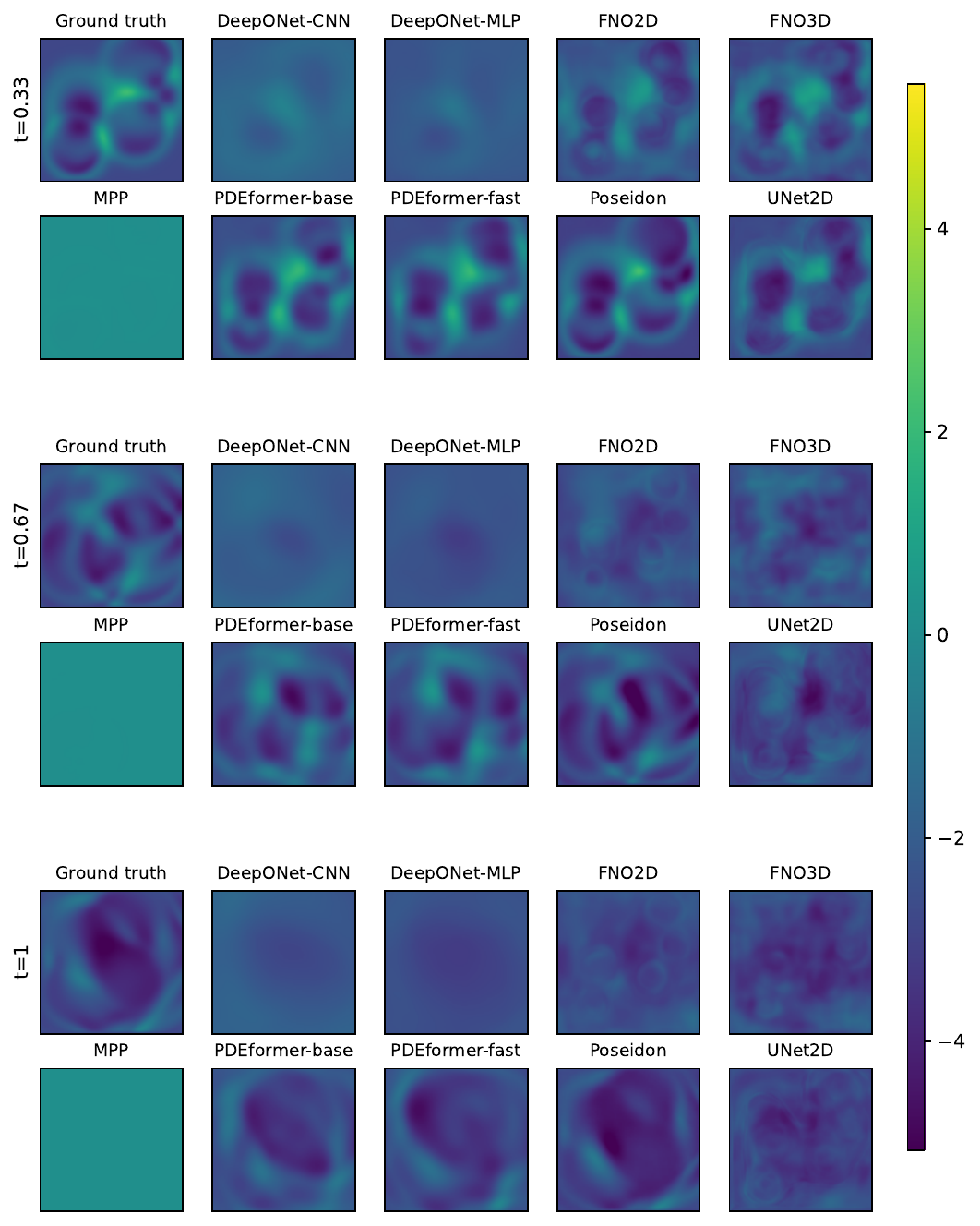}
	\begin{equation*}
		u_{tt}-c(\boldsymbol{r})^2 \Delta u=0
	\end{equation*}
	\caption{Visualization of the wave equation with Gaussian initial conditions.}
	\label{fig:wave-gauss}
\end{figure}

We consider the wave equation with Gaussian initial conditions on a square domain with absorbing boundaries, and train the models using 80 data samples.
Visualization of the solution is shown in Fig.~\ref{fig:wave-gauss}.

\end{document}